\documentclass[3p]{elsarticle}
\usepackage[T1]{fontenc}
\usepackage[latin9]{inputenc}
\usepackage{array}
\usepackage{float}
\usepackage{booktabs}
\usepackage{units}
\usepackage{multirow}
\usepackage{amsmath}
\usepackage{amsthm}
\usepackage{amssymb}
\usepackage{graphicx}
\usepackage{esint}
\usepackage{lineno}
\usepackage{hyperref}
\usepackage{epstopdf}
\makeatletter

\providecommand{\tabularnewline}{\\}

  \theoremstyle{definition}
  \newtheorem{defn}{\protect\definitionname}
 \theoremstyle{definition}
  \newtheorem{example}{\protect\examplename}

\@ifundefined{showcaptionsetup}{}{%
 \PassOptionsToPackage{caption=false}{subfig}}
\usepackage{subfig}
\AtBeginDocument{
  
}

\makeatletter
\def\ps@pprintTitle{%
   \let\@oddhead\@empty
   \let\@evenhead\@empty
   \def\@oddfoot{\reset@font\hfil\thepage\hfil}
   \let\@evenfoot\@oddfoot
}

\makeatother

\providecommand{\definitionname}{Definition}
\providecommand{\examplename}{Example}
\renewcommand{\vec}[1]{\mathbf{#1}}

\begin{document}

\begin{frontmatter}

\title{Automated Computation of Autonomous \\ Spectral Submanifolds for Nonlinear Modal Analysis}

\author[mymainaddress]{Sten Ponsioen}
\ead{stenp@ethz.ch}
\author[mymainaddress]{Tiemo Pedergnana}
\ead{ptiemo@student.ethz.ch}
\author[mymainaddress]{George Haller\corref{mycorrespondingauthor}}
\ead{georgehaller@ethz.ch}

\address[mymainaddress]{Institute for Mechanical Systems \\ ETH Zürich, Leonhardstrasse 21, 8092 Zürich, Switzerland}
\cortext[mycorrespondingauthor]{Corresponding author.}

\begin{abstract}
We discuss an automated computational methodology for computing two-dimensional
spectral submanifolds (SSMs) in autonomous nonlinear mechanical systems
of arbitrary degrees of freedom. In our algorithm, SSMs, the smoothest
nonlinear continuations of modal subspaces of the linearized system,
are constructed up to arbitrary orders of accuracy, using the parameterization
method. An advantage of this approach is that the construction of
the SSMs does not break down when the SSM folds over its underlying
spectral subspace. A further advantage is an automated a posteriori error estimation
feature that enables a systematic increase in the orders of the SSM
computation until the required accuracy is reached. We find that the
present algorithm provides a major speed-up, relative to numerical
continuation methods, in the computation of backbone curves, especially
in higher-dimensional problems. We illustrate the accuracy and speed
of the automated SSM algorithm on lower- and higher-dimensional mechanical
systems.
\end{abstract}


\end{frontmatter}

\section{Introduction}

A fundamental notion in decomposing nonlinear mechanical oscillations,
is the \textit{nonlinear normal mode} (NNM) concept of Rosenberg \cite{Rosenberg1962},
who defined a nonlinear normal mode as a synchronous periodic oscillation
that reaches its maximum in all modal coordinates at the same time.
An alternative definition of a NNM, proposed by Shaw and Pierre \cite{Shaw1993},
is an invariant manifold that serves as the nonlinear continuation
of two-dimensional subspaces formed by normal modes of the linearized
system. Shaw and Pierre seek such invariant manifolds as graphs over
those two-dimensional subspaces. For several extensive discussions
about these two NNM definitions, we refer the reader to the work of
Kerschen et al. \cite{Kerschen2009}, Peeters et al. \cite{Peeters2009},
Mikhlin and Avramov \cite{mikhlin2010nonlinears} and Vakakis et al.
\cite{Vakakis2001}. 

If one relaxes the synchronicity requirement of Rosenberg, a clear
relationship between the above two views on NNMs emerges for conservative
oscillatory systems by the Lyapunov subcenter-manifold theorem \cite{Lyapunov1992,Kelley1969}.
Indeed, under appropriate non-resonance conditions, these references
guarantee the existence of a unique, analytic and robust Shaw\textemdash Pierre-type
invariant manifold tangent to each two-dimensional modal subspace
of the linearized system. This manifold, in turn, is filled with Rosenberg-type
periodic orbits. 

In a non-conservative setting, this geometrical relationship between
the two classic NNM concepts no longer holds, as periodic orbits become
rare and isolated in the phase space, whereas infinitely many invariant
manifolds tangent to each two-dimensional modal subspace will exist.
A unified approach has been proposed by Haller and Ponsioen \cite{Haller2016}
to clarify the relationship between the Rosenberg and Shaw\textemdash Pierre
NNM concepts. Specifically, \cite{Haller2016} defines a nonlinear
normal mode simply as a recurrent motion with finitely many frequencies.
Included in this theory is a trivial NNM or fixed point (no frequencies),
a periodic NNM (the frequencies are rationally commensurate, as for
a Rosenberg-type periodic orbit) and a quasiperiodic NNM (the frequencies
are rationally incommensurate, with the orbit filling an invariant
torus). 

Using this NNM definition, Haller and Ponsioen \cite{Haller2016}
define a spectral submanifold (SSM) as the smoothest invariant manifold
tangent to a modal subspace of a NNM. They then invoke rigorous existence,
uniqueness and persistence results for autonomous and non-autonomous
SSMs, providing an exact mathematical foundation for constructing
nonlinear reduced-order models over appropriately chosen spectral
subspaces. These models are obtained by reducing the full dynamics
to the exactly invariant SSM surfaces, tangent to those subspaces.

More recently, Szalai et al. \cite{Szalai2017} have shown that the
dynamics on a single-mode SSM can be seen as a nonlinear extension
of the linear dynamics of the underlying modal subspace, making it
possible to extract the \textit{backbone curve}, defined as a graph
plotting the instantaneous amplitude of vibration as a function of
the instantaneous frequency of vibration. This approach to backbone-curve
computations assumes that the Lyapunov subcenter-manifold perturbs
smoothly to a unique SSM under appropriate non-resonance conditions
and under small enough damping, which is consistent with the numerical
observations as shown by Kerschen et al. \cite{Kerschen2006}, Peeters
et al. \cite{Peeters2011} and Szalai et al. \cite{Szalai2017}. 

Computing invariant manifolds tangent to modal subspaces in realistic
applications has been a challenge. Prior approaches have mostly focussed
on solving the invariance equations that such manifolds have to satisfy
(Blanc et al. \cite{Blanc2013}, Pesheck et al. \cite{Pesheck2002}
and Renson et al. \cite{Renson2014}). These invariance equations
have infinitely many solutions, out of which the numerical approaches
employed by different authors selected one particular solution. In
contrast, the SSM theory employed here guarantees a unique solution
that can be approximated with arbitrary high precision via the parameterization
method of Cabré et al. \cite{Cabre2003,Cabre2003b,Cabre2005}. In
the present work, we describe an automated computational algorithm
for two-dimensional SSMs constructed over two-dimensional modal subspaces.
This algorithm\footnote{SSMtool is available at: \href{http://www.georgehaller.com}{www.georgehaller.com}}
allows us to compute the SSMs, their reduced dynamics and associated
backbone curves to arbitrary orders of accuracy, limited only by available
memory. An important feature of the algorithm is a direct a posteriori estimation
of the error in computing the SSM at a given approximation order.
This error estimate measures directly the extend to which the SSM
is invariant. If the error is unsatisfactory, the user can select
higher order approximations until the error falls below a required
bound.

In technical terms, we construct the SSMs as embeddings of the modal
subspaces into the phase space of the mechanical system, as required
by the parameterization method (Cabré et al. \cite{Cabre2003,Cabre2003b,Cabre2005}).
A major advantage compared to most earlier calculations (Haller and
Ponsioen \cite{Haller2016}) is that the parameterized construction
of SSMs does not break down when the SSM folds over the underlying
modal subspace. Another advantage of the method is its suitability
for algorithmic implementations for arbitrary orders of accuracy in
arbitrary dimensions. For applications of the parameterization method
to other types of dynamical systems, we refer the reader to the work
of Haro et al. \cite{Haro2016}, van den Berg and Mireles James \cite{Berg2016}
and Mireles James \cite{James2015}.

\section{System set-up }

We consider $n$-degree-of-freedom, autonomous mechanical systems
of the form

\begin{equation}
\vec{M}\ddot{\vec{y}}+\vec{C}\dot{\vec{y}}+\vec{K}\vec{y}+\vec{f}(\vec{y},\dot{\vec{y}})=\vec{0},\qquad \vec{f}(\vec{y},\dot{\vec{y}})=\mathcal{O}\left(\left|\vec{y}\right|^{2},\left|\vec{y}\right|\left|\dot{\vec{y}}\right|,\left|\dot{\vec{y}}\right|^{2}\right),\label{eq:mech_sys}
\end{equation}
where $\vec{y}\in\mathbb{R}^{n}$ is the generalized position vector; $\vec{M}=\vec{M}^{T}\in\mathbb{R}^{n\times n}$
is the positive definite mass matrix; $\vec{C}=\vec{C}^{T}\in\mathbb{R}^{n\times n}$
is the damping matrix; $\vec{K}=\vec{K}^{T}\in\mathbb{R}^{n\times n}$ is the
stiffness matrix and $\vec{f}(\vec{x},\dot{\vec{x}})$ denotes all the nonlinear terms
in the system. These nonlinearities are assumed to be of class $C^{r}$
in $(\vec{x},\dot{\vec{x}})$, with $r\in\mathbb{N}^{+}\cup\left\{ \infty,a\right\} $.
Here $r\in\mathbb{N}^{+}$ refers to finite differentiability, $r=\infty$
refers to infinite differentiability, and $r=a$ refer to analyticity,
all three of which are allowed in our treatment. 

System (\ref{eq:mech_sys}) can be transformed into a set of $2n$
first-order ordinary differential equations by introducing a change
of variables $\vec{x}_{1}=\vec{y}$, $\vec{x}_{2}=\dot{\vec{y}}$, with $\vec{x}=(\vec{x}_{1},\vec{x}_{2})\in\mathbb{R}^{2n}$,
which gives,
\begin{gather}
\dot{\vec{x}}=\left(\begin{array}{cc}
\vec{0} & \vec{I}\\
-\vec{M}^{-1}\vec{K} & -\vec{M}^{-1}\vec{C}
\end{array}\right)\vec{x}+\left(\begin{array}{c}
\vec{0}\\
-\vec{M}^{-1}\vec{f}(\vec{x}_{1},\vec{x}_{2})
\end{array}\right)=\vec{A}\vec{x}+\vec{F}(\vec{x}),\label{eq:dyn_sys}\\
\vec{x}\in\mathbb{R}^{2n},\qquad \vec{F}(\vec{x})=\mathcal{O}\left(\left|\vec{x}\right|^{2}\right).\nonumber 
\end{gather}
The transformed system (\ref{eq:dyn_sys}) has a fixed point at $\vec{x}=\vec{0}$,
$\vec{A}\in\mathbb{R}^{2n\times2n}$ is a constant matrix and $\vec{F}(\vec{x})$ is
a class $C^{r}$ function containing all the nonlinearities. Note
that the inverse of the mass matrix is well-defined because $\vec{M}$ is
assumed positive definite. 

The linearized part of (\ref{eq:dyn_sys}) is 

\begin{equation}
\dot{\vec{x}}=\vec{A}\vec{x},\label{eq:lin_dyn_sys}
\end{equation}
where the matrix $\vec{A}$ has $2n$ eigenvalues $\lambda_{k}\in\mathbb{C}$
for $k=1,\ldots,2n$. Counting multiplicities, we sort these eigenvalues
based on their real parts in the decreasing order,

\begin{equation}
\text{\text{Re}}(\lambda_{2n})\leq\text{\text{Re}}(\lambda_{2n-1})\leq\ldots\leq\text{\text{Re}}(\lambda_{1})<0,
\end{equation}
assuming that the real part of each eigenvalue is less than zero and
hence the fixed point is asymptotically stable. We further assume
that the constant matrix $\vec{A}$ is semisimple, which implies that the
algebraic multiplicity of each $\lambda_{k}$ is equal to its geometric
multiplicity, i.e.\ $\text{alg}(\lambda_{k})=\text{geo}(\lambda_{k})$.
We can therefore identify $2n$ linearly independent eigenvectors
$\vec{v}_{k}\in\mathbb{C}^{2n}$, with $k=1,\ldots,2n$, each spanning a
real eigenspace $E_{k}\subset\mathbb{R}^{2n}$ with $\text{dim}(E_{k})=2\times\text{alg}(\lambda_{k})$
in case $\text{Im}(\lambda_{k})\neq0$, or $\text{dim}(E_{k})=\text{alg}(\lambda_{k})$
in case $\text{Im}(\lambda_{k})=0$.

\section{Spectral submanifolds for continuous mechanical systems }

As $\vec{A}$ is semisimple, the linear part of system (\ref{eq:dyn_sys})
can be diagonalized by introducing a linear change of coordinates
$\vec{x}=\vec{T}\vec{q}$, with $\vec{T}=\left[\vec{v}_{j_{1}},\vec{v}_{j_{2}},\ldots,\vec{v}_{j_{2n}}\right]\in\mathbb{C}^{2n\times2n}$
and $\vec{q}\in\mathbb{C}^{2n}$,

\begin{gather}
\dot{\vec{q}}=\vec{T}^{-1}\vec{A}\vec{T}\vec{q}+\vec{T}^{-1}\vec{F}(\vec{T}\vec{q})=\underbrace{\text{diag}(\lambda_{j_{1}},\lambda_{j_{2}}\ldots,\lambda_{j_{2n}})}_{\vec{\Lambda}}\vec{q}+\underbrace{\vec{T}^{-1}\vec{F}(\vec{T}\vec{q})}_{\vec{G}(\vec{q})}=\vec{\Lambda}\vec{q}+\vec{G}(\vec{q}).\label{eq:ds_diag}
\end{gather}
We now seek a two-dimensional modal subspace $\mathcal{E}=\text{span}\left\{ \vec{v}_{j_{1}},\vec{v}_{j_{2}}\right\} \subset\mathbb{C}^{2n}$,
with $\vec{v}_{j_{2}}=\bar{\vec{v}}_{j_{1}}.$ Note that $\vec{v}_{j_{1}}$ and $\vec{v}_{j_{2}}$
are purely real if $\lambda_{j_{1}},\lambda_{j_{2}}\in\mathbb{R}$,
in which case $\mathcal{E}$ corresponds to either a single critically
damped mode ($\lambda_{j_{1}}=\lambda_{j_{2}}$), or to two overdamped
modes ($\lambda_{j_{1}}\neq\lambda_{j_{2}}$). In contrast, if $\lambda_{j_{1}},\bar{\lambda}_{j_{2}}\in\mathbb{C}$
with $\text{Im }\lambda_{j_{1}}\neq0$, then $\mathcal{E}$ corresponds
to a single underdamped mode. 

The remaining linearly independent eigenvectors $\vec{v}_{j_{3}},\ldots,\vec{v}_{j_{2n}}$
span a complex subspace $\mathcal{C}\subset\mathbb{C}^{2n}$ such
that the full phase space of (\ref{eq:ds_diag}) can be expressed
as the direct sum

\begin{equation}
\mathbb{C}^{2n}=\mathcal{E}\oplus\mathcal{C}.
\end{equation}
The diagonal matrix $\vec{\Lambda}$ is the representation of the linear
matrix $\vec{A}$ with respect to this decomposition, which we can also
write as

\begin{equation}
\vec{\Lambda}=\left[\begin{array}{cc}
\vec{\Lambda}_{\mathcal{E}} & 0\\
0 & \vec{\Lambda}_{\mathcal{C}}
\end{array}\right],\quad\text{Spect}\left(\vec{\Lambda}_{\mathcal{E}}\right)=\left\{ \lambda_{j_{1}},\lambda_{j_{2}}\right\} ,\quad\text{Spect}\left(\vec{\Lambda}_{\mathcal{C}}\right)=\left\{ \lambda_{j_{3}},\ldots,\lambda_{j_{2n}}\right\} ,\label{eq:lin_decomp}
\end{equation}
with $\vec{\Lambda}_{\mathcal{E}}=\text{diag}(\lambda_{j_{1}},\lambda_{j_{2}})$
and $\vec{\Lambda}_{\mathcal{C}}=\text{diag}(\lambda_{j_{3}},\ldots,\lambda_{j_{2n}})$.

Following the work of Haller and Ponsioen \cite{Haller2016}, we now
define a spectral submanifold of the nonlinear system (\ref{eq:ds_diag})
as an invariant manifold tangent to the spectral subspace $\mathcal{E}$.
\begin{defn}
A \textit{spectral submanifold} (SSM), $\mathcal{W}(\mathcal{E})$,
corresponding to a spectral subspace $\mathcal{E}$ of $\vec{\Lambda}$
is an invariant manifold of the dynamical system (\ref{eq:ds_diag})
such that
\begin{description}
\item [{(i)}] $\mathcal{W}(\mathcal{E})$ is tangent to $\mathcal{E}$
and has the same dimension as $\mathcal{E}$.
\item [{(ii)}] $\mathcal{W}(\mathcal{E})$ is strictly smoother than any
other invariant manifold satisfying (i).
\end{description}
\end{defn}
We define the \textit{outer spectral quotient} $\sigma_{\text{out}}(\mathcal{E})$
as the positive integer 

\begin{equation}
\sigma_{\text{out}}(\mathcal{E})=\text{Int}\left[\frac{\min_{\lambda\in\text{Spect}(\vec{\Lambda}_{\mathcal{C}})}\text{Re}\lambda}{\max_{\lambda\in\text{Spect}(\vec{\Lambda}_{\mathcal{E}})}\text{Re}\lambda}\right]\in\mathbb{N}^{+}.\label{eq:outer_spect_quo}
\end{equation}
Which is the integer part of the ratio between the strongest decay
rate of the linearized oscillations outside $\mathcal{E}$ and the
weakest decay rate of the linearized oscillations inside $\mathcal{E}$.
As we will see shortly, $\sigma_{\text{out}}(\mathcal{E})$ determines
the smoothness class in which $\mathcal{W}(\mathcal{E})$ turns out
to be unique. 

To state the main results on SSMs from Haller and Ponsioen \cite{Haller2016},
we need the following two assumptions:
\begin{description}
\item [{(A1)}] $\sigma_{\text{out}}(\mathcal{E})\leq r$,
\item [{(A2)}] The \textit{outer} non-resonance conditions 

\begin{equation}
a\lambda_{j_{1}}+b\lambda_{j_{2}}\neq\lambda_{l},\quad\forall\lambda_{l}\in\text{Spect}(\vec{\Lambda}_{\mathcal{C}}).\label{eq:ext_res}
\end{equation}

\end{description}
hold for all positive integers $a$ and $b$ satisfying $2\leq a+b\leq\sigma_{\text{out}}(\mathcal{E})$.

Under these assumptions, we have the following main result on an SSM
tangent to the modal subspace $\mathcal{E}$ in system (\ref{eq:ds_diag}).
\begin{description}
\item [{Theorem}] \textbf{1} 
\item [{(i)}] There exist a two-dimensional SSM, $\mathcal{W}(\mathcal{E})$,
that is tangent to the spectral subspace $\mathcal{E}$ at the fixed
point $\vec{q}=\vec{0}$.
\item [{(ii)}] $\mathcal{W}(\mathcal{E})$ is of class $C^{r}$ and is
unique among all class $C^{\sigma_{\text{out}}(\mathcal{E})+1}$ two-dimensional
invariant manifolds that are tangent to $\mathcal{E}$ at $\vec{q}=\vec{0}$.
\item [{(iii)}] $\mathcal{W}(\mathcal{E})$ can be parameterized over an
open set $\mathcal{U}\subset\mathbb{C}^{2}$ via the map

\begin{equation}
\vec{W}:\mathcal{U}\subset\mathbb{C}^{2}\rightarrow\mathbb{C}^{2n},\label{eq:W_map}
\end{equation}
into the phase space of system (\ref{eq:ds_diag}).
\item [{(iv)}] There exist a $C^{r}$ polynomial function $\vec{R}:\mathcal{U}\rightarrow\mathcal{U}$
satisfying the following invariance relationship 

\begin{equation}
\vec{\Lambda} \vec{W}+\vec{G}\circ \vec{W}=\nabla \vec{W}\vec{R},\label{eq:invar}
\end{equation}
such that the reduced dynamics on the SSM can be expressed as

\begin{equation}
\dot{\vec{z}}=\vec{R}(\vec{z}),\quad \vec{R}(\vec{0})=\vec{0},\quad\nabla \vec{R}(\vec{0})=\vec{\Lambda}_{\mathcal{E}}=\left[\begin{array}{cc}
\lambda_{j_{1}} & 0\\
0 & \lambda_{j_{2}}
\end{array}\right],\quad \vec{z}=(z_{j_{1}},z_{j_{2}})\in\mathcal{U}.\label{eq:map_R}
\end{equation}

\item [{(v)}] If the \textit{inner} non-resonance conditions 

\begin{gather}
a\lambda_{j_{1}}+b\lambda_{j_{2}}\neq\lambda_{j_{i}},\quad i=1,2\label{eq:int_res}
\end{gather}

hold for all positive integers $a$ and $b$ with

\[
2\leq a+b\leq\sigma_{\text{in}}(\mathcal{E})=\text{Int}\left[\frac{\min_{\lambda\in\text{Spect}(\vec{\Lambda}_{\mathcal{E}})}\text{Re}\lambda}{\max_{\lambda\in\text{Spect}(\vec{\Lambda}_{\mathcal{E}})}\text{Re}\lambda}\right]\in\mathbb{N}^{+},
\]
then the mapping $\vec{W}$ can be chosen in such a way that $\vec{R}(\vec{z})$ does
not contain the terms $z_{j_{1}}^{a}z_{j_{2}}^{b}$. In particular,
if no inner resonances arise, then the reduced dynamics on the SSM
can be linearized. 
\end{description}
The proof of Theorem 1 can be found in the work of Haller and Ponsioen
\cite{Haller2016}, which is based on the more abstract results of
Cabré et al. \cite{Cabre2003,Cabre2003b,Cabre2005} for mappings on
Banach spaces. 

\section{SSM computation \label{sec:ssmcomp}}

In this section we show how the parameterized spectral submanifolds
are approximated, around a fixed point, using polynomials. We express
$\vec{W}(\vec{z})$, $\vec{R}(\vec{z})$ and $\vec{G}(\vec{q})$ as multivariate polynomial functions,
which is done by using the Kronecker product. Substituting the expressions
in the invariance equation (\ref{eq:invar}), we obtain, for each
order, a linear system of equations that can be solved under appropriate
non-resonance conditions. For a different application of the parameterization
method, we refer to the work of Mireles James \cite{James2015}. Where
the parameterization method is used for approximating (un)stable manifolds
of one parameter families of analytic dynamical systems, by using
polynomials. For a more elaborate discussion about the numerical computation
of the coefficients of higher order power series expansions of parameterized
invariant manifolds around a fixed point of an elementary vector field,
where the coefficients of the power series expansions are obtained
using methods of Automatic Differentiation, we refer to Haro et al.
\cite{Haro2016}.

\subsection{The Kronecker product}

We now describe a computational algorithm for constructing the mapping
$\vec{W}(\vec{z})$ in (\ref{eq:W_map}) that maps $\mathcal{U}\subset\mathbb{C}^{2}$
into the phase space of system (\ref{eq:dyn_sys}), and the reduced
dynamics $\vec{R}(\vec{z})$ in (\ref{eq:map_R}). To handle the polynomial calculations
arising in the algorithm efficiently, we first need to recall the
notion and some properties of the Kronecker product \cite{Laub2005}. 
\begin{defn}
\label{def:kronecker}Let $\vec{A}\in\mathbb{C}^{m\times n}$, $\vec{B}\in\mathbb{C}^{p\times q}$.
Then we define the Kronecker product of $\vec{A}$ and $\vec{B}$ as

\begin{equation}
\vec{A}\otimes \vec{B}=\left[\begin{array}{ccc}
a_{11}\vec{B} & \cdots & a_{1n}\vec{B}\\
\vdots & \ddots & \vdots\\
a_{m1}\vec{B} & \cdots & a_{mn}\vec{B}
\end{array}\right]\in\mathbb{C}^{mp\times nq}.
\end{equation}
\end{defn}
In accordance with Definition \ref{def:kronecker}, the Kronecker
product of two vectors $\vec{x}\in\mathbb{C}^{m}$ and $\vec{y}\in\mathbb{C}^{n}$
is given by the vector
\begin{align}
\vec{x}\otimes \vec{y} & =\left[\begin{array}{c}
x_{1}\vec{y}\\
\vdots\\
x_{m}\vec{y}
\end{array}\right]=\left[x_{1}y_{1}\cdots x_{1}y_{n}\cdots x_{m}y_{1}\cdots x_{m}y_{n}\right]^{T}\in\mathbb{C}^{mn},
\end{align}
or equivalently, written in index notation,

\begin{equation}
\vec{x}\otimes \vec{y}=\sum_{i=1}^{m}\sum_{j=1}^{n}x_{i}y_{j}\vec{e}_{i}^{x}\otimes\vec{e}_{j}^{y},\label{eq:index_kron}
\end{equation}
where $\vec{e}_{i}^{x}\in\mathbb{C}^{m}$ and $\vec{e}_{j}^{y}\in\mathbb{C}^{n}$
are basis vectors containing a one in their $i^{\text{th}}$ and $j^{\text{th}}$
entries, respectively, and zeros elsewhere. Differentiating equation
(\ref{eq:index_kron}) with respect to time yields

\begin{equation}
\frac{d}{dt}\left(\vec{x}\otimes \vec{y}\right)=\frac{d}{dt}\left(x_{i}y_{j}\vec{e}_{i}^{x}\otimes\vec{e}_{j}^{y}\right)=\dot{x}_{i}y_{j}\vec{e}_{i}^{x}\otimes\vec{e}_{j}^{y}+x_{i}\dot{y}_{j}\vec{e}_{i}^{x}\otimes\vec{e}_{j}^{y}=\dot{\vec{x}}\otimes \vec{y}+\vec{x}\otimes\dot{\vec{y}}.\label{eq:kron_derivative}
\end{equation}
which is simply the product rule. By using the same reasoning, one
shows that the product rule also applies to the time derivative of
the Kronecker product of multiple vectors. We will use the shorthand
notation $\vec{z}^{\otimes i}$ defined as

\begin{equation}
\vec{z}^{\otimes i}=\begin{cases}
\vec{z} & \quad\text{for }i=1.\\
\underbrace{\vec{z}\otimes \vec{z}\otimes\cdots\otimes \vec{z}}_{i\text{ times}} & \quad\text{for }i>1.
\end{cases}
\end{equation}
For subsequent derivations, we will make use of several properties
of the Kronecker product, which we list in \ref{sec:Appendix-A}
for convenience. 

We Taylor expand $\vec{W}(\vec{z})$ and $\vec{R}(\vec{z})$ and express them as multivariate
polynomial functions

\begin{align}
\vec{W}(\vec{z}) & =\sum_{i=1}^{n_{w}}\vec{W}_{i}\vec{z}^{\otimes i}=\vec{W}_{1}\vec{z}+\vec{W}_{2}\vec{z}\otimes \vec{z}+\vec{W}_{3}\vec{z}\otimes \vec{z}\otimes \vec{z}+\ldots,\quad \vec{W}_{i}\in\mathbb{C}^{2n\times2^{i}},\quad \vec{z}\in\mathbb{C}^{2},\label{eq:Wpoly}\\
\vec{R}(\vec{z}) & =\sum_{i=1}^{n_{w}}\vec{R}_{i}\vec{z}^{\otimes i}=\vec{R}_{1}\vec{z}+\vec{R}_{2}\vec{z}\otimes \vec{z}+\vec{R}_{3}\vec{z}\otimes \vec{z}\otimes \vec{z}+\ldots,\quad \vec{R}_{i}\in\mathbb{C}^{2\times2^{i}},\quad \vec{z}\in\mathbb{C}^{2},\label{eq:Rpoly}
\end{align}
with $n_{w}\geq\sigma_{\text{out}}(\mathcal{E})+1$ denoting the order
of the SSM expansion. We also Taylor expand the nonlinear part of
our dynamical system (\ref{eq:ds_diag}) up to order $n_{w}$, around
the fixed point $\vec{q}=\vec{0}$, such that we can represent the nonlinearities,
in a fashion similar to equations (\ref{eq:Wpoly}) and (\ref{eq:Rpoly}),
as

\begin{equation}
\vec{G}(\vec{q})=\sum_{i=2}^{\Gamma}\vec{G}_{i}\vec{q}^{\otimes i}=\vec{G}_{2}\vec{q}\otimes \vec{q}+\vec{G}_{3}\vec{q}\otimes \vec{q}\otimes \vec{q}+\ldots,\quad \vec{G}_{i}\in\mathbb{C}^{2n\times(2n)^{i}},\quad \vec{q}\in\mathbb{C}^{2n},\label{eq:nonlin_poly}
\end{equation}
with $\Gamma$ denoting the maximum order of nonlinearity considered. 

By construction, the vector $\vec{q}^{\otimes i}$ will have redundant terms
along its elements, and hence $\vec{G}_{i}$ will have infinitely many possible
representations at the $i^{\text{th}}$ order in equation (\ref{eq:ds_diag}).
The redundancy in $\vec{q}^{\otimes i}$ allows us to introduce constraints
between the different coefficients that are related to the same monomial
term. We can always set these constraints such that $\vec{G}(\vec{q})$ will represent
the $2n$-dimensional polynomial vector $\vec{T}^{-1}\vec{F}(\vec{T}\vec{q})$. For more detail,
we refer the reader to \ref{sec:Multiple-representations}.

\subsection{The coefficient equations \label{sec:coeffeq}}

We recall here the diagonalized dynamical system (\ref{eq:ds_diag})

\begin{equation}
\dot{\vec{q}}=\vec{\Lambda} \vec{q}+\vec{G}(\vec{q}).\label{eq:ds_diag_repeat}
\end{equation}
Substituting $\vec{q}=\vec{W}(\vec{z})$ on the right-hand side of equation (\ref{eq:ds_diag_repeat}),
then differentiating $\vec{q}=\vec{W}(\vec{z})$ with respect to time and substituting
the result $\dot{\vec{q}}=\nabla\vec{W}(\vec{z})\dot{\vec{z}}$ on the left-hand side of equation
(\ref{eq:ds_diag_repeat}), we obtain the invariance relation (\ref{eq:invar})
in the form

\begin{align}
 & \vec{W}_{1}\vec{R}(\vec{z})+\vec{W}_{2}\left(\vec{R}(\vec{z})\otimes \vec{z}+\vec{z}\otimes \vec{R}(\vec{z})\right)\nonumber \\
 & +\vec{W}_{3}\left(\vec{R}(\vec{z})\otimes \vec{z}\otimes \vec{z}+\vec{z}\otimes \vec{R}(\vec{z})\otimes \vec{z}+\vec{z}\otimes \vec{z}\otimes \vec{R}(\vec{z})\right)+\ldots\nonumber \\
 & +\vec{W}_{k}\left(\vec{R}(\vec{z})\otimes \vec{z}^{\otimes k-1}+\Sigma_{j=1}^{k-2}\left(\vec{z}^{\otimes j}\otimes \vec{R}(\vec{z})\otimes \vec{z}^{\otimes k-j-1}\right)+\vec{z}^{\otimes k-1}\otimes \vec{R}(\vec{z})\right)\nonumber \\
 & =\vec{\Lambda} \vec{W}(\vec{z})+\vec{G}_{2}\vec{W}(\vec{z})^{\otimes2}+\ldots+\vec{G}_{\Gamma}\vec{W}(\vec{z})^{\otimes\Gamma},\label{eq:full_eq}
\end{align}
for $k=\{2,\ldots,n_{w}\}$. The time derivative of $\vec{q}=\vec{W}(\vec{z})=\sum_{i}\vec{W}_{i}\vec{z}^{\otimes i}$
can be expressed, using the product rule for the Kronecker product
of vectors, as 

\begin{equation}
\dot{\vec{q}}=\vec{W}_{1}\dot{\vec{z}}+\vec{W}_{2}\left(\dot{\vec{z}}\otimes \vec{z}+\vec{z}\otimes\dot{\vec{z}}\right)+\vec{W}_{3}\left(\dot{\vec{z}}\otimes \vec{z}\otimes \vec{z}+\vec{z}\otimes\dot{\vec{z}}\otimes \vec{z}+\vec{z}\otimes \vec{z}\otimes\dot{\vec{z}}\right)+\ldots.\label{eq:W_der}
\end{equation}
Substituting $\dot{\vec{z}}=\vec{R}(\vec{z})$ into (\ref{eq:W_der}), we precisely
obtain the left-hand side of equation (\ref{eq:full_eq}). Rewriting
equation (\ref{eq:full_eq}) and collecting terms of equal power $\vec{z}^{\otimes i}$
for $i=\{1,\ldots n_{w}\}$, we obtain, for $i=1$,

\begin{equation}
\underbrace{\left[\begin{array}{cc}
\vec{\Lambda}_{\mathcal{E}} & \vec{0}\\
\vec{0} & \vec{\Lambda}_{\mathcal{C}}
\end{array}\right]}_{\vec{\Lambda}}\vec{W}_{1}=\vec{W}_{1}\vec{R}_{1}.\label{eq:linear_part}
\end{equation}
From (\ref{eq:map_R}), we then require that $\nabla\vec{R}(\vec{0})=\vec{\Lambda}_{\mathcal{E}}=\vec{R}_{1}$.
Therefore, equation (\ref{eq:linear_part}) will be satisfied if we
set $\vec{W}_{1}\in\mathcal{C}^{2n\times2}$ equal to 

\begin{equation}
\vec{W}_{1}=\left[\begin{array}{cc}
1 & 0\\
0 & 1\\
0 & 0\\
\vdots & \vdots\\
0 & 0
\end{array}\right].\label{eq:W1}
\end{equation}
For $2\leq i\leq n_{w}$ we have 

\begin{align}
\vec{\Lambda} \vec{W}_{i}-\vec{W}_{i}\overbrace{\sum_{\left|\vec{s}\right|=1}\vec{\Lambda}_{\mathcal{E}}^{s_{1}}\otimes\cdots\otimes\vec{\Lambda}_{\mathcal{E}}^{s_{i}}}^{\tilde{\vec{\Lambda}}_{\mathcal{E},i}} & =\vec{W}_{1}\vec{R}_{i}+\sum_{m=2}^{i-1}\vec{W}_{m}\sum_{\left|\vec{p}\right|=1}\vec{R}_{i+1-m}^{p_{1}}\otimes\ldots\otimes \vec{R}_{i+1-m}^{p_{m}}\label{eq:cohomological}\\
 & -\vec{G}_{i}\vec{W}_{1}^{\otimes i}-\sum_{m=2}^{i-1}\vec{G}_{m}\sum_{\left|\vec{r}\right|=i}\vec{W}_{r_{1}}\otimes\ldots\otimes \vec{W}_{r_{m}},\nonumber 
\end{align}
where we make use of multi-index notation for $\vec{s}=\left\{ s_{1},\ldots,s_{i}\right\} \in\mathbb{N}^{i}$,
$\vec{p}=\left\{ p_{1},\ldots,p_{m}\right\} \in\mathbb{N}^{m}$ and $\vec{r}=\left\{ r_{1},\ldots,r_{m}\right\} \in\mathbb{N}^{m}$.
The notation $\vec{\Lambda}_{\mathcal{E}}^{s_{j}}$ is used to indicate
that the matrix $\vec{\Lambda}_{\mathcal{E}}$ is taken to the power $s_{j}\in\mathbb{N}$,
where the zeroth power will simply return the identity matrix of the
same dimension as $\vec{\Lambda}_{\mathcal{E}}$. Note that we also adopt the same notation for $\vec{R}_{i+1-m}\in\mathbb{C}^{2\times2^{i+1-m}}$, where we set $\vec{R}_{i+1-m}^0\triangleq \vec{I} \in \mathbb{R}^{2\times2}$.
 
The right hand side of equation (\ref{eq:cohomological}) consists
of the lower-order terms $\vec{W}_{j}$ for $2\leq j<i-1$, which are known
for the current order $i$. The term $\vec{R}_{i}$ represents the coefficient
matrix corresponding to the $i^{\text{th}}$-order of the polynomial
$\vec{R}(\vec{z})$. This polynomial depends on the preferred style of parameterization and will be chosen to remove near-inner resonances from the SSM expressions, as explained later in section \ref{sec:The-backbone-curve}. The matrices
$\vec{G}_{j}$, for $2\leq j\leq i$, are known by definition because they
represent the nonlinearities of system (\ref{eq:ds_diag}). 

\subsubsection{Partitioning the coefficient equations}

Due to the diagonal structure of $\vec{\Lambda}$, we can partition equation
(\ref{eq:cohomological}) into the two separate matrix equations 

\begin{align}
\vec{\Lambda}_{\mathcal{E}}\vec{W}_{i}^{\mathcal{E}}-\vec{W}_{i}^{\mathcal{E}}\tilde{\vec{\Lambda}}_{\mathcal{E},i} & =\vec{R}_{i}+\vec{B}_{i}^{\mathcal{E}},\label{eq:parti_co_1}\\
\vec{\Lambda}_{\mathcal{C}}\vec{W}_{i}^{\mathcal{C}}-\vec{W}_{i}^{\mathcal{C}}\tilde{\vec{\Lambda}}_{\mathcal{E},i} & =\vec{B}_{i}^{\mathcal{C}},\label{eq:parti_co_2}
\end{align}
where $\vec{W}_{i}$ is partitioned as 

\begin{equation}
\vec{W}_{i}=\left[\begin{array}{c}
\vec{W}_{i}^{\mathcal{E}}\\
\vec{W}_{i}^{\mathcal{C}}
\end{array}\right]\in\mathbb{C}^{2n\times2^{i}},\quad \vec{W}_{i}^{\mathcal{E}}\in\mathbb{C}^{2\times2^{i}},\quad \vec{W}_{i}^{\mathcal{C}}\in\mathbb{C}^{(2n-2)\times2^{i}}.
\end{equation}
The matrices $\vec{B}_{i}^{\mathcal{E}}\in\mathcal{C}^{2\times2^{i}}$ and
$\vec{B}_{i}^{C}\in\mathcal{C}^{(2n-2)\times2^{i}}$ are such that

\begin{equation}
\left[\begin{array}{c}
\vec{B}_{i}^{\mathcal{E}}\\
\vec{B}_{i}^{\mathcal{C}}
\end{array}\right]=\sum_{m=2}^{i-1}\vec{W}_{m}\sum_{\left|\vec{p}\right|=1}\vec{R}_{i+1-m}^{p_{1}}\otimes\ldots\otimes \vec{R}_{i+1-m}^{p_{m}}-\vec{G}_{i}\vec{W}_{1}^{\otimes i}-\sum_{m=2}^{i-1}\vec{G}_{m}\sum_{\left|\vec{r}\right|=i}\vec{W}_{r_{1}}\otimes\ldots\otimes \vec{W}_{r_{m}}.\label{eq:Bmat}
\end{equation}

Equations (\ref{eq:parti_co_1}) and (\ref{eq:parti_co_2}) are also
known as the Sylvester equations \cite{Laub2005}, having the unknown
coefficient matrices $\vec{W}_{i}^{\mathcal{E}}$ and $\vec{W}_{i}^{\mathcal{C}}$.
Using the Kronecker product and the vectorization operation

\begin{equation}
\text{vec}(\vec{A})=\text{vec}\left(\left[\vec{a}_{1}\vec{a}_{2}\ldots \vec{a}_{n}\right]\right)=\left[\begin{array}{c}
\vec{a}_{1}\\
\vdots\\
\vec{a}_{n}
\end{array}\right]\in\mathbb{C}^{mn},\quad \vec{A}\in\mathbb{C}^{m\times n},
\end{equation}
with $\vec{a}_{1},\ldots,\vec{a}_{n}$ denoting the column vectors of $\vec{A}$, we
rewrite equation (\ref{eq:parti_co_1}) and equation (\ref{eq:parti_co_2})
as 

\begin{align}
\underbrace{\left(\vec{I}_{2^{i}\times2^{i}}\otimes\vec{\Lambda}_{\mathcal{E}}-\tilde{\vec{\Lambda}}_{\mathcal{E},i}^{T}\otimes \vec{I}_{2\times2}\right)}_{\vec{\Theta}_{i}^{\mathcal{E}}}\text{vec}\left(\vec{W}_{i}^{\mathcal{E}}\right) & =\text{vec}\left(\vec{R}_{i}\right)+\text{vec}\left(\vec{B}_{i}^{\mathcal{E}}\right).\label{eq:cohomological_vec_1}\\
\underbrace{\left(\vec{I}_{2^{i}\times2^{i}}\otimes\vec{\Lambda}_{\mathcal{C}}-\tilde{\vec{\Lambda}}_{\mathcal{E},i}^{T}\otimes \vec{I}_{(2n-2)\times(2n-2)}\right)}_{\vec{\Theta}_{i}^{\mathcal{C}}}\text{vec}\left(\vec{W}_{i}^{\mathcal{C}}\right) & =\text{vec}\left(\vec{B}_{i}^{\mathcal{C}}\right).\label{eq:cohomological_vec_2}
\end{align}

\subsubsection{Invertibility of $\vec{\Theta}_{i}^{\mathcal{E}}$ and $\vec{\Theta}_{i}^{\mathcal{C}}$}

Finding a unique solution for $\vec{W}_{i}^{\mathcal{C}}$ in equation (\ref{eq:cohomological_vec_2})
for a nonzero right-hand side requires the matrix $\vec{\Theta}_{i}^{\mathcal{C}}$
to be non-singular. If, however, $\vec{\Theta}_{i}^{\mathcal{E}}$ is
singular, which arises from exact inner resonances, it suffices for
the vectorized solution for $\vec{W}_{i}^{\mathcal{E}}$ to be in the kernel
of $\vec{\Theta}_{i}^{\mathcal{E}}$. We can ensure this by choosing $\text{vec}\left(\vec{R}_{i}\right)$
such that the right-hand side of equation (\ref{eq:cohomological_vec_1})
is zero. To carry out all this, we need to find the eigenvalues of
$\vec{\Theta}_{i}^{\mathcal{E}}$ and $\vec{\Theta}_{i}^{\mathcal{C}}$. 

It can be shown \cite{Laub2005} that for a matrix $\vec{A}\in\mathbb{C}^{n\times n}$,
with the eigenvalues $\lambda_{i}$ and a matrix $\vec{B}\in\mathbb{C}^{m\times m}$
with eigenvalues $\mu_{j}$, the matrix $\vec{I}_{m\times m}\otimes \vec{A}-\vec{B}\otimes \vec{I}_{n\times n}$
has the $mn$ eigenvalues 

\begin{equation}
\lambda_{1}-\mu_{1},\ldots,\lambda_{1}-\mu_{m},\lambda_{2}-\mu_{1},\ldots,\lambda_{2}-\mu_{m},\ldots,\lambda_{n}-\mu_{1},\ldots,\lambda_{n}-\mu_{m}.\label{eq:direct_sum_prop}
\end{equation}
In the current setting of (\ref{eq:cohomological_vec_1}) and (\ref{eq:cohomological_vec_2}),
we know that the eigenvalues of $\vec{\Lambda}_{\mathcal{E}}$ and $\vec{\Lambda}_{\mathcal{C}}$
are $\lambda_{j_{1}},\lambda_{j_{2}}$ and $\lambda_{j_{3}}\ldots,\lambda_{j_{2n}}$,
respectively. By further exploiting the structure of 

\begin{align}
\tilde{\vec{\Lambda}}_{\mathcal{E},i}^{T}=\tilde{\vec{\Lambda}}_{\mathcal{E},i} & =\sum_{\left|\vec{s}\right|=1}\vec{\Lambda}_{\mathcal{E}}^{s_{1}}\otimes\cdots\otimes\vec{\Lambda}_{\mathcal{E}}^{s_{i}}=\vec{\Lambda}_{\mathcal{E}}\otimes \vec{I}\otimes\cdots\otimes \vec{I}+\ldots+\vec{I}\otimes \vec{I}\otimes\cdots\otimes\vec{\Lambda}_{\mathcal{E}},
\end{align}
or, equivalently,

\begin{equation}
\tilde{\vec{\Lambda}}_{\mathcal{E},i}=\text{diag}\left(a_{1}\lambda_{1}+b_{1}\lambda_{2},\ldots,a_{2^{i}}\lambda_{1}+b_{2^{i}}\lambda_{2}\right)\in\mathbb{C}^{2^{i}\times2^{i}},\label{eq:Ae_diag}
\end{equation}
we observe that each diagonal term of the matrix $\tilde{\vec{\Lambda}}_{\mathcal{E}}$
for a given order $i$ will consist of a linear combination of $\lambda_{1}$
and $\lambda_{2}$, i.e., $a_{j}\lambda_{1}+b_{j}\lambda_{2}$ for
$j=\left\{ 1,\ldots,2^{i}\right\} $, with $a_{j},b_{j}\in\mathbb{N}$.
Now let $\vec{\Omega}_{i}$ be a $2^{i}$-dimensional vector containing
all possible lexicographically ordered $i$-tuples made out of elements
of the set $\left\{ 1,2\right\} $, in which repetition is allowed.
The multiplicity corresponding to the numbers $1$ and $2$, in the
$j^{\text{th}}$ element of $\vec{\Omega}_{i}$, will represent $a_{j}$
and $b_{j}$ respectively. To illustrate this, we give an example 
\begin{example}
\label{ex:Ae} {[}\textit{Constructing the matrix }$\tilde{\vec{\Lambda}}_{\mathcal{E},2}$
{]} For $i=2$, the diagonal matrix $\tilde{\vec{\Lambda}}_{\mathcal{E},2}$
is equal to

\begin{align}
\tilde{\vec{\Lambda}}_{\mathcal{E},2}=\vec{\Lambda}_{\mathcal{E}}\otimes \vec{I}+\vec{I}\otimes\vec{\Lambda}_{\mathcal{E}} & =\left[\begin{array}{cc}
\lambda_{j_{1}} & 0\\
0 & \lambda_{j_{2}}
\end{array}\right]\otimes\left[\begin{array}{cc}
1 & 0\\
0 & 1
\end{array}\right]+\left[\begin{array}{cc}
1 & 0\\
0 & 1
\end{array}\right]\otimes\left[\begin{array}{cc}
\lambda_{j_{1}} & 0\\
0 & \lambda_{j_{2}}
\end{array}\right]\nonumber \\
 & =\left[\begin{array}{cccc}
2\lambda_{j_{1}} & 0 & 0 & 0\\
0 & \lambda_{j_{1}}+\lambda_{j_{2}} & 0 & 0\\
0 & 0 & \lambda_{j_{2}}+\lambda_{j_{1}} & 0\\
0 & 0 & 0 & 2\lambda_{j_{2}}
\end{array}\right].\label{eq:example_kron2}
\end{align}
The four-dimensional array $\vec{\Omega}_{2}$ can be expressed as

\[
\vec{\Omega}_{2}=\left(11,12,21,22\right),
\]
from which we obtain the coefficients $(a_{j},b_{j}$) by determining
the multiplicity of the numbers 1 and 2 for each element $j$ in $\vec{\Omega}_{2}$.
\end{example}
As $\tilde{\vec{\Lambda}}_{\mathcal{E},i}$ is diagonal by construction,
its eigenvalues are positioned on the diagonal and take the form of

\begin{equation}
a\lambda_{j_{1}}+b\lambda_{j_{2}},\quad a,b\in\mathbb{N}:a+b=i.
\end{equation}
We, therefore, conclude from (\ref{eq:direct_sum_prop}) that the
eigenvalues of $\vec{\Theta}_{i}^{\mathcal{E}}$ and $\vec{\Theta}_{i}^{\mathcal{C}}$
can be written as

\begin{align}
\lambda_{l}-(a\lambda_{j_{1}}+b\lambda_{j_{2}}),\quad a,b & \in\mathbb{N}:a+b=i,\quad\forall\lambda_{l}\in\text{Spect}(\Lambda_{\mathcal{E}}),\label{eq:int_non_res}\\
\lambda_{l}-(a\lambda_{j_{1}}+b\lambda_{j_{2}}),\quad a,b & \in\mathbb{N}:a+b=i,\quad\forall\lambda_{l}\in\text{Spect}(\Lambda_{\mathcal{C}}).\label{eq:ext_non_res}
\end{align}
Equations (\ref{eq:int_non_res}) and (\ref{eq:ext_non_res}) lead
precisely to the inner and outer non-resonance conditions (\ref{eq:int_res})
and (\ref{eq:ext_res}), respectively, related to the $i^{\text{th}}$-order
of the SSM expansion. If there exists an inner resonance in (\ref{eq:int_non_res}),
for a particular order $i$, the matrix $\vec{\Theta}_{i}^{\mathcal{E}}$
will be singular. This means that for a nonzero right-hand side of
equation (\ref{eq:cohomological_vec_1}), there will be no solution
for $\text{vec}(\vec{W}_{i}^{\mathcal{E}})$. However, we can then set $\vec{R}_{i}$
in equation (\ref{eq:cohomological_vec_1}) equal to $\text{-}\vec{B}_{i}^{\mathcal{E}}$,
which gives a zero right-hand side. As a consequence, the solution
$\text{vec}(\vec{W}_{i}^{\mathcal{E}})$ has to be in the kernel of $\vec{\Theta}_{i}^{\mathcal{E}}$,
creating an opportunity to remove resonant terms in the expression
for $\vec{W}(\vec{z})$. The presence of an outer resonance in (\ref{eq:ext_non_res})
will result in a breakdown of the SSM. In this case, we do not have
the freedom to alter the right-hand side of equation (\ref{eq:cohomological_vec_2}). 

\section{Reduced dynamics on the SSM \label{sec:The-backbone-curve}}

\subsection{Near-inner-resonances }

Based on eq. (\ref{eq:cohomological_vec_1}), the polynomial dynamics
on the SSM must be parameterized as nonlinear when an inner resonance
arises in the spectral subspace $\mathcal{E}$ over which the SSM
is constructed. When the eigenvalues $\lambda_{j_{1}}$ and $\lambda_{j_{2}}$
are complex conjugate, the inner non-resonance conditions (\ref{eq:int_res})
will never be violated. However, as explained by Szalai et al. \cite{Szalai2017},
for a lightly damped spectral subspace corresponding to a complex
pair of eigenvalues, the following near-inner-resonance conditions
will always hold: 

\begin{equation}
2\lambda_{j_{1}}+\bar{\lambda}_{j_{1}}\approx\lambda_{j_{1}},\quad\lambda_{j_{1}}+2\bar{\lambda}_{j_{1}}\approx\bar{\lambda}_{j_{1}}.\label{eq:third-near-inner-res}
\end{equation}
These near-inner-resonances, in turn, will lead to small denominators
in the coefficients related to the monomial terms $z_{j_{1}}^{2}\bar{z}_{j_{1}}$
and $z_{j_{1}}\bar{z}_{j_{1}}^{2}$ in the third-order coefficient
matrix $\vec{W}_{3}^{\mathcal{E}}$ (cf. eqs. (\ref{eq:cohomological_vec_1}-\ref{eq:cohomological_vec_2})).
Such small denominators generally reduce the domain of convergence
of the Taylor series we compute for $\vec{W}(\vec{z})$. 

Luckily, we have the freedom to remove these resonant terms in $\vec{W}_{3}^{\mathcal{E}}$
by setting $\vec{R}_{3}$ on the right-hand side of equation (\ref{eq:parti_co_1})
equal to $\text{-}\vec{B}_{3}^{\mathcal{E}}$. However, due to the particular
diagonal structure of equation (\ref{eq:cohomological_vec_1}), it
is possible to specifically remove the resonant terms $z_{j_{1}}^{2}\bar{z}_{j_{1}}$
and $z_{j_{1}}\bar{z}_{j_{1}}^{2}$ in $\vec{W}_{3}^{\mathcal{E}}$ by only
setting the coefficients in $\vec{R}_{3}$ related to the the resonant terms
equal to the coefficients in $\text{-}\vec{B}_{3}^{\mathcal{E}}$. This
corresponds to a mixed parameterization style, as explained in Haro
et al. \cite{Haro2016}, which can also be applied to higher orders. 

The third order near-inner-resonance condition (\ref{eq:third-near-inner-res})
can be extended to higher-order near-inner-resonances by introducing
an appropriate resonance-closeness measure

\begin{equation}
I(a,b,\lambda_{l})=\left|\frac{\left\langle \vec{c}(a,b),\vec{v}_{\mathcal{E}}(\lambda_{l})\right\rangle }{\left\Vert \vec{c}(a,b)\right\Vert \left\Vert \vec{v}_{\mathcal{E}}(\lambda_{l})\right\Vert }\right|<\delta,\quad0<\delta\ll1,\label{eq:int_res_measure}
\end{equation}
for $\delta$ sufficiently small, with

\begin{gather*}
\vec{c}(a,b)=\left[\begin{array}{c}
a\\
b\\
-1
\end{array}\right],\quad \vec{v}_{\mathcal{E}}(\lambda_{l})=\left[\begin{array}{c}
\lambda_{j_{1}}\\
\lambda_{j_{2}}\\
\lambda_{l}
\end{array}\right],\quad a,b\in\mathbb{N},\quad\forall\lambda_{l}\in\left\{ \lambda_{j_{1}},\lambda_{j_{2}}\right\} .
\end{gather*}
The resonance-closeness measure $I$ takes values between $I_\text{min}=0$ and $I_\text{max}=1$. We consider $\delta$ to be small when $\delta$ is at least one order of magnitude smaller than $I_\text{max}$. In the presence of near-inner resonances, the choice of $\delta$ affects the accuracy of the SSM and the reduced dynamics. If the observed invariance error of the SSM (cf. section \ref{sec:invarselect}) is unsatisfactory, $\delta$ can be increased in order to account for even weaker near-inner resonances. 

In the case of an exact inner resonance, $I(a,b,\lambda_{l})$ in
equation (\ref{eq:int_res_measure}) will be zero. Using the same
measure, we can also quantify closeness to outer resonances by substituting
all possible $\lambda_{l}\notin\left\{ \lambda_{j_{1}},\lambda_{j_{2}}\right\} $
into (\ref{eq:int_res_measure}). 

\subsection{Instantaneous amplitude and frequency }

When the chosen spectral subspace $\mathcal{E}$ is spanned by a complex
pair of eigenvectors, which in turn corresponds to a complex conjugate
pair of eigenvalues $\lambda_{j_{1}}$ and $\lambda_{j_{2}}$, the
complex conjugate pair of coordinates $z_{j_{1}}$ and $z_{j_{2}}=\bar{z}_{j_{1}}$
in the reduced dynamics $\vec{R}(\vec{z})$ can be expressed in real amplitude-phase
coordinates $(\rho,\theta)$ as

\begin{align}
z_{j_{1}} & =\rho \mathrm{e}^{\mathrm{i}\theta},\quad\bar{z}_{j_{1}}=\rho \mathrm{e}^{-\mathrm{i}\theta}.\label{eq:polar_cor}
\end{align}
Assume now that the spectral subspace $\mathcal{E}$ has higher-order
near-inner-resonances, i.e.

\begin{gather*}
I(a,b,\lambda_{j_{1}})<\delta,\quad I(b,a,\lambda_{j_{2}})<\delta,\quad(a,b)\in\left\{ (2,1),(3,2),(4,3),\ldots\right\} =S,
\end{gather*}
and, additionally, the coefficients in $\vec{B}_{i}^{\mathcal{E}}$ for
$i=3,5,7\ldots$, corresponding to the monomial terms $z_{j_{1}}^{a}\bar{z}_{j_{1}}^{b}$
and $z_{j_{1}}^{s}\bar{z}_{j_{1}}^{r}$, on the right hand side of
(\ref{eq:parti_co_1}) are nonzero. We then obtain the following expression
for the reduced dynamics on the spectral submanifold $\mathcal{W}(\mathcal{E})$:

\begin{equation}
\dot{\vec{z}}=\vec{R}(\vec{z})=\left[\begin{array}{c}
\lambda_{j_{1}}z_{j_{1}}+{\displaystyle \sum_{\forall(a,b)\in S}\gamma_{a,b}z_{j_{1}}^{a}\bar{z}_{j_{1}}^{b}}\\
\bar{\lambda}_{j_{1}}\bar{z}_{j_{1}}+{\displaystyle \sum_{\forall(a,b)\in S}\bar{\gamma}_{a,b}z_{j_{1}}^{b}\bar{z}_{j_{1}}^{a}}
\end{array}\right].\label{eq:red_dyn_mixed}
\end{equation}
Here $\gamma_{a,b}$ depends directly on $\vec{B}_{i}^{\mathcal{E}}$, which
is known for the current order $i$. Alternatively, $\gamma_{a,b}$ will be equal to the sum of all nonzero coefficients in $\vec{R}_i$, with $i=a+b$ corresponding to the monomial term $z_{j_{1}}^{a}\bar{z}_{j_{1}}^{b}$.
Substituting equation (\ref{eq:polar_cor}) into the left- and right-hand side of equation (\ref{eq:red_dyn_mixed}) gives

\begin{align}
\dot{\rho} & =\text{Re}(\lambda_{j_{1}})\rho+\sum_{\forall(a,b)\in S}\text{Re}(\gamma_{a,b})\rho^{(a+b)},\\
\omega & =\dot{\theta}=\text{Im}(\lambda_{j_{1}})+\sum_{\forall(a,b)\in S}\text{Im}(\gamma_{a,b})\rho^{(a+b-1)}.\label{eq:instant_freq}
\end{align}
for $\rho\neq0$. Equation (\ref{eq:instant_freq}) determines an
instantaneous frequency $\dot{\theta}$ that depends solely upon $\rho$.
To any $\rho$ value, we assign an instantaneous physically observable
amplitude by defining 

\begin{equation}
A(\rho)=\frac{1}{2\pi}\intop_{0}^{2\pi}\left|\vec{T}_{\vec{y}}\vec{W}(\vec{z}(\rho,\theta))\right|d\theta,\label{eq:instant_amp}
\end{equation}
where the transformation matrix $\vec{T}_{\vec{y}}$ acts on $\vec{W}(\vec{z}(\rho,\theta))$
and hence returns physical position coordinates $\vec{y}\in\mathbb{R}^{n}$
of our mechanical system (\ref{eq:dyn_sys}). Then, following the
definition of Szalai et al. \cite{Szalai2017}, we define a backbone
curve for the reduced dynamics on the SSM to be the parameterized
curve 

\begin{equation}
\mathcal{B}=\left\{ \omega(\rho),A(\rho)\right\} _{\rho\in\mathbb{R}^{+}}.
\end{equation}

\begin{figure}[H]
\begin{centering}
\includegraphics[scale=0.45]{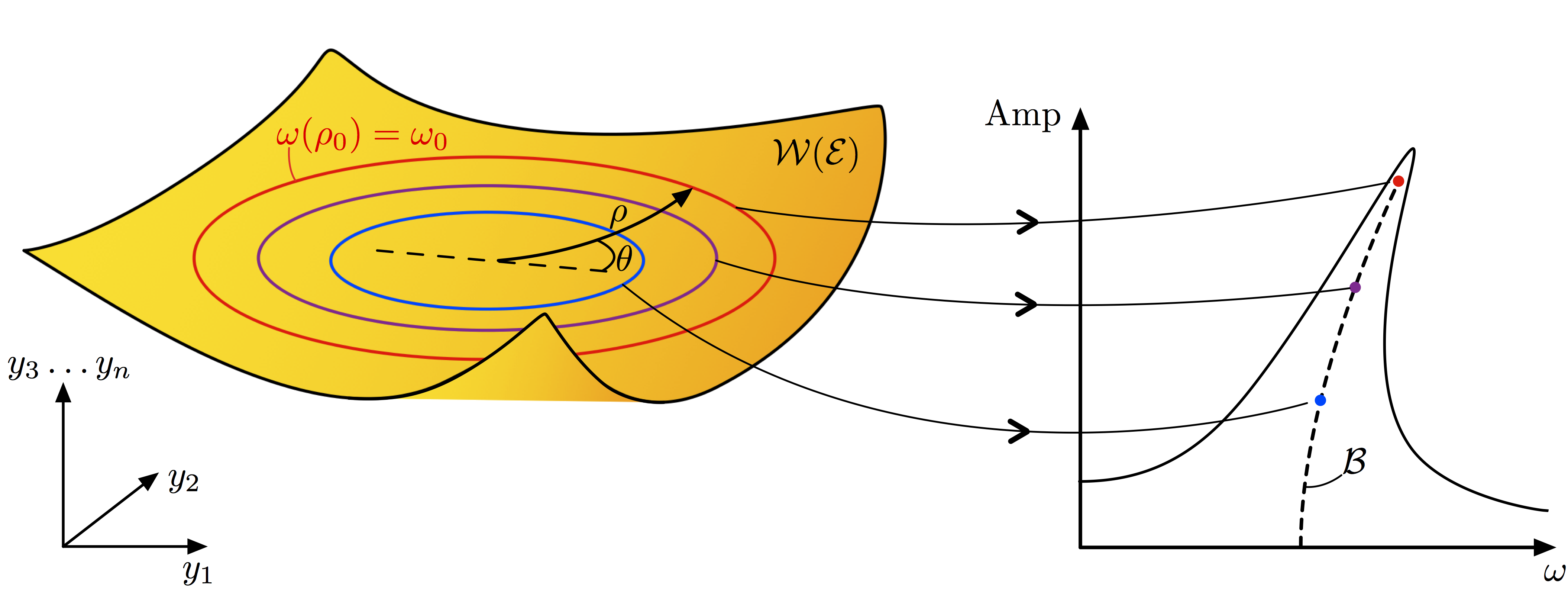}
\par\end{centering}
\caption{Illustration of the backbone curve construction using the parameterized
SSM. For each fixed radius $\rho_{0}$ on the SSM, we can identify
an instantaneous frequency (\ref{eq:instant_freq}). By averaging
the physical coordinates over one period, evaluated on the SSM constraint
to the fixed radius $\rho_{0}$, we obtain the instantaneous amplitude
(\ref{eq:instant_amp}). For each $\rho=\rho_{0}\in\mathbb{R}^{+}$,
we obtain a point of the curve $\mathcal{B}$, shown in the amplitude-frequency
plot. The continuous black line in the amplitude-frequency plot consists
out of periodic orbits of the periodically forced system for varying
forcing frequency. \label{fig:Illustration-of-the-backbone}}
\end{figure}

An illustration of how the parametrized SSM, constructed over a lightly
damped spectral subspace $\mathcal{E}$, can be used to construct
the backbone curve $\mathcal{B}$ is shown in figure \ref{fig:Illustration-of-the-backbone}. 

\section{Invariance measure and order selection \label{sec:invarselect}}

\begin{figure}[H]
\begin{centering}
\includegraphics[scale=0.35]{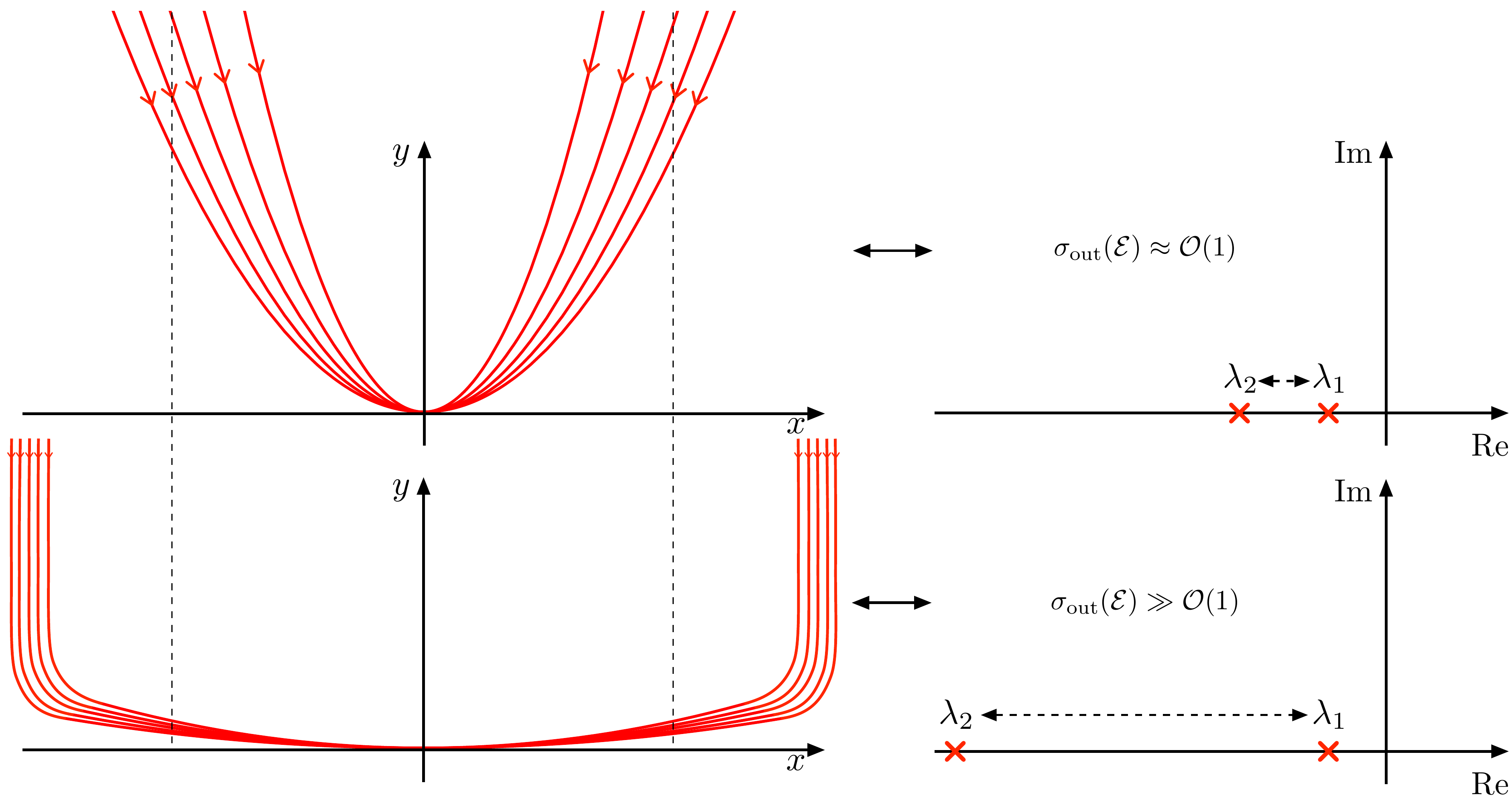}
\par\end{centering}
\caption{Illustration of a system having a low outer spectral quotient (top)
and a high outer spectral quotient (bottom).\label{fig:order_taylor_outer}}
\end{figure}

As stated in Theorem 1, the unique SSM is captured and approximated
by Taylor expanding up to order $\sigma_{\text{out}}(\mathcal{E})+1$.
The outer spectral quotient $\sigma_{\text{out}}$, is defined as
the integer part of the ratio between the strongest decay rate of
the linearized oscillations outside $\mathcal{E}$ and the weakest
decay rate of the linearized oscillations inside $\mathcal{E}$. As
explained by Géradin and Rixen \cite{Geradin2014}, a first-order
approximation of the real part of each eigenvalue of a lightly damped
mechanical system, of the form (\ref{eq:mech_sys}), scales with the
square of its natural eigenfrequency. This means that for a discretized
non-conservative mechanical system with a high number of degrees of
freedom, the order of the SSM needed to be unique can become large.
An illustration of this is shown in figure \ref{fig:order_taylor_outer}.

As the outer spectral quotient $\sigma_{\text{out}}(\mathcal{E})$
increases, trajectories transverse to the slow SSM die out fast compared
to trajectories on the SSM, as indicated in the lower part of figure
\ref{fig:order_taylor_outer}. This collapsing nature of the transverse
trajectories makes it harder to distinguish between the unique SSM
and any other two-dimensional invariant manifold tangent the same
modal subspace. In order to approximate the SSMs with a large outer
spectral quotient, without having to compute the SSMs up to extremely
high orders, we introduce an invariance error measure, $\delta_{\text{inv}}$,
that quantifies the accuracy of the computed invariant manifolds and
the reduced dynamics on them. 

\begin{figure}[H]
\begin{centering}
\includegraphics[scale=0.42]{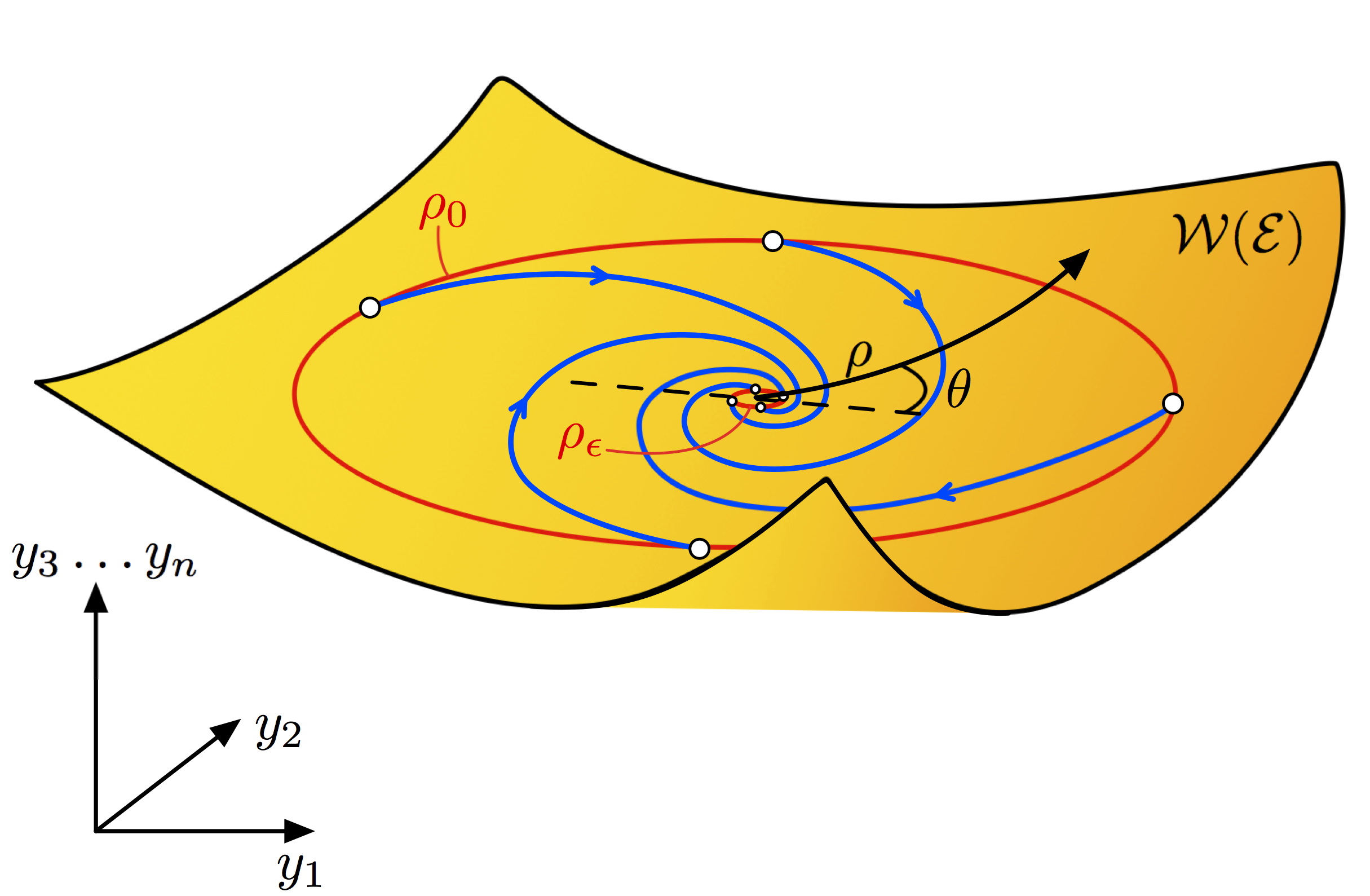}
\par\end{centering}
\caption{Illustration of the invariance error measure. \label{fig:illu-invar-measure}}
\end{figure}

The invariance error measure compares trajectories of the full system
$\vec{x}_{i}$, with trajectories of the reduced system $\tilde{\vec{x}}_{i}$.
Trajectories from the full and reduced system are launched from a
circle with fixed radius $\rho_{0}$, from the origin, and integrated
until the reduced trajectories cross the inner circle of radius $\rho_{\epsilon}<\rho_{0}$,
therefore removing the time dependency. An illustration of this is
shown in figure \ref{fig:illu-invar-measure}. We mathematically formalize
the invariance error as follows 

\begin{equation}
\delta_{\text{inv}}=\frac{1}{N}\sum_{i=1}^{N}\frac{\text{dist}(i)}{\underset{\theta\in S^{1}}{\text{max}}\left\Vert \tilde{\vec{x}}(\rho_{0},\theta)\right\Vert _{2}},\quad\text{dist}(i)=\text{max}\left\Vert \left.\vec{x}_{i}\right|_{\rho_{0}}^{\rho_{\epsilon}}-\left.\tilde{\vec{x}}_{i}\right|_{\rho_{0}}^{\rho_{\epsilon}}\right\Vert _{2},\label{eq:inv_error_measure-1}
\end{equation}
where we take the average of the maximum Euclidean distance between
$N$ trajectories $\vec{x}_{i}$ and $\tilde{\vec{x}}_{i}$, for $i=1,\ldots,N$,
traveling from a circle with radius $\rho_{0}$ to an inner circle
with radius $\rho_{\epsilon}$, and normalize the result by the maximum
Euclidean distance from the origin to the circle with fixed radius
$\rho_{0}$.

If the invariance error $\delta_{\text{inv}}$, for a given order
of the approximated SSM, is above a certain pre-specified bound, then
the order of the SSM approximation has to be further increased. 

\section{Applications }

We now apply our computational algorithm to three different mechanical
systems. The numerical results and figures we show have all been generated
directly by SSMtool. 

\subsection{The modified Shaw\textendash Pierre example: Inner resonances \label{subsec:The-modified-Shaw=002013Pierre-internal}}

\begin{figure}[H]
\begin{centering}
\includegraphics[scale=0.42]{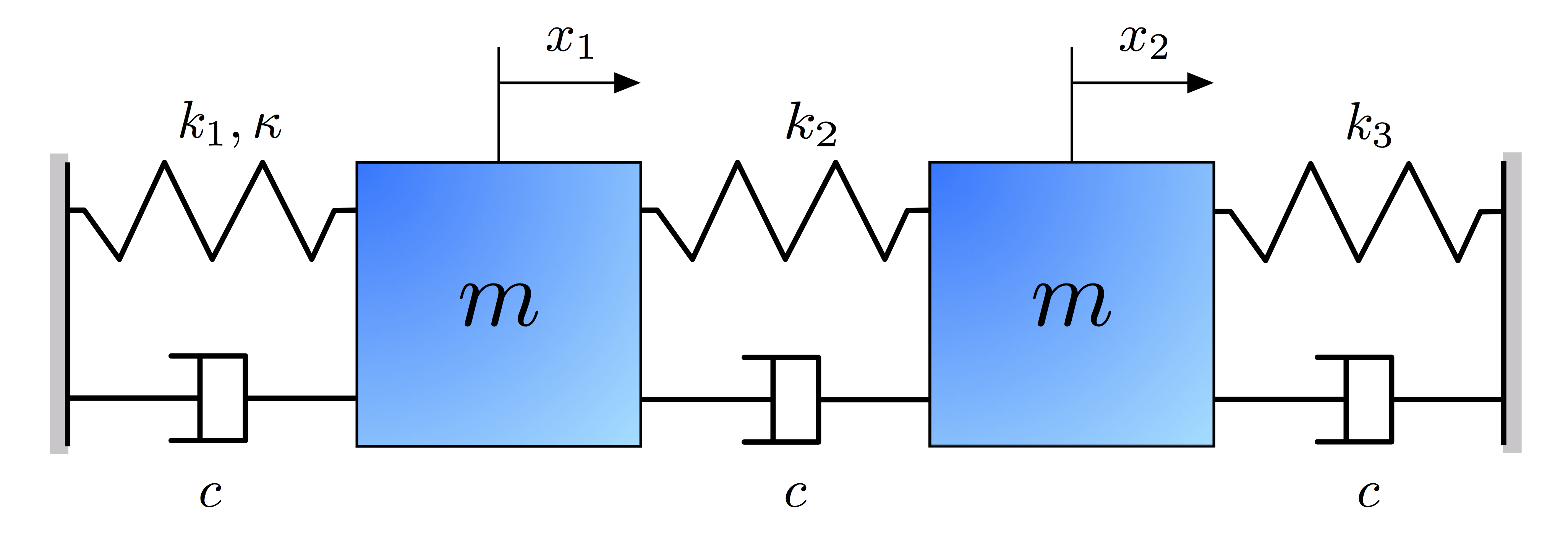}
\par\end{centering}
\caption{Two-degree-of-freedom modified Shaw\textendash Pierre example.\label{fig:Two-degree-of-freedom-modified}}

\end{figure}

We first consider a slightly modified version of the example of Shaw
and Pierre \cite{Shaw1994}, shown in figure \ref{fig:Two-degree-of-freedom-modified}.
The original Shaw\textendash Pierre example involves a two-degree-of-freedom
mechanical oscillator, which is modified in the current setting such
that the damping matrix is proportional to the mass and stiffness
matrices (also known as Rayleigh damping, see, e.g., Géradin and Rixen
\cite{Geradin2014}). For this problem, the SSM coefficients have
been explicitly calculated in Szalai et al. \cite{Szalai2017}, up
to third order.

For $k_{1}=k_{2}=k_{3}=k$, the equations of motion of the system
are 

\begin{equation}
\left[\begin{array}{cc}
m & 0\\
0 & m
\end{array}\right]\left[\begin{array}{c}
\ddot{x}_{1}\\
\ddot{x}_{2}
\end{array}\right]+\left[\begin{array}{cc}
2c & -c\\
-c & 2c
\end{array}\right]\left[\begin{array}{c}
\dot{x}_{1}\\
\dot{x}_{2}
\end{array}\right]+\left[\begin{array}{cc}
2k & -k\\
-k & 2k
\end{array}\right]\left[\begin{array}{c}
x_{1}\\
x_{2}
\end{array}\right]+\left[\begin{array}{c}
\kappa x_{1}^{3}\\
0
\end{array}\right]=\left[\begin{array}{c}
0\\
0
\end{array}\right],\label{eq:sys_mod_sp}
\end{equation}
with the linear part having the eigenvalue pairs

\begin{equation}
\lambda_{1,2}=-\frac{c}{2}\pm \mathrm{i}\sqrt{k\left(1-\frac{c^{2}}{4k}\right)},\quad\lambda_{3,4}=-\frac{3c}{2}\pm \mathrm{i}\sqrt{3k\left(1-\frac{3c^{2}}{4k}\right)},
\end{equation}
when both linear normal modes are underdamped ($c<2\sqrt{k/3}$) and
the mass $m$ is equal to $1\text{ kg}$. As noted in Szalai et al.
\cite{Szalai2017}, the two spectral subspaces, $\mathcal{E}_{1}$
and $\mathcal{E}_{2}$, corresponding to the eigenvalues $\lambda_{1,2}$
and $\lambda_{3,4}$ respectively, have the outer and inner spectral
quotients

\begin{align}
\sigma_{\text{out}}(\mathcal{E}_{1}) & =\text{Int}\left[\frac{\text{Re}\lambda_{3}}{\text{Re}\lambda_{1}}\right]=3,\quad\sigma_{\text{out}}(\mathcal{E}_{2})=\left[\frac{\text{Re}\lambda_{1}}{\text{Re}\lambda_{3}}\right]=0,\\
\sigma_{\text{in}}(\mathcal{E}_{1}) & =\text{Int}\left[\frac{\text{Re}\lambda_{1}}{\text{Re}\lambda_{1}}\right]=1,\quad\sigma_{\text{in}}(\mathcal{E}_{2})=\left[\frac{\text{Re}\lambda_{3}}{\text{Re}\lambda_{3}}\right]=1.
\end{align}
The non-resonance conditions (\ref{eq:ext_res}) and (\ref{eq:int_res})
are satisfied for both of these spectral subspaces, thus, there exist
two two-dimensional analytic SSMs, $\vec{W}(\mathcal{E}_{1})$
and $\vec{W}(\mathcal{E}_{2})$, that are unique among all $C^{4}$
and $C^{1}$ invariant manifolds tangent to $\mathcal{E}_{1}$ and
$\mathcal{E}_{2}$, respectively. 

Rewriting the equations of motion (\ref{eq:sys_mod_sp}) in first-order
form, we obtain

\begin{equation}
\frac{d}{dt}\left[\begin{array}{c}
x_{1}\\
x_{2}\\
\dot{x}_{1}\\
\dot{x}_{2}
\end{array}\right]=\underbrace{\left[\begin{array}{cccc}
0 & 0 & 1 & 0\\
0 & 0 & 0 & 1\\
-2k & k & -2c & c\\
k & -2k & c & -2c
\end{array}\right]}_{\vec{A}}\left[\begin{array}{c}
x_{1}\\
x_{2}\\
\dot{x}_{1}\\
\dot{x}_{2}
\end{array}\right]+\underbrace{\left[\begin{array}{c}
0\\
0\\
-\kappa x_{1}^{3}\\
0
\end{array}\right]}_{\vec{F}(\vec{x})}.\label{eq:mod_SP_first_order}
\end{equation}

\subsubsection{Computing $\vec{W}(\mathcal{E}_{1})$ and $\vec{W}(\mathcal{E}_{2})$ \label{subsec:Computing_W(E1)}}

The spectral submanifolds, $\vec{W}(\mathcal{E}_{1})$ and $\vec{W}(\mathcal{E}_{2})$,
will be tangent to their corresponding spectral subspaces, $\mathcal{E}_{1}$
and $\mathcal{E}_{2}$. To compute $\vec{W}(\mathcal{E}_{1})$, we diagonalize
(\ref{eq:mod_SP_first_order}) by introducing a linear change of coordinates
$\vec{x}=\vec{T}\vec{q}$, where the columns of $\vec{T}$ contain the eigenvectors of (\ref{eq:mod_SP_first_order}),
i.e.,

\begin{equation}
\vec{T}=\left[\vec{v}_{1},\bar{\vec{v}}_{1},\vec{v}_{3},\bar{\vec{v}}_{3}\right]=\left[\begin{array}{ccrr}
1 & 1 & 1 & 1\\
1 & 1 & -1 & -1\\
\lambda_{1} & \bar{\lambda}_{1} & \lambda_{3} & \bar{\lambda}_{3}\\
\lambda_{1} & \bar{\lambda}_{1} & -\lambda_{3} & -\bar{\lambda}_{3}
\end{array}\right].
\end{equation}
We can rewrite equation (\ref{eq:mod_SP_first_order}) in the form
of (\ref{eq:ds_diag_repeat})

\begin{equation}
\dot{\vec{q}}=\text{diag}(\lambda_{1},\bar{\lambda}_{1},\lambda_{3},\bar{\lambda}_{3})\vec{q}+\vec{T}^{-1}\vec{F}(\vec{T}\vec{q})=\vec{\Lambda} \vec{q}+\vec{G}(\vec{q}).
\end{equation}
To compute $\vec{W}(\mathcal{E}_{2})$, equation (\ref{eq:mod_SP_first_order})
must be diagonalized via a similar linear change of coordinates $\vec{x}=\tilde{\vec{T}}\tilde{\vec{q}}$,
where the columns of $\tilde{\vec{T}}$ now contain the eigenvectors of
(\ref{eq:mod_SP_first_order}) in the following order

\begin{equation}
\vec{T}=\left[\vec{v}_{3},\bar{\vec{v}}_{3},\vec{v}_{1},\bar{\vec{v}}_{1}\right]=\left[\begin{array}{rrcc}
1 & 1 & 1 & 1\\
-1 & -1 & 1 & 1\\
\lambda_{3} & \bar{\lambda}_{3} & \lambda_{1} & \bar{\lambda}_{1}\\
-\lambda_{3} & -\bar{\lambda}_{3} & \lambda_{1} & \bar{\lambda}_{1}
\end{array}\right].
\end{equation}
Similarly, equation (\ref{eq:mod_SP_first_order}) can be written
in the form of (\ref{eq:ds_diag_repeat}), i.e.

\begin{equation}
\dot{\tilde{\vec{q}}}=\text{diag}(\lambda_{3},\bar{\lambda}_{3},\lambda_{1},\bar{\lambda}_{1})\tilde{\vec{q}}+\tilde{\vec{T}}^{-1}\vec{F}(\tilde{\vec{T}}\tilde{\vec{q}})=\tilde{\vec{\Lambda}}\tilde{\vec{q}}+\tilde{\vec{G}}(\tilde{\vec{q}}).
\end{equation}
The polynomial expressions for the nonlinearities $\vec{G}(\vec{q})$ and $\tilde{\vec{G}}(\tilde{\vec{q}})$
only contain cubic nonlinearities and therefore only the nonlinear
coefficient matrices $\vec{G}_{3}$ and $\tilde{\vec{G}}_{3}$ will be non-zero.
We will compute $\vec{W}(\mathcal{E}_{1})$ and $\vec{W}(\mathcal{E}_{2})$ for
the following parameter values

\begin{equation}
k=1 \text{ N }\text{m}^{-1},\quad c=0.03\text{ N s }\text{m}^{-1},\quad\kappa=0.5\text{ N }\text{m}^{-3},\quad\delta=0.05.\label{eq:par-setting}
\end{equation}

We justify the choice up to which order we have to approximate the
SSMs to get an accurate reduced order model, by evaluating the invariance
error (\ref{eq:inv_error_measure-1}) for different approximation
orders. For a given fixed radius $\rho_{0}=0.35$ we take $50$ initial
points, each corresponding to an angle $\theta_{0}$, uniformly distributed
in $S^{1}$. 

\begin{figure}[H]
\begin{centering}
\subfloat[\label{fig:invar_mode_1}]{\centering{}\includegraphics[scale=0.25]{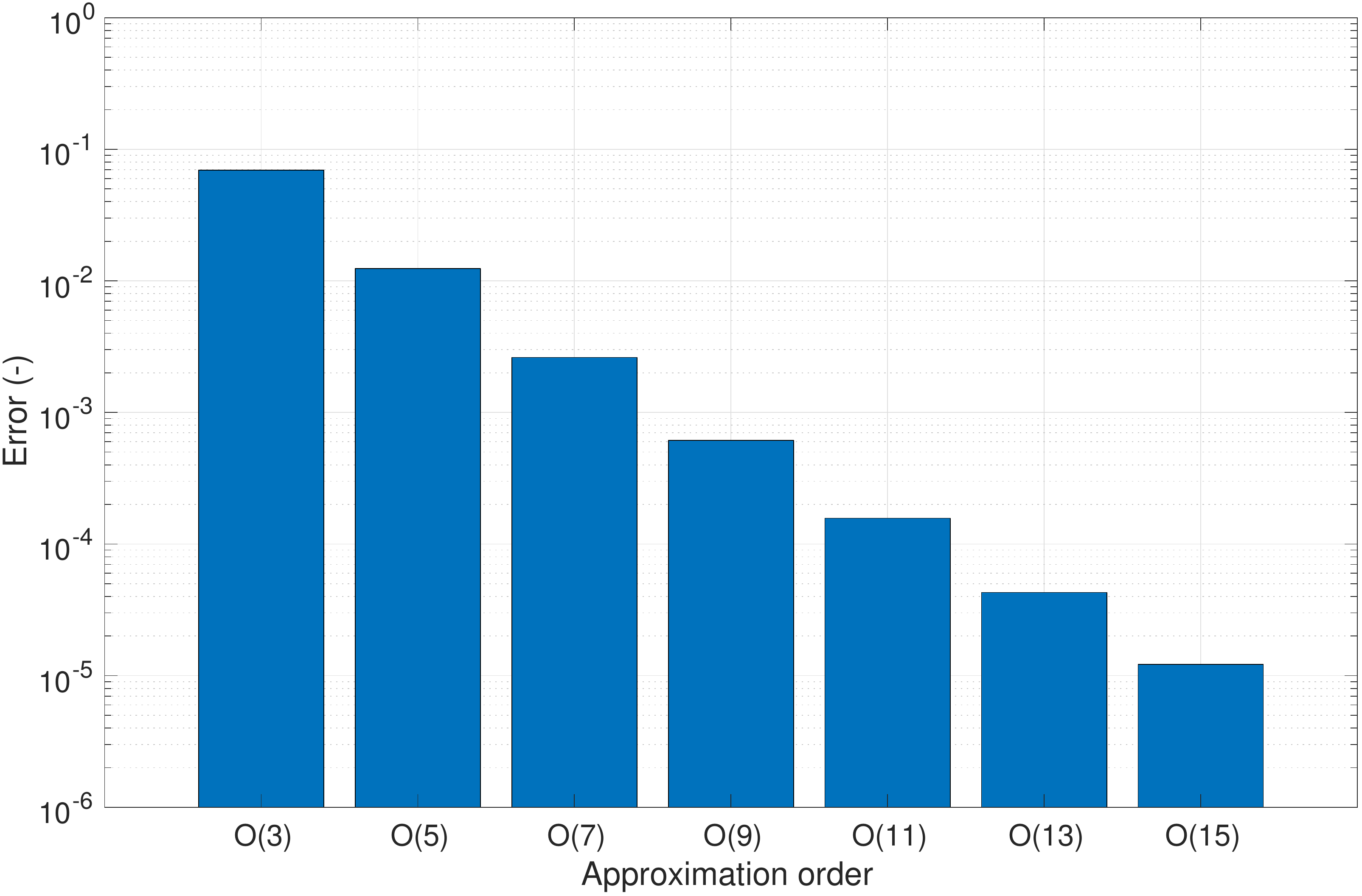}}\hfill{}\subfloat[\label{fig:invar_mode_2}]{\begin{centering}
\includegraphics[scale=0.25]{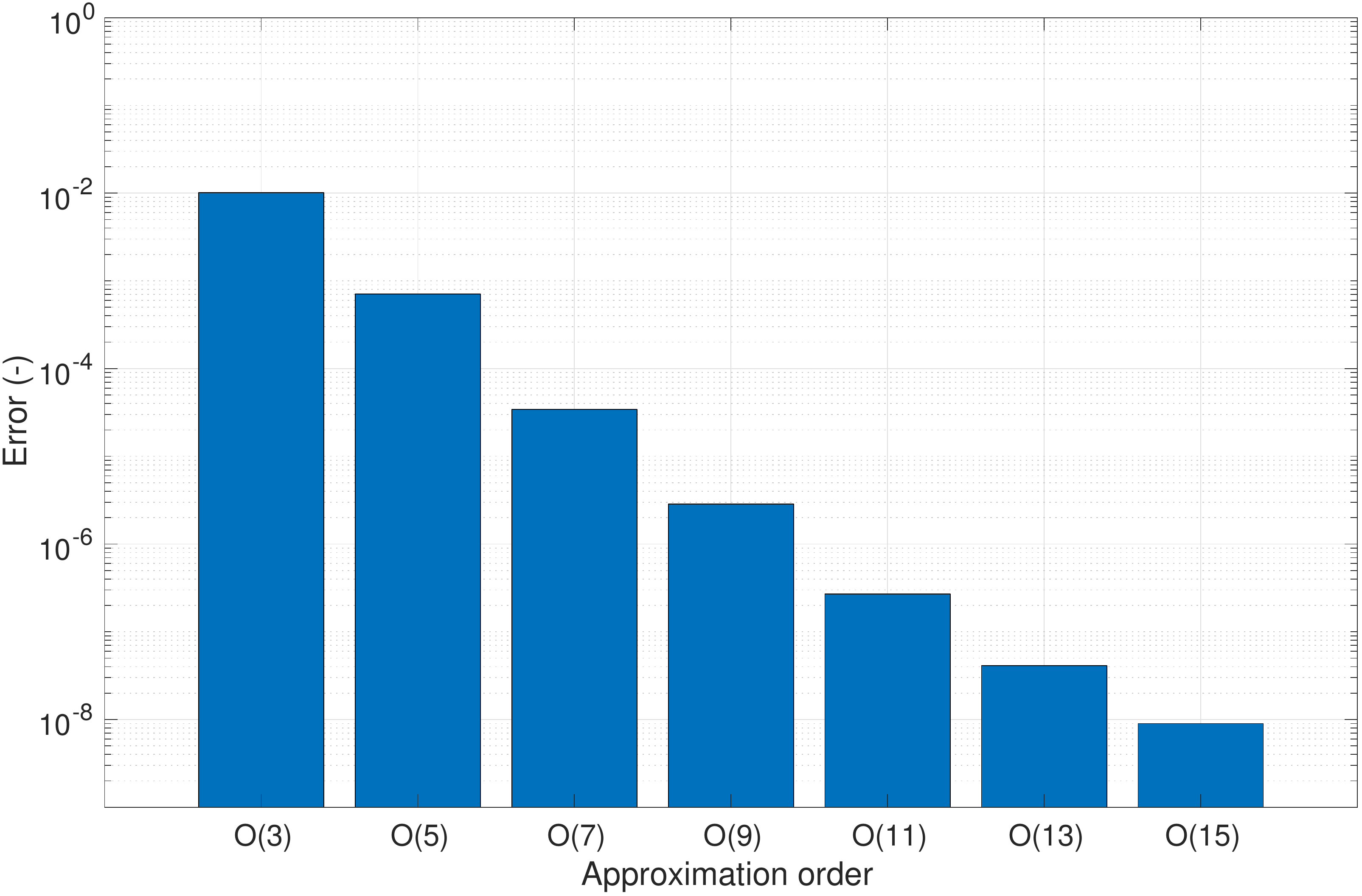}
\par\end{centering}
}
\par\end{centering}
\caption{Normalized error for the $3^{\text{th}}$-$15^{\text{th}}$ order
approximation of $\vec{W}(\mathcal{E}_{1})$ (\ref{fig:invar_mode_1}) and
$\vec{W}(\mathcal{E}_{2})$ (\ref{fig:invar_mode_2}) for 50 evenly distributed
initial positions lying on a fixed radius $\rho_{0}=0.35$ and $\theta_{0}\in S^{1}$.
For each trajectory traveling between $\rho_{0}$ and $\rho_{\epsilon}=0.01$,
we identify the maximum error and take the average over all trajectories.
\label{fig:SSM_invar_error_2DOF}}
\end{figure}

As the order of the approximation of $\vec{W}(\mathcal{E}_{1})$ and $\vec{W}(\mathcal{E}_{2})$
is increased, the error $\delta_{\text{inv}}$ is substantially reduced,
as expected. We conclude that the $15^{\text{th}}$ order approximation
for both spectral submanifolds is high enough to guarantee them to
be accurate for oscillation amplitudes up to $\rho_{0}=0.35$, which
corresponds to a physical maximum displacement of $\left|x_{1}\right|\approx0.66\text{ m}$
and $\left|x_{2}\right|\approx0.71\text{ m}$ for $\vec{W}(\mathcal{E}_{1})$
and $\left|x_{1}\right|\approx0.73\text{ m}$ and $\left|x_{2}\right|\approx0.66\text{ m}$
for $\vec{W}(\mathcal{E}_{2})$. 

We observe that the following near-inner-resonances conditions are
satisfied within the spectral subspaces $\mathcal{E}_{1}$ and $\mathcal{E}_{2}$,
see table (\ref{tab:int_res_E_1}) and (\ref{tab:int_res_E_2}) respectively.

\begin{table}[H]
\begin{centering}
\subfloat[\label{tab:int_res_E_1}]{\begin{centering}
\begin{tabular}{ccccc}
\toprule 
\multicolumn{5}{c}{$\mathcal{E}_{1}$}\tabularnewline
\midrule 
 & $a$ & $b$ & $\lambda_{l}$ & $I$\tabularnewline
\midrule
\multirow{2}{*}{$\mathcal{O}(\left|\vec{z}\right|^{3})$} & 2 & 1 & $\lambda_{1}$ & 0.00707\tabularnewline
 & 1 & 2 & $\bar{\lambda}_{1}$ & 0.00707\tabularnewline
\midrule
\multirow{2}{*}{$\mathcal{O}(\left|\vec{z}\right|^{5})$} & 3 & 2 & $\lambda_{1}$ & 0.00926\tabularnewline
 & 2 & 3 & $\bar{\lambda}_{1}$ & 0.00926\tabularnewline
\midrule
\multirow{2}{*}{$\mathcal{O}(\left|\vec{z}\right|^{7})$} & 4 & 3 & $\lambda_{1}$ & 0.01019\tabularnewline
 & 3 & 4 & $\bar{\lambda}_{1}$ & 0.01019\tabularnewline
\midrule
\multirow{2}{*}{$\mathcal{O}(\left|\vec{z}\right|^{9})$} & 5 & 4 & $\lambda_{1}$ & 0.01069\tabularnewline
 & 4 & 5 & $\bar{\lambda}_{1}$ & 0.01069\tabularnewline
\midrule
\multirow{2}{*}{$\mathcal{O}(\left|\vec{z}\right|^{11})$} & 6 & 5 & $\lambda_{1}$ & 0.01100\tabularnewline
 & 5 & 6 & $\bar{\lambda}_{1}$ & 0.01100\tabularnewline
\midrule
\multirow{2}{*}{$\mathcal{O}(\left|\vec{z}\right|^{13})$} & 7 & 6 & $\lambda_{1}$ & 0.01121\tabularnewline
 & 6 & 7 & $\bar{\lambda}_{1}$ & 0.01121\tabularnewline
\midrule
\multirow{2}{*}{$\mathcal{O}(\left|\vec{z}\right|^{15})$} & 8 & 7 & $\lambda_{1}$ & 0.01136\tabularnewline
 & 7 & 8 & $\bar{\lambda}_{1}$ & 0.01136\tabularnewline
\bottomrule
\end{tabular}
\par\end{centering}
\begin{centering}
\par\end{centering}
}\,\subfloat[\label{tab:int_res_E_2}]{\begin{centering}
\begin{tabular}{ccccc}
\toprule 
\multicolumn{5}{c}{$\mathcal{E}_{2}$}\tabularnewline
\midrule 
 & $a$ & $b$ & $\lambda_{l}$ & $I$\tabularnewline
\midrule
\multirow{2}{*}{$\mathcal{O}(\left|\vec{z}\right|^{3})$} & 2 & 1 & $\lambda_{3}$ & 0.01225\tabularnewline
 & 1 & 2 & $\bar{\lambda}_{3}$ & 0.01225\tabularnewline
\midrule
\multirow{2}{*}{$\mathcal{O}(\left|\vec{z}\right|^{5})$} & 3 & 2 & $\lambda_{3}$ & 0.01604\tabularnewline
 & 2 & 3 & $\bar{\lambda}_{3}$ & 0.01604\tabularnewline
\midrule
\multirow{2}{*}{$\mathcal{O}(\left|\vec{z}\right|^{7})$} & 4 & 3 & $\lambda_{3}$ & 0.01765\tabularnewline
 & 3 & 4 & $\bar{\lambda}_{3}$ & 0.01765\tabularnewline
\midrule
\multirow{2}{*}{$\mathcal{O}(\left|\vec{z}\right|^{9})$} & 5 & 4 & $\lambda_{3}$ & 0.01852\tabularnewline
 & 4 & 5 & $\bar{\lambda}_{3}$ & 0.01852\tabularnewline
\midrule
\multirow{2}{*}{$\mathcal{O}(\left|\vec{z}\right|^{11})$} & 6 & 5 & $\lambda_{3}$ & 0.01905\tabularnewline
 & 5 & 6 & $\bar{\lambda}_{3}$ & 0.01905\tabularnewline
\midrule
\multirow{2}{*}{$\mathcal{O}(\left|\vec{z}\right|^{13})$} & 7 & 6 & $\lambda_{3}$ & 0.01941\tabularnewline
 & 6 & 7 & $\bar{\lambda}_{3}$ & 0.01941\tabularnewline
\midrule
\multirow{2}{*}{$\mathcal{O}(\left|\vec{z}\right|^{15})$} & 8 & 7 & $\lambda_{3}$ & 0.01967\tabularnewline
 & 7 & 8 & $\bar{\lambda}_{3}$ & 0.01967\tabularnewline
\bottomrule
\end{tabular}
\par\end{centering}
\centering{}}
\par\end{centering}
\caption{Near-inner-resonances for $\mathcal{E}_{1}$ and $\mathcal{E}_{2}$
with $\delta=0.05$.\label{tab:Near-internal-resonances_E1}}

\end{table}

\begin{figure}[H]
\begin{centering}
\subfloat[\label{subfig:SSM_1_q3}]{
\includegraphics[scale=0.16]{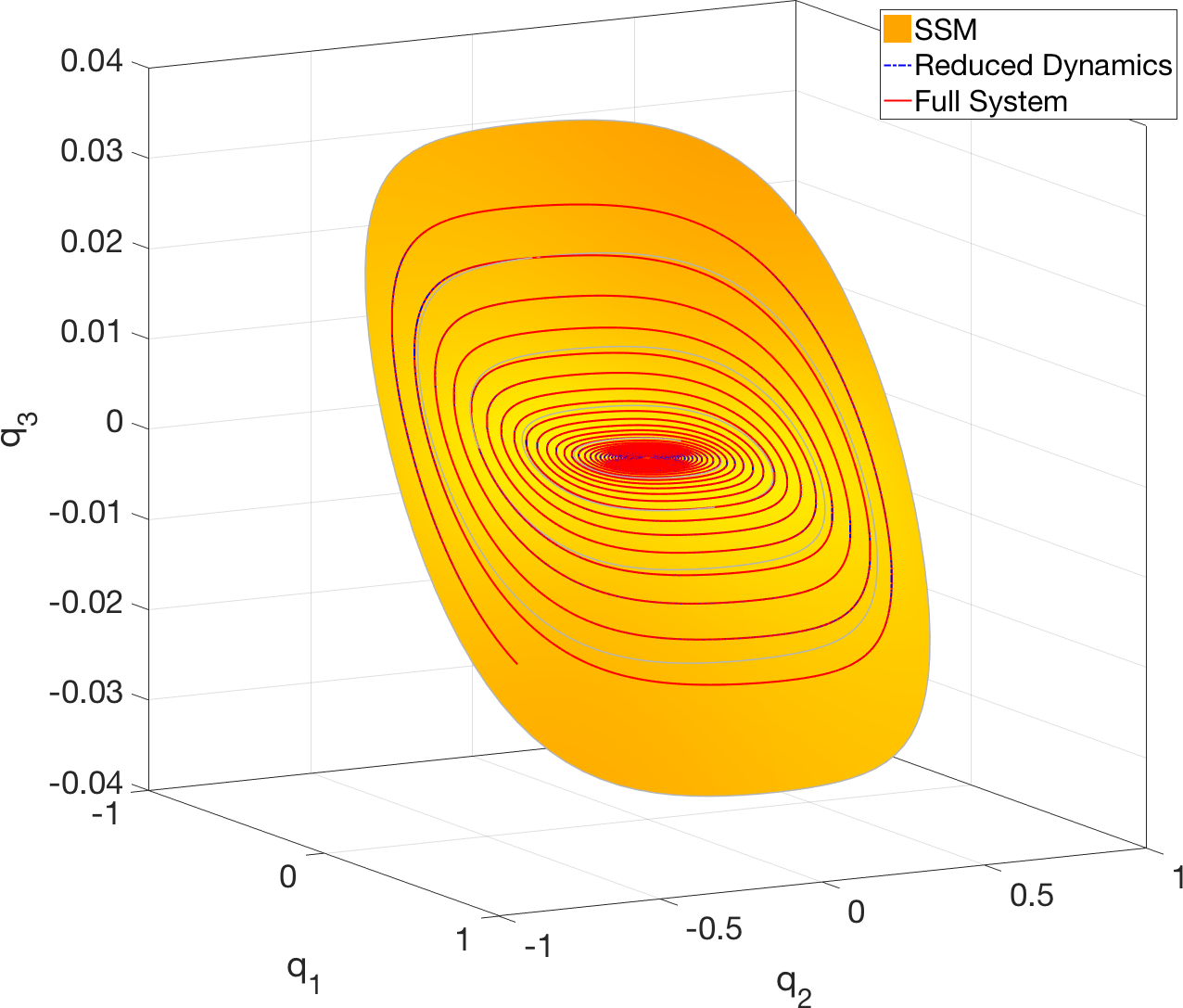}}
\subfloat[\label{subfig:SSM_1_q4}]{
\includegraphics[scale=0.17]{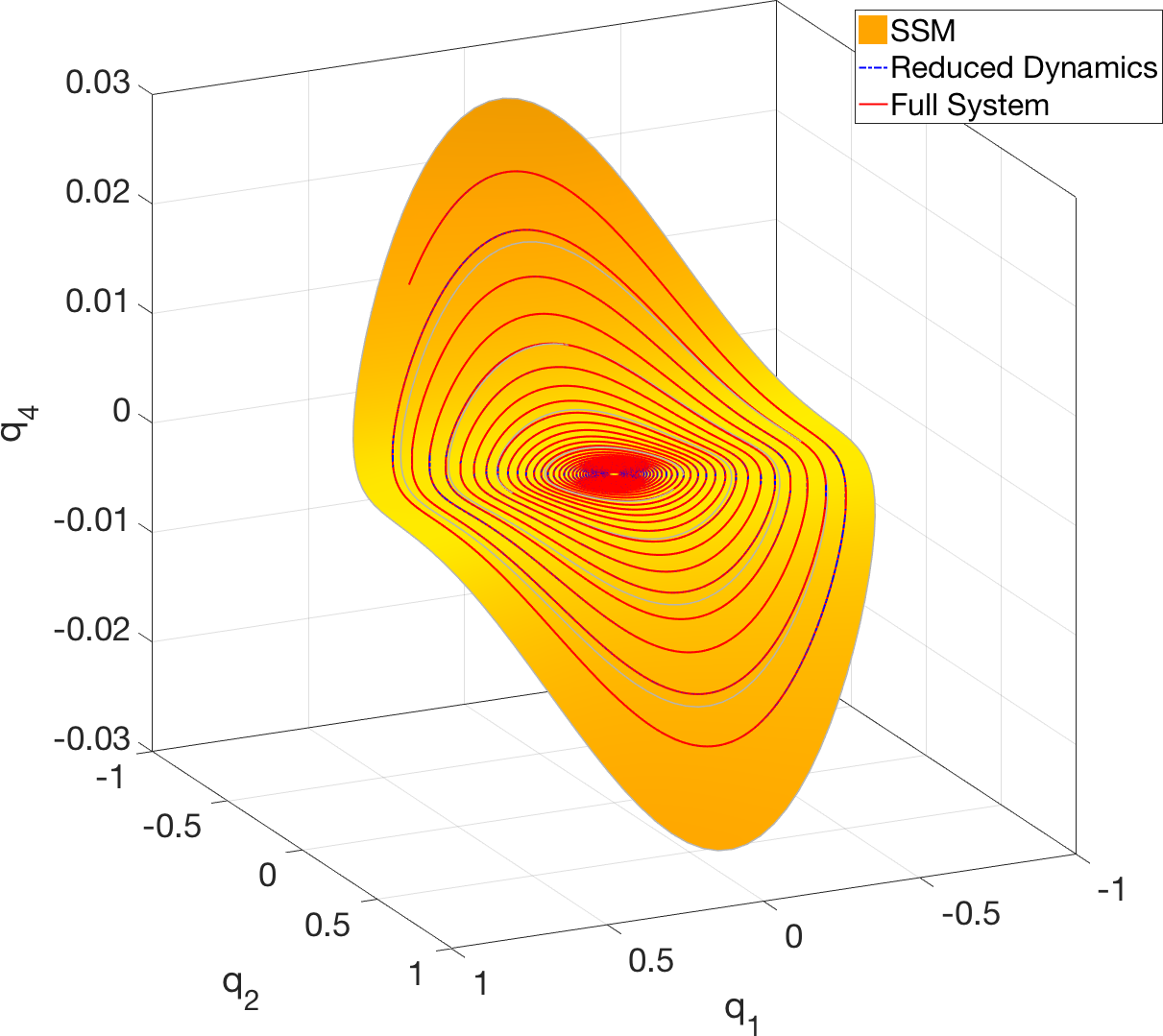}}
\end{centering}
\vfill{}
\begin{centering}
\subfloat[\label{subfig:SSM_2_q3}]{
\includegraphics[scale=0.17]{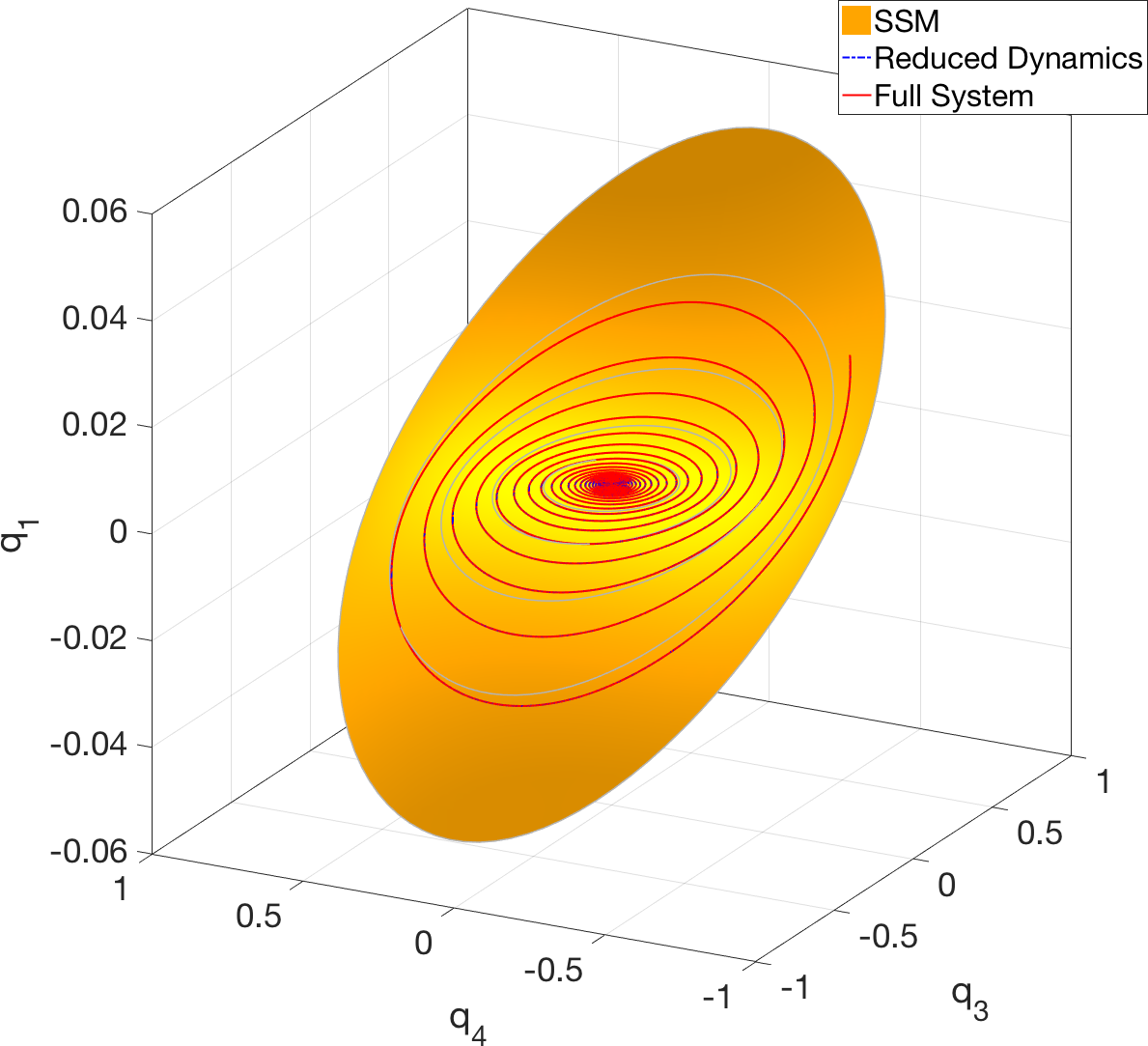}}
\subfloat[\label{subfig:SSM_2_q4}]{
\includegraphics[scale=0.17]{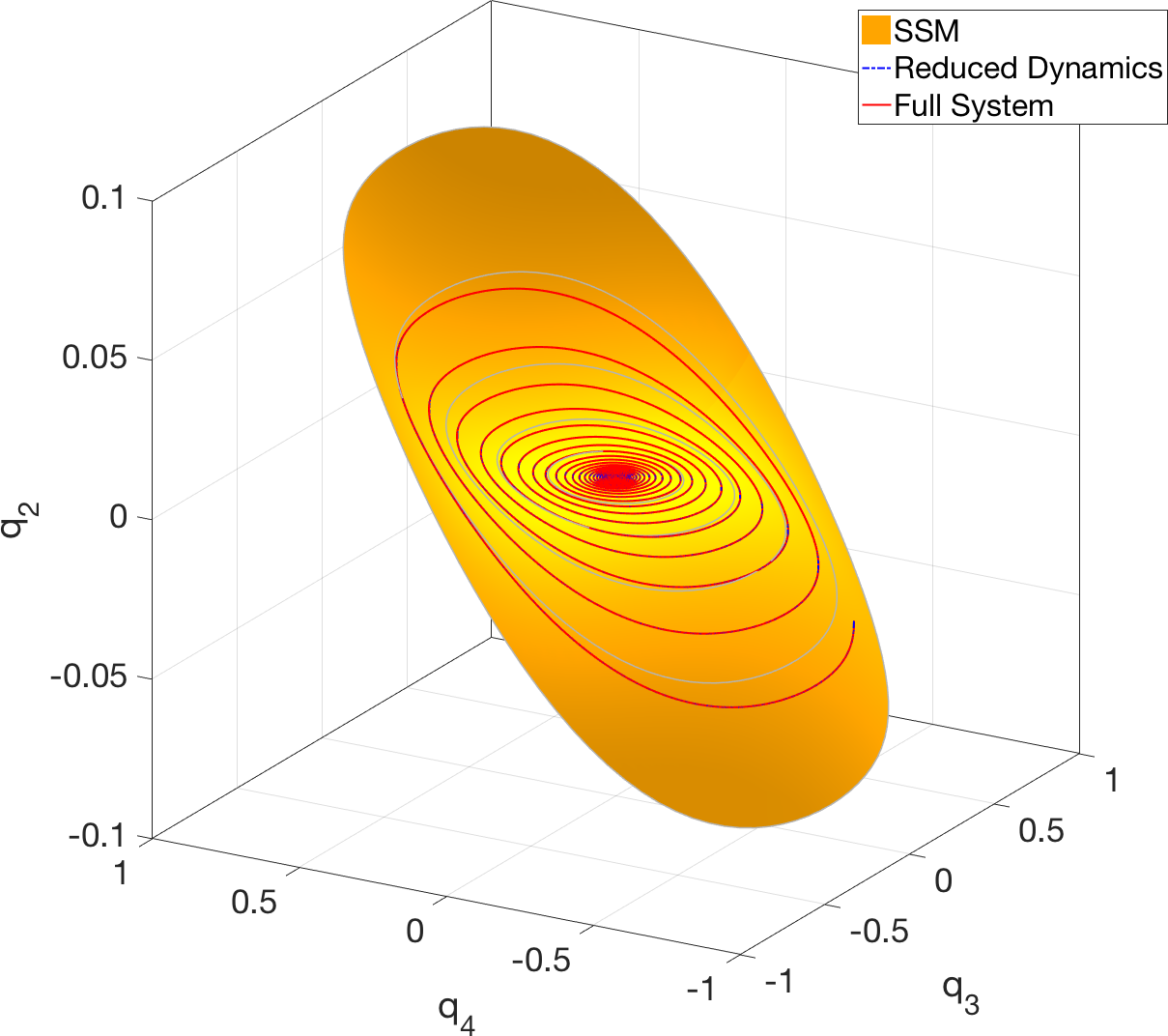}}
\end{centering}
\caption{Lower-dimensional projections of the full phase space, showing the
$15^{\text{th}}$ order approximations of $\vec{W(\mathcal{E}_{1})}$ and
$\vec{W}(\mathcal{E}_{2})$. Figures \ref{subfig:SSM_1_q3} and \ref{subfig:SSM_1_q4}
show the spectral submanifold $\vec{W}(\mathcal{E}_{1})$ tangent to $\mathcal{E}_{1}$.
Figures \ref{subfig:SSM_2_q3} and \ref{subfig:SSM_2_q4} show the
spectral submanifold $\vec{W}(\mathcal{E}_{2})$ tangent to $\mathcal{E}_{2}$.
The dashed curves indicate different projections of a trajectory of
the reduced system $\vec{R}(\vec{z})$, starting from the initial position $\rho=0.35$
and $\theta=1$. The solid curves represent trajectories of the full
system for the same initial position. The solid gray curves represent contour lines of equal 
parameterized distance $\rho$.  \label{fig:SSM_fig_internal}}
\end{figure}

We intend to remove the near-inner resonant terms $z_{1}^{2}z_{2}$,
$z_{1}z_{2}^{2}$, $z_{1}^{3}z_{2}^{2}$, $z_{1}^{2}z_{2}^{3}$, $z_{1}^{4}z_{2}^{3}$,
$z_{1}^{3}z_{2}^{4}$, $z_{1}^{5}z_{2}^{4}$, $z_{1}^{4}z_{2}^{5}$,
$z_{1}^{6}z_{2}^{5}$, $z_{1}^{5}z_{2}^{6}$, $z_{1}^{7}z_{2}^{6}$,
$z_{1}^{6}z_{2}^{7}$, $z_{1}^{8}z_{2}^{7}$ and $z_{1}^{7}z_{2}^{8}$
in the expressions of $\vec{W}(\mathcal{E}_{1})$ and $\vec{W}(\mathcal{E}_{2})$
and add them to the polynomial expressions for the reduced dynamics
on the spectral submanifolds. Due to the choice of nonlinearities,
all the coefficients of $\vec{W}(\mathcal{E}_{1})$ and $\vec{W}(\mathcal{E}_{2})$
corresponding to even powers in $\left|\vec{z}\right|$ are zero. Solving
the partitioned Sylvester equations (\ref{eq:cohomological_vec_1})
and (\ref{eq:cohomological_vec_2}) for $\vec{W}_{i}^{\mathcal{E}}$ and
$\vec{W}_{i}^{\mathcal{C}}$ for orders $i=2,\ldots,15$, we obtain the
lower-dimensional projections of the full phase space for the $15^{\text{th}}$
order approximations of $\vec{W}(\mathcal{E}_{1})$ and $\vec{W}(\mathcal{E}_{2})$,
shown in figure \ref{fig:SSM_fig_internal}. The images are directly
obtained from SSMtool, which detects resonant terms and adds them
to the reduced dynamics $\vec{R}(\vec{z})$ when solving the partitioned Sylvester
equations (\ref{eq:cohomological_vec_1}) and (\ref{eq:cohomological_vec_2}). 

\subsubsection{Reduced Dynamics}

\begin{figure}[H]
\begin{centering}
\subfloat[\label{fig:backbone_curve_1}]{\begin{centering}
\includegraphics[scale=0.21]{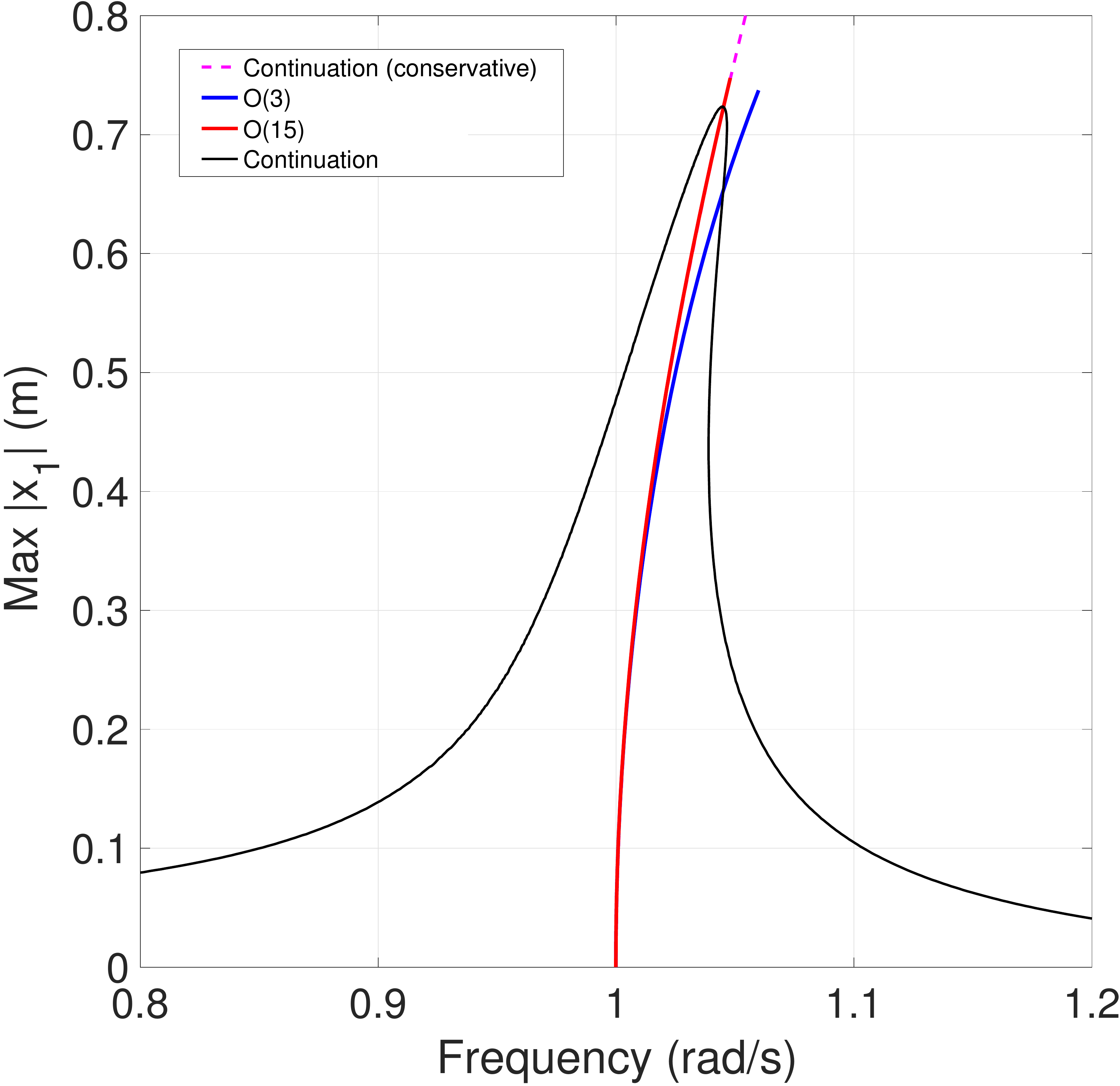}
\par\end{centering}
}\hfill{}\subfloat[\label{fig:backbone_curve_2}]{\begin{centering}
\includegraphics[scale=0.21]{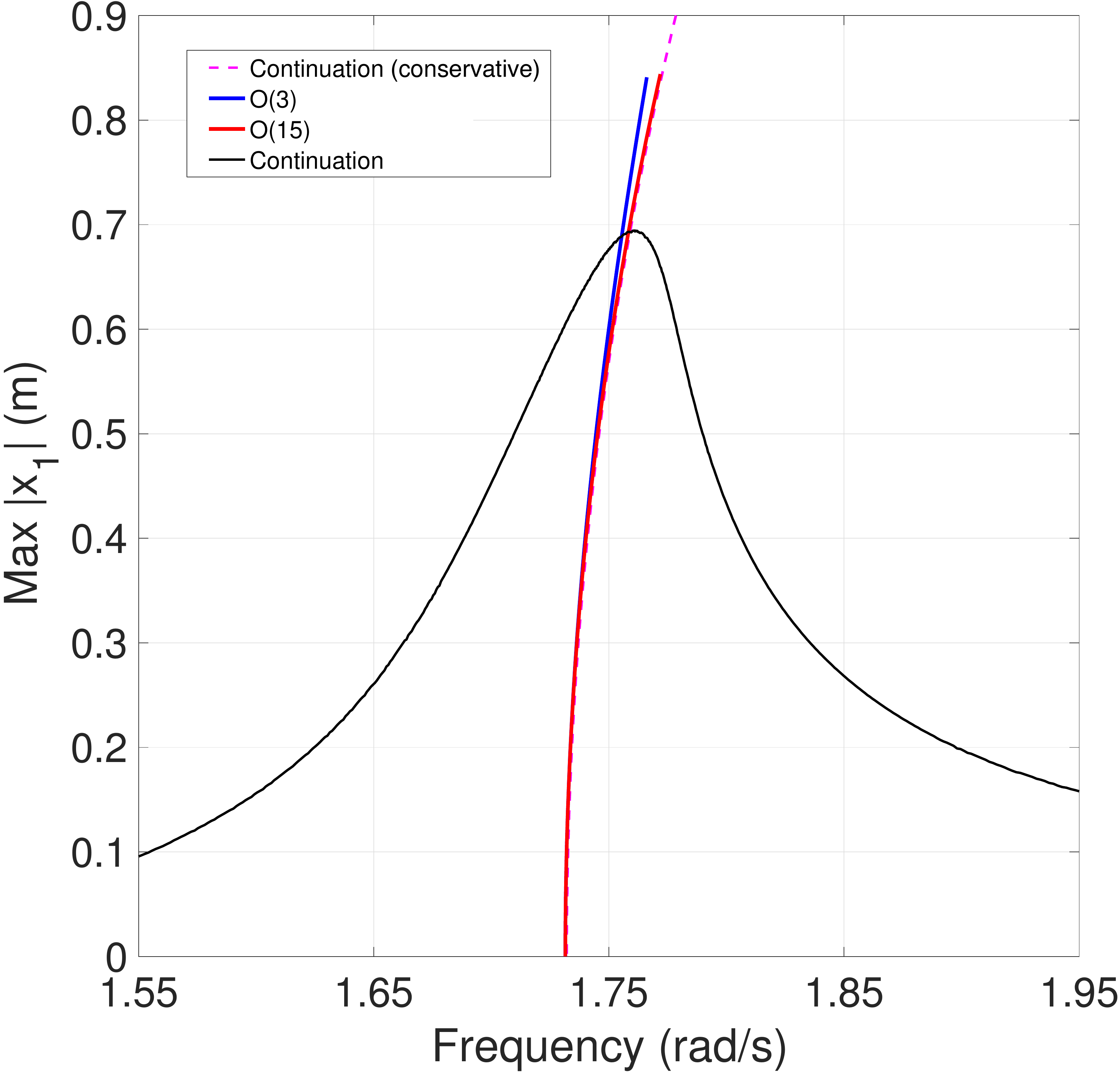}
\par\end{centering}
}
\par\end{centering}
\caption{Backbone curves and periodically forced responses for different amplitudes
of the mechanical system (\ref{eq:sys_mod_sp}). Figure \ref{fig:backbone_curve_1}
shows the $\mathcal{O}(15)$ (red) and $\mathcal{O}(3)$ (blue) approximations
of the backbone curves for the in-phase mode of the system. Figure
\ref{fig:backbone_curve_2} shows the $\mathcal{O}(15)$ (red) and
$\mathcal{O}(3)$ (blue) approximations of the backbone curves for
the out-of-phase mode of the system. Black lines mark amplitudes of
periodic orbits of the periodically forced system for different forcing
amplitudes for varying forcing frequency. The dashed lines (magenta)
represent the backbone curves extracted from COCO in the conservative
limit of the mechanical system, without any forcing. \label{fig:Backbone-curves}}

\end{figure}

The near-inner resonances within the spectral subspaces $\mathcal{E}_{1}$
and $\mathcal{E}_{2}$ introduce nonlinear terms in the reduced dynamics
on the spectral submanifolds. The reduced dynamics on $\vec{W}(\mathcal{E}_{1})$
and $\vec{W}(\mathcal{E}_{2})$ is of the general form (\ref{eq:red_dyn_mixed}).
After transforming to polar coordinates, we obtain the following reduced
equations for the in-phase mode of the system from SSMtool:

\begin{align}
\dot{\rho}= & -0.015\rho-0.00079121\rho^{5}-0.0012708\rho^{7}\label{eq:instant_amp_mode_1}\\
 & +0.0090446\rho^{9}-0.03569\rho^{11}+0.12918\rho^{13}-0.45878\rho^{15},\nonumber \\
\omega= & 0.99989+0.37504\rho^{2}-0.60592\rho^{4}+1.1713\rho^{6}\label{eq:instant_freq_mode_1}\\
 & -2.5137\rho^{8}+5.7885\rho^{10}-14.01\rho^{12}+35.159\rho^{14}.\nonumber 
\end{align}
The reduced dynamics for the out-of-phase mode of the system is obtained
from SSMtool as

\begin{align}
\dot{\rho}= & -0.045\rho+0.016267\rho^{5}+0.02614\rho^{7}\\
 & +0.015714\rho^{9}-0.012768\rho^{11}-0.03437\rho^{13}-0.0308\rho^{15},\nonumber \\
\omega= & 1.7315+0.21658\rho^{2}+0.19904\rho^{4}+0.14858\rho^{6}\label{eq:instant_freq_mode_2}\\
 & +0.072849\rho^{8}+0.017657\rho^{10}+0.004087\rho^{12}-0.011824\rho^{14},\nonumber 
\end{align}
where we set the order of computations to $\mathcal{O}(15)$. Both
instantaneous frequencies (\ref{eq:instant_freq_mode_1}) and (\ref{eq:instant_freq_mode_2})
depend on $\rho$ only. The two red curves in figure \ref{fig:backbone_curve_1}
and figure \ref{fig:backbone_curve_2} represent the $\mathcal{O}(15)$
backbone curves for the in-phase and out-of-phase mode of the mechanical
system, whereas the blue curves display the $\mathcal{O}(3)$ approximations
of the backbone curves. 

We used the numerical continuation software COCO \cite{Dankowicz2013}
to find periodic orbits of the periodically forced system for a fixed
forcing amplitude while varying the forcing frequency. We have extensively optimized the continuation parameters to ensure accurate but fast results. In the current work, COCO is used as an off-the-shelf open-source benchmark to which we compare SSMtool as a stand-alone package. The promising techniques of Blanc et al. \citep{Blanc2013} and Renson et al. \citep{Renson2014} would be expected to perform better than COCO in these computations, but have no available open-source implementations at this point.   

The equations of motion of the forced system are 
\begin{equation}
\left[\begin{array}{cc}
m & 0\\
0 & m
\end{array}\right]\left[\begin{array}{c}
\ddot{x}_{1}\\
\ddot{x}_{2}
\end{array}\right]+\left[\begin{array}{cc}
2c & -c\\
-c & 2c
\end{array}\right]\left[\begin{array}{c}
\dot{x}_{1}\\
\dot{x}_{2}
\end{array}\right]+\left[\begin{array}{cc}
2k & -k\\
-k & 2k
\end{array}\right]\left[\begin{array}{c}
x_{1}\\
x_{2}
\end{array}\right]+\left[\begin{array}{c}
\kappa x_{1}^{3}\\
0
\end{array}\right]=\left[\begin{array}{c}
A\text{ cos }\omega t\\
0
\end{array}\right],
\end{equation}
where we use the same parameter values (\ref{eq:par-setting}) and
introduced a forcing term with amplitude $A$ and forcing frequency
$\omega$. The resulting periodic response amplitudes are shown in
figure \ref{fig:backbone_curve_1} and figure \ref{fig:backbone_curve_2}
in black for a forcing amplitude of $A=0.05\text{ N}$ and $A=0.2\text{ N}$,
respectively. As shown, the $\mathcal{O}(15)$ approximations for
both backbone curves fit the forced peak responses well. The computational
time for the continuation curve in figure \ref{fig:backbone_curve_1}
takes 8 minutes and 6 seconds on a Mac Pro $2\times3.06$ GHz 6-Core
Intel Xeon, which technically corresponds to a single point on the
backbone curve. The computational time for the backbone curve, extracted
from the $15^{\text{th}}$ order approximation of $\vec{W}(\mathcal{E}_{1})$,
is approximately 3 minutes, resulting in a parameterized curve that
can be subsequently evaluated at any required frequency. The computational
time for the continuation curve in figure \ref{fig:backbone_curve_2}
takes a total of 11 minutes and 16 seconds, whereas the computational
time for the backbone curve, extracted from the $15^{\text{th}}$
order approximation of $\vec{W}(\mathcal{E}_{2})$, is also approximately
3 minutes. 

\subsection{The modified Shaw\textendash Pierre example: Outer resonances \label{subsec:The-modified-Shaw=002013Pierre-external}}

We now consider the same analytic example as in section \ref{subsec:The-modified-Shaw=002013Pierre-internal},
but with the linear springs $k_{1}$, $k_{2}$, $k_{3}$ and the damping
$c$ tuned such that there are near-outer resonances and no near-inner
resonances. Here we only compute $\vec{W}(\mathcal{E}_{1})$, the slow SSM
arising from the slow complex pair of eigenvalues. By definition,
it is impossible to obtain an outer resonance for the spectral subspace
corresponding to the remaining fast complex pair of eigenvalues.

The current example is taken from work of Cirillo et al. \cite{cirillo2016spectral},
where a global parameterization method is proposed for the computation
of invariant manifolds in a domain where the stringent non-resonance
conditions of analytic linearization hold. 

For the parameter values
\begin{equation}
k_{1}=k_{3}=1\text{ N }\text{m}^{-1},\quad k_{2}=4.005 \text{ N }\text{m}^{-1},\quad c=0.4 \text{ N s }\text{m}^{-1},\quad\kappa=0.5\text{ N }\text{m}^{-3},\quad\delta=0.05,
\end{equation}
the system has two complex conjugate pairs of eigenvalues with

\begin{align*}
\lambda_{1} & =-0.2+0.9798\mathrm{i},\\
\lambda_{2} & =-0.6+2.9411\mathrm{i},
\end{align*}
and their conjugates. We construct the two-dimensional spectral subspace
$\mathcal{E}_{1}$ corresponding to the first conjugate pair of slow
eigenvalues $\lambda_{1}$ and $\bar{\lambda}_{1}$, whose inner and
outer spectral quotients are

\begin{align}
\sigma_{\text{in}}(\mathcal{E}_{1}) & =\text{Int}\left[\frac{\text{Re}\lambda_{1}}{\text{Re}\lambda_{1}}\right]=1.\\
\sigma_{\text{out}}(\mathcal{E}_{1}) & =\text{Int}\left[\frac{\text{Re}\bar{\lambda}_{2}}{\text{Re}\lambda_{1}}\right]=3,
\end{align}
The exact inner and outer non-resonance conditions, (\ref{eq:ext_res})
and (\ref{eq:int_res}) respectively, again are satisfied for the
spectral subspace $\mathcal{E}_{1}$, i.e., there exists a two-dimensional
analytic SSMs, $\vec{W}(\mathcal{E}_{1})$, that is unique among all $C^{4}$
invariant manifolds tangent to $\mathcal{E}_{1}$. Due to the higher
choice of damping, there are no near-inner resonances and hence the
reduced dynamics on the manifold can be expressed as linear. However,
the SSM is close to having two third-order outer resonances which
in turn leads to the two near-outer resonances shown in table \ref{tab:Near-external-resonances_E1}.
For $k_{2}=4\text{ N }\text{m}^{-1}$, the SSM construction will break down as
$\vec{\Theta}_{3}^{\mathcal{C}}$ becomes singular while equation (\ref{eq:cohomological_vec_2})
has a nonzero right-hand side. 

\begin{figure}[h]
\begin{centering}
\includegraphics[scale=0.3]{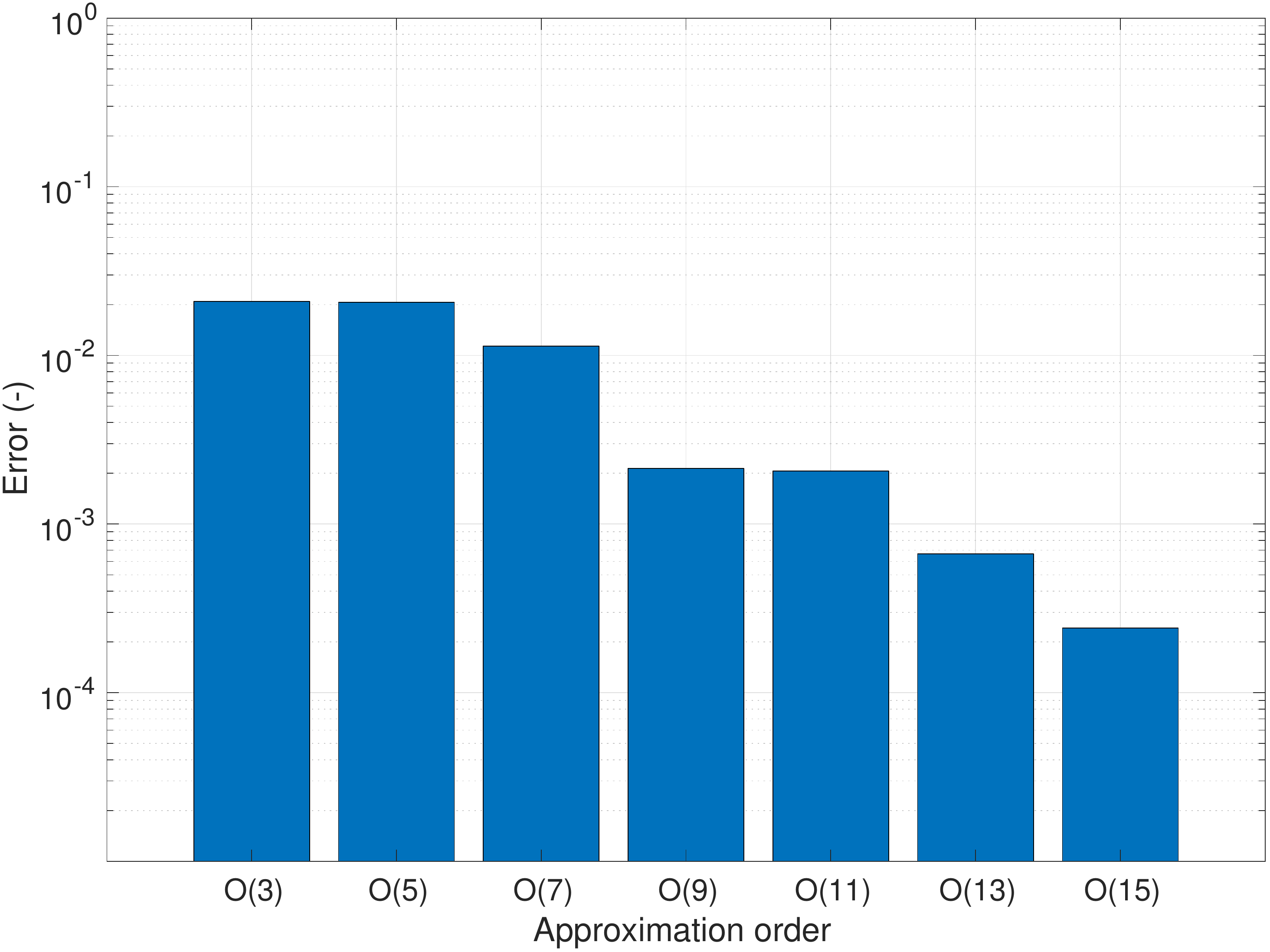}
\par\end{centering}
\begin{centering}
\caption{Invariance error for the $3^{\text{th}}$-$15^{\text{th}}$ order
approximations of $\vec{W}(\mathcal{E}_{1})$ for 50 evenly distributed
initial positions lying on a fixed radius $\rho_{0}=0.28$ and $\theta_{0}\in S^{1}$.
For each trajectory traveling between $\rho_{0}$ and $\rho_{\epsilon}=0.01$,
we identify the maximum error and take the average over all trajectories.
\label{fig:SSM_invar_error_2DOF_ext_res}}
\par\end{centering}
\end{figure}

\begin{table}[H]
\begin{centering}
\begin{tabular}{ccccc}
\toprule 
\multicolumn{5}{c}{$\mathcal{E}_{1}$}\tabularnewline
\midrule 
 & $a$ & $b$ & $\lambda_{l}$ & $I$\tabularnewline
\midrule
\multirow{2}{*}{$\mathcal{O}(\left|\vec{z}\right|^{3})$} & 3 & 0 & $\lambda_{2}$ & 0.000162\tabularnewline
 & 0 & 3 & $\bar{\lambda}_{2}$ & 0.000162\tabularnewline
\bottomrule
\end{tabular}
\par\end{centering}
\caption{Near-outer-resonances for $\mathcal{E}_{1}$ with $\delta=0.05$.\label{tab:Near-external-resonances_E1}}
\end{table}
As has been done in section \ref{subsec:Computing_W(E1)}, we would
like to identify the order to which we have to approximate the SSM
to obtain an accurate reduced order model. Using the invariance measure
defined in equation (\ref{eq:inv_error_measure-1}), we test the invariance
of $\vec{W}(\mathcal{E}_{1})$ for different approximation orders. In figure
\ref{fig:SSM_invar_error_2DOF_ext_res}, we show the invariance error
for seven different approximations of $\vec{W}(\mathcal{E}_{1})$. For a
given fixed radius $\rho_{0}=0.28$ we take $50$ initial points,
each corresponding to an angle $\theta_{0}$, uniformly distributed
in $S^{1}$. 

\begin{figure}[H]
\centering{}\subfloat[\label{fig:ssm_ext_res_c1}]{
\begin{centering}
\includegraphics[scale=0.07]{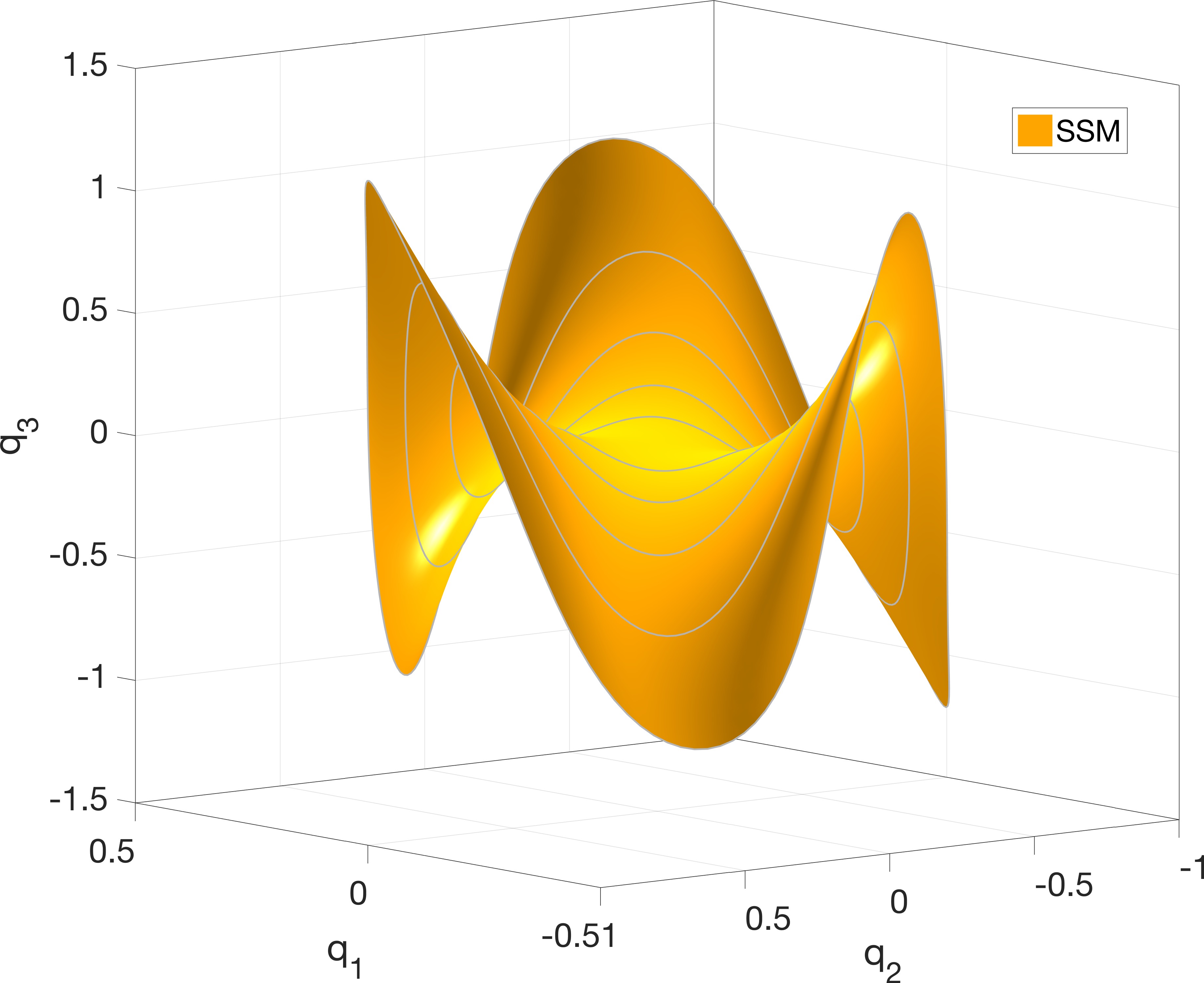}
\par\end{centering}
}\vfill{}
\subfloat[\label{fig:ssm_ext_res_c2}]{\begin{centering}
\includegraphics[scale=0.05]{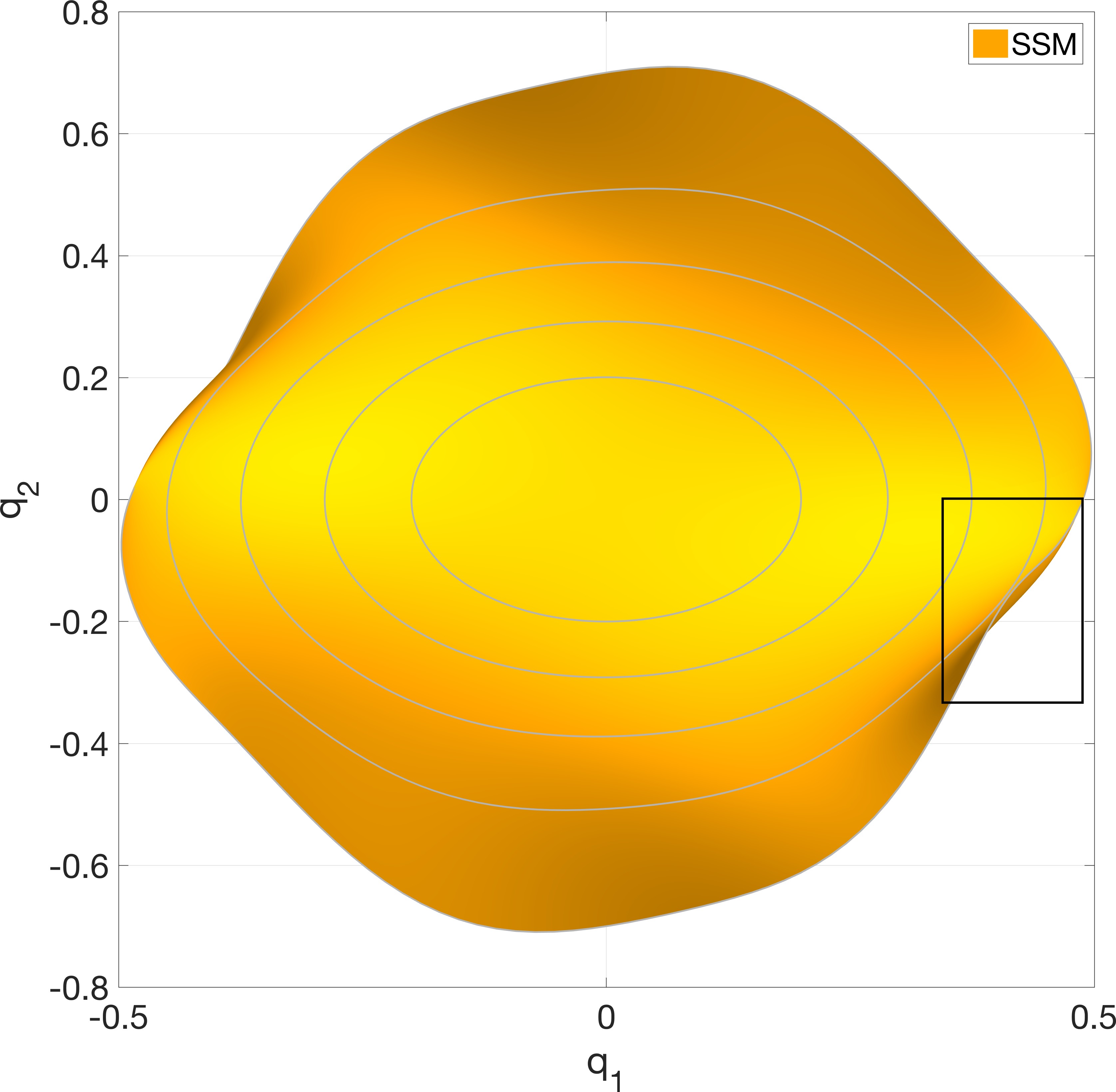}
\par\end{centering}
}\hfill{}\subfloat[\label{fig:ssm_ext_res_c3}]{
\begin{centering}
\includegraphics[scale=0.19]{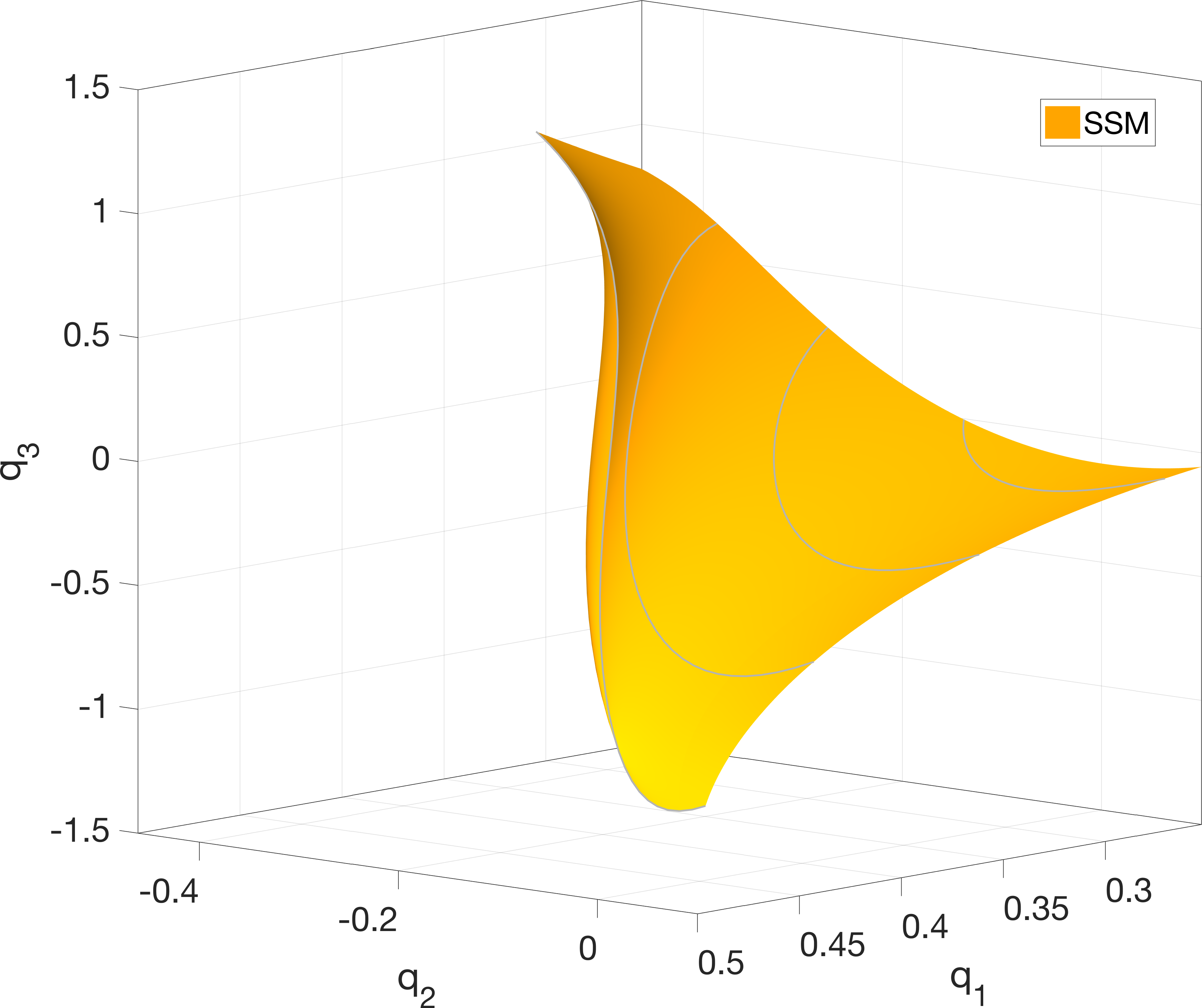}
\par\end{centering}
}\caption{Lower-dimensional projections of the full phase space showing the
$15^{\text{th}}$ order approximation of $\vec{W}(\mathcal{E}_{1})$
in modal coordinates. Figures \ref{fig:ssm_ext_res_c1}, \ref{fig:ssm_ext_res_c2}
and \ref{fig:ssm_ext_res_c3} show the spectral submanifold $\vec{W}(\mathcal{E}_{1})$
tangent to $\mathcal{E}_{1}$, being close to outer resonance. Figure
\ref{fig:ssm_ext_res_c2}, a top-view of $\vec{W}(\mathcal{E}_{1})$,
shows the development of a fold as indicated by the rectangle. Figure
\ref{fig:ssm_ext_res_c3} is a zoomed-in version of the fold. \label{fig:ssm_ext_res_all}}
\end{figure}

As a the SSM is near an outer resonance, a folding of the SSM over
its underlying modal subspace is more likely to occur. Such a folding
is illustrated in figure \ref{fig:ssm_ext_res_all}, showing a lower-dimensional
projection of the full phase space of the $15^{\text{th}}$ order
approximation of $\vec{W}(\mathcal{E}_{1})$. This example brings out the
power of the parameterization method, as constructing the SSM as a
graph over its modal subspace would break down at the point of folding. 

\begin{figure}[H]
\begin{centering}
\subfloat[\label{subfig:SSM_Cir_1_x2}]{\begin{centering}
\includegraphics[scale=0.053]{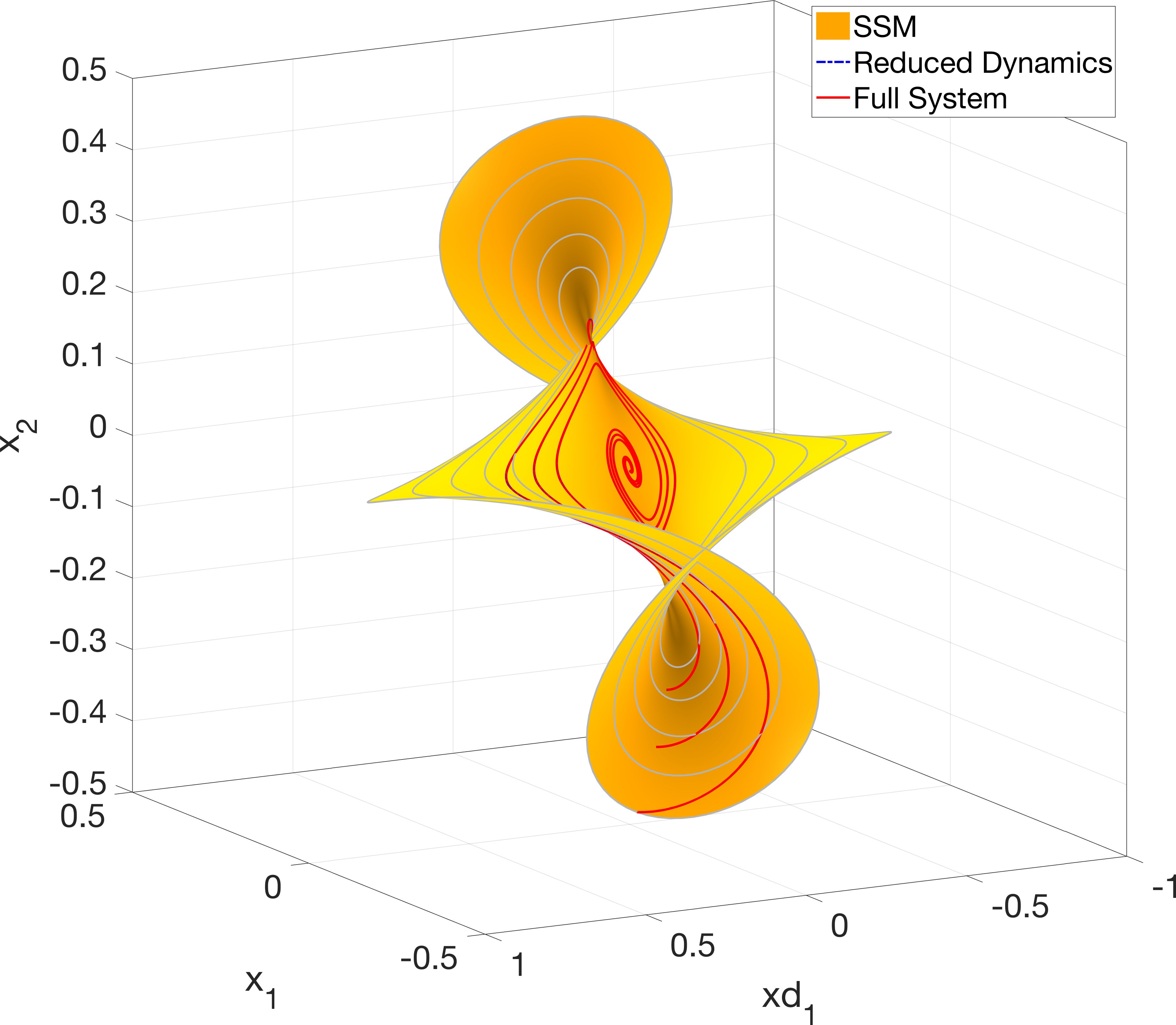}
\par\end{centering}
}\hfill{}\subfloat[\label{subfig:SSM_Cir_2_xd2}]{\begin{centering}
\includegraphics[scale=0.053]{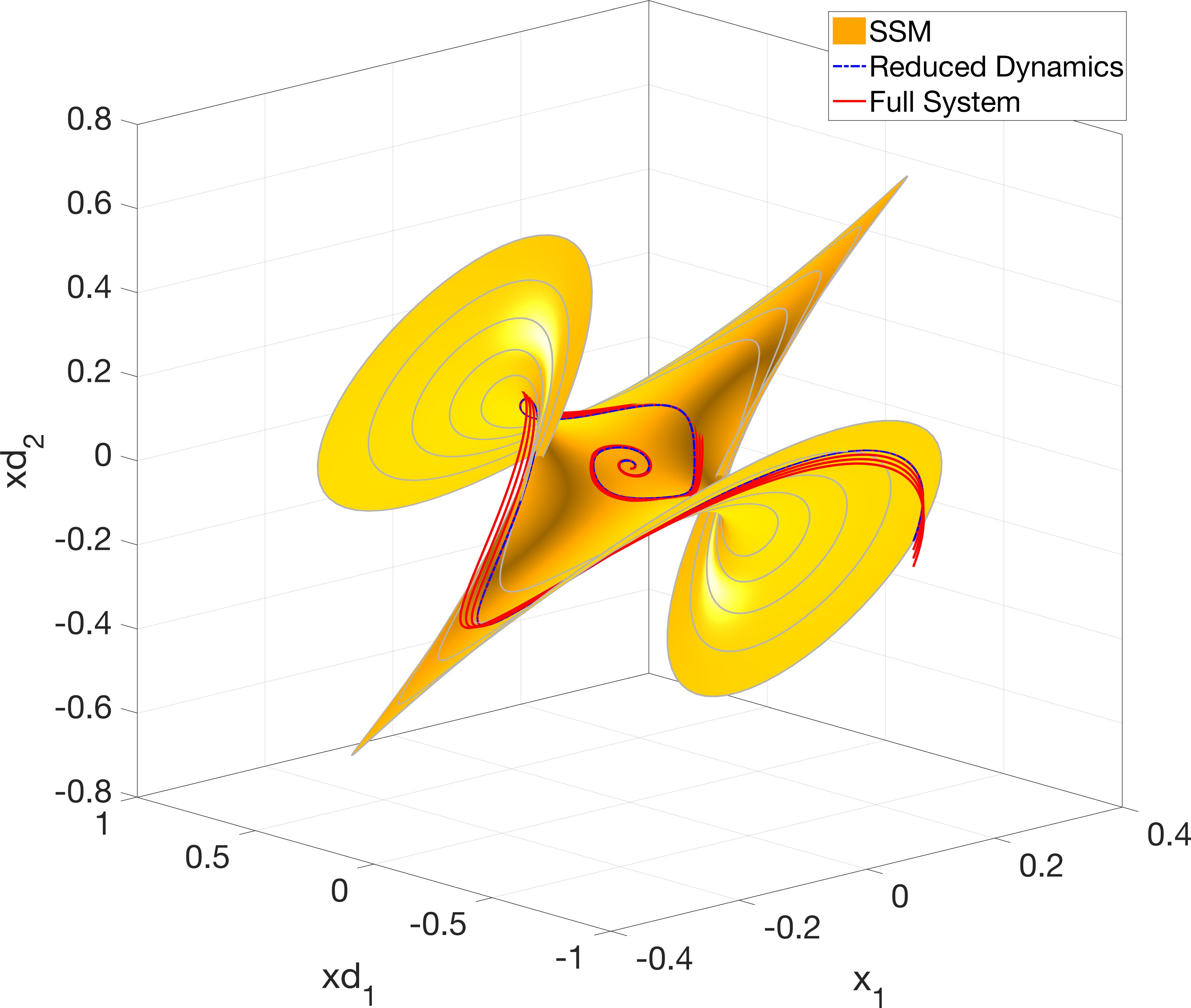}
\par\end{centering}
}
\par\end{centering}
\caption{Lower-dimensional projections of the full phase space showing the
$15^{\text{th}}$ order approximation of $\vec{W}(\mathcal{E}_{1})$,
transformed to physical coordinates. The dashed curves in figure \ref{subfig:SSM_Cir_1_x2}
correspond to trajectories of the reduced system $\vec{R}(\vec{z})$ corresponding
to the initial positions $\rho=\left\{ 0.15,0.13,0.11\right\} $ and
$\theta=3$. The solid curves represent trajectories of the full system
for the same initial positions. In figure \ref{subfig:SSM_Cir_2_xd2},
the dashed curves corresponds to a trajectory of the reduced system
$\vec{R}(\vec{z})$ for the initial position $\rho=0.15$ and $\theta=1$. The
solid lines represent trajectories of the full system having an initial
position off the manifold, showing the convergence towards $\vec{W}(\mathcal{E}_{1})$.
\label{fig:SSM_Cir_fig_external}}
\end{figure}

In figure \ref{fig:SSM_Cir_fig_external}, we show the SSM transformed
to physical coordinates (also an option in SSMtool), where we demonstrate
the invariance of $\vec{W}(\mathcal{E}_{1})$ (figure \ref{subfig:SSM_Cir_1_x2})
and that different trajectories converge towards $\vec{W}(\mathcal{E}_{1})$
(figure \ref{subfig:SSM_Cir_2_xd2}), when starting close to $\vec{W}(\mathcal{E}_{1})$.

\subsection{The discretized nonlinear Timoshenko beam \label{subsec:application-beam}}

In this section, we construct a reduced order model for a discretized
nonlinear Timoshenko beam by computing the reduced dynamics on the
two-dimensional SSM arising from the slowest modal subspace. We will
briefly outline the steps leading to the derivation of the partial
differential equations (PDEs) governing the dynamics of the beam.
Our reasoning largely follows the presentation given by Reddy \cite{Reddy2013}.
The problem considered here is a square 2D beam placed in a cartesian
coordinate system with coordinates $(x,y,z)$ and basis $(\vec{e}_{x},\vec{e}_{y},\vec{e}_{z})$.
Initially, the beam is straight, with its main axis parallel to the
$x$-axis, while its cross section lies in the $y-z$ plane. The relevant
beam parameters are listed in table \ref{tab:system_par_beam}. 

\begin{table}[H]
\begin{centering}
\begin{tabular}{|c|c|}
\hline 
Symbol & Meaning {(}unit{)}\tabularnewline
\hline 
\hline 
$L$ & Length of beam {(}mm{)}\tabularnewline
\hline 
$h$ & Height of beam {(}mm{)}\tabularnewline
\hline 
$b$ & Width of beam {(}mm{)}\tabularnewline
\hline 
$\varrho$ & Density {(}kg mm\textsuperscript{-3}{)}\tabularnewline
\hline 
$E$ & Young's Modulus {(}MPa{)}\tabularnewline
\hline 
$G$ & Shear Modulus {(}MPa{)}\tabularnewline
\hline 
$\eta$ & Axial material damping constant {(}MPa s{)}\tabularnewline
\hline 
$\mu$ & Shear material damping constant {(}MPa s{)}\tabularnewline
\hline 
$\lambda$ & External damping constant {(}MPa s mm\textsuperscript{-2}{)}\tabularnewline
\hline 
$A=bh$ & Cross-section of beam {(}mm\textsuperscript{2}{)}\tabularnewline
\hline 
\end{tabular}
\par\end{centering}
\caption{Notation used in subsequent derivations \label{tab:system_par_beam}}
\end{table}

\begin{figure}[H]
\begin{centering}
\includegraphics[scale=0.7]{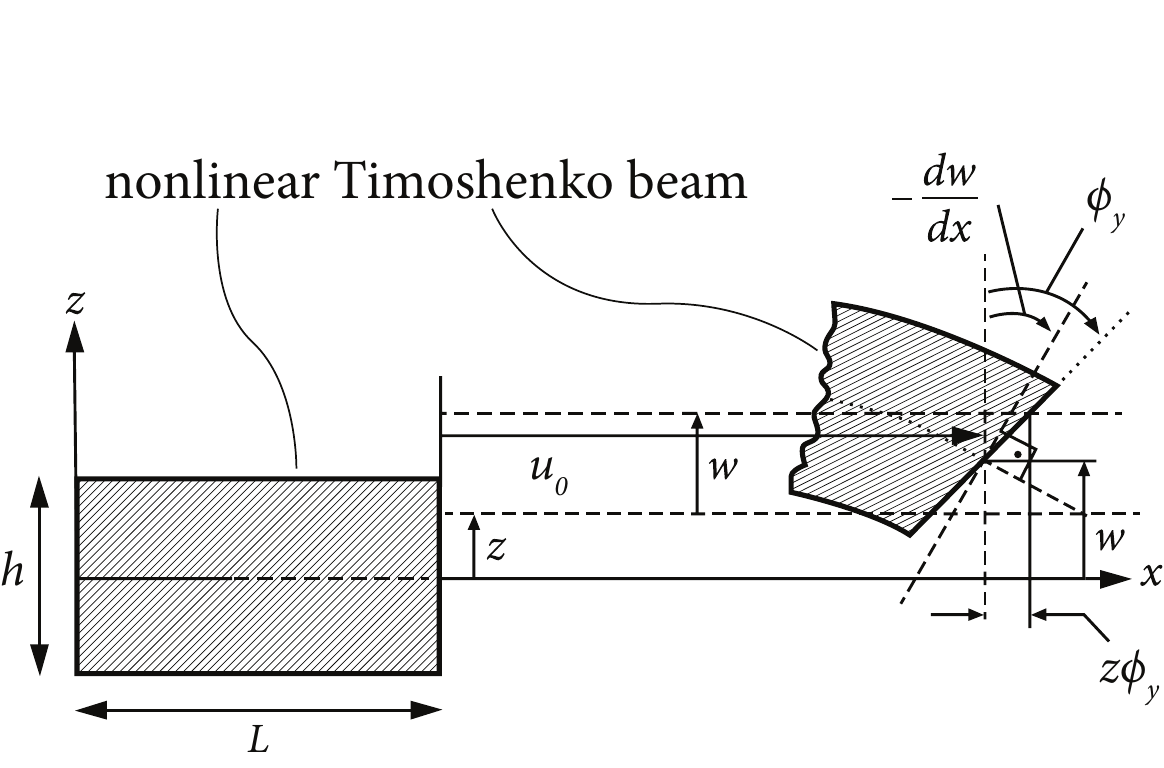}
\par\end{centering}
\caption{Kinematics of the nonlinear Timoshenko beam.\label{fig:kinematics-beam}}
\end{figure}

We call the line that initially coincides with the $x$-axis the beam's
neutral axis. The kinematic assumptions underlying the Timoshenko
beam model can be obtained by relaxing the restrictions of the Bernoulli
hypothesis which is the basis of the more classical and well-known
Euler-Bernoulli beam theory. The Bernoulli hypothesis states (cf.
Reddy \cite{Reddy2013}) that initially straight material lines normal
to the neutral axis remain (a) straight and (b) inextensible after
deformation, and (c) rotate as rigid lines to remain perpendicular
to the beam's neutral axis after deformation. We relax (c) by allowing
for rigid rotation of the cross section about the $y$-axis. These
kinematic assumptions are satisfied by the following displacement
field

\begin{align}
u_{x}(x,y,z) & =u_{0}(x)+z\phi_{y}(x),\\
u_{y}(x,y,z) & =0,\\
u_{z}(x,y,z) & =w(x).
\end{align}
Here $(u_{x},u_{y},u_{z})$ are the components of the displacement
field $\vec{u}(x,y,z)$ for a material point in the $(x,y,z)$ directions,
respectively. The functions $u_{0}(x)$ and $w(x)$ represent the
displacements of a material point with initial coordinates on the
beam's neutral axis, given by $z=0$. The rotation of a normal section
about the y-axis is denoted by $\phi_{y}(x)$. We illustrate the kinematics
in figure \ref{fig:kinematics-beam}.

Following Reddy \cite{Reddy2013}, we neglect all $(u_{0,x})^{2}$
terms in the Green-Lagrange strain tensor, where we use the shorthand
notation $\left(\cdot\right)_{,x}=\partial_{x}(\cdot)$. This approximation
accounts for geometric nonlinearities due to moderately large rotations
while assuming small membrane strains. The relevant non-zero components
of the simplified Green-Lagrange strain tensor $\varepsilon_{ij}$
take the form

\begin{align}
\varepsilon_{xx} & =\varepsilon_{xx}^{0}+z\varepsilon_{xx}^{1},\\
\gamma_{xz} & =2\varepsilon_{xy}=\gamma_{xz}^{0}+z\gamma_{xz}^{1},
\end{align}
with

\begin{align}
\varepsilon_{xx}^{0} & =\partial_{x}u_{0}+\frac{1}{2}(\partial_{x}w)^{2},\ \varepsilon_{xx}^{1}=\partial_{x}\phi_{y},\\
\gamma_{xz}^{0} & =\phi_{y}+\partial_{x}w+\phi_{y}\partial_{x}u_{0},\ \gamma_{xz}^{1}=\phi_{y}\partial_{x}\phi_{y}.
\end{align}
We assume a linear viscoelastic constitutive relation between the
stresses and strains of the following form

\begin{align}
\sigma_{xx} & =E\varepsilon_{xx}+\eta\dot{\varepsilon}_{xx},\label{eq:cauchy_xx}\\
\sigma_{xz} & =G\gamma_{xz}+\mu\dot{\gamma}_{xz}.\label{eq:cauchy_xz}
\end{align}
Here $\sigma_{xx}$ and $\sigma_{xz}$ are the components of the Cauchy
stress tensor $\boldsymbol\sigma$ (see, e.g., Lai et al. \cite{Lai2009}). The
relations given in (\ref{eq:cauchy_xx}) and (\ref{eq:cauchy_xz})
are a special case of the more general linear viscoelastic material
models that can be found, for example, in Skrzypek and Ganczarski
\cite{Skrzypek2015} and is also used, e.g., by Lesieutre and Kauffman
\cite{Lesieutre2013}. We explain the derivation and discretization
of the equations of motion of the beam in \ref{sec:Equations-of-motion-beam}. 

In the following computations, we will consider a beam that is clamped
on one end and free on the other, which means that on the clamped
end all displacements $(u_{0}$, $w$, $\phi_{y})$ are zero, while
on the free end no restrictions are placed on the displacements. After
implementation of the essential boundary conditions, the number of
degrees of freedom $n$ of our system is given by 

\begin{equation}
n=5m+1,
\end{equation}
where $m$ is the number of finite beam elements used in the discretization.
Additionally, we set the external damping parameter $c$, discussed
in \ref{sec:Equations-of-motion-beam}, to zero and therefore
the damping of our beam only enters through the viscoelastic constitutive
relation.

We continue by constructing the slowest single-mode SSM for a specific
beam, for which we will use the SSMtool to reduce the beam dynamics
to a two-dimensional system of ordinary differential equations. The
chosen geometric and material parameters are listed in Table \ref{tab:system_par_beam_ex}. 

\begin{table}[H]
\begin{centering}
\begin{tabular}{|c|c|}
\hline 
Parameter & Value\tabularnewline
\hline 
\hline 
$L$ & $\unit[1000]{mm}$\tabularnewline
\hline 
$h$ & $\unit[100]{mm}$\tabularnewline
\hline 
$b$ & $\unit[100]{mm}$\tabularnewline
\hline 
$\varrho$ & $\unit[7850\cdot10^{-9}]{kg\text{ }mm^{-3}}$\tabularnewline
\hline 
$E$ & $\unit[90]{GPa}$\tabularnewline
\hline 
$G$ & $\unit[34.6]{GPa}$\tabularnewline
\hline 
$\eta$ & $\unit[33.6]{MPa\text{ }s} $\tabularnewline
\hline 
$\mu$ & $\unit[20.9]{MPa\text{ }s}$\tabularnewline
\hline 
\end{tabular}
\par\end{centering}
\caption{Geometric and material parameters. \label{tab:system_par_beam_ex}}
\end{table}
We simulate the beam with three elements, resulting in a 32-dimensional
phase space. For the chosen parameter values, the eigenvalues corresponding
to the slowest modal subspace $\mathcal{E}$ are $\lambda_{1,2}=-0.02286\pm11.03\mathrm{i}$.
In terms of its exponential decay rate, the eigenspace $\mathcal{E}$
is about 50 times slower compared to the second slowest eigenspace
of the system. This spectral ratio indicates that trajectories of
the system transverse to the slow SSM die out fast, making the slowest
SSM an excellent choice for model reduction, because it will contain
trajectories that remain active for the longest time. 

\begin{figure}[H]
\begin{centering}
\includegraphics[scale=0.25]{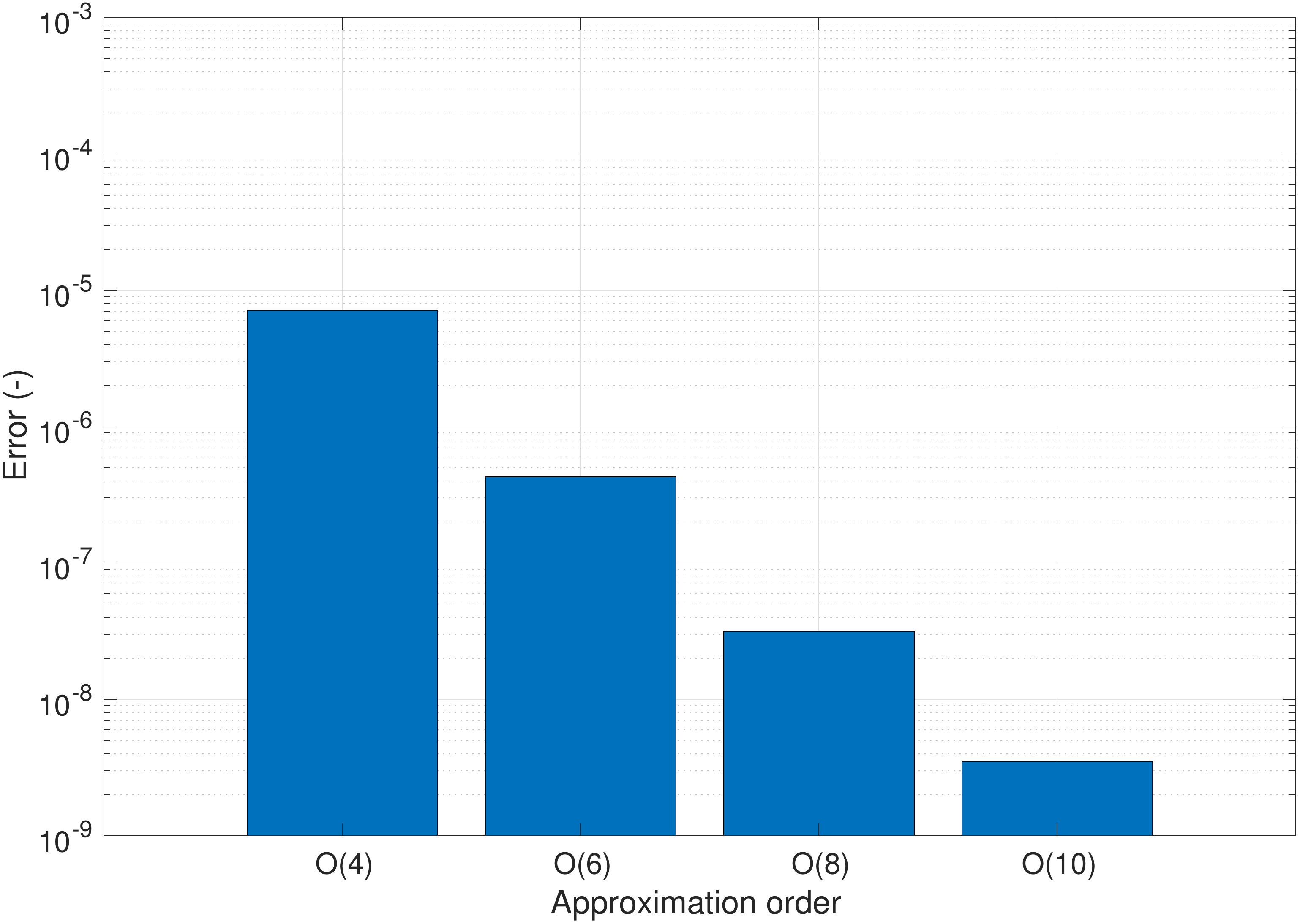}
\par\end{centering}
\caption{Invariance error for the $4^{\text{th}}$, $6^{\text{th}}$, $8^{\text{th}}$
and $10^{\text{th}}$ order approximations of $\vec{W}(\mathcal{E})$ for
50 evenly distributed initial positions lying on a fixed radius $\rho_{0}=1.5$
and $\theta_{0}\in S^{1}$. The chosen value of $\rho_{0}=1.5$ corresponds
to a maximum physical vertical displacement of the endpoint of the
beam of $160\text{ mm}$. For each trajectory traveling between $\rho_{0}$
and $\rho_{\epsilon}=0.2$, we identify the maximum error and take
the average over all trajectories. \label{fig:SSM_invar_error}}
\end{figure}

In figure \ref{fig:SSM_invar_error}, we show the invariance error
(\ref{eq:inv_error_measure-1}) for four different orders of approximations
of $\vec{W}(\mathcal{E})$. 

\begin{figure}[H]
\begin{centering}
\subfloat[\label{fig:beam_ssm_1}]{\centering{}\includegraphics[scale=0.05]{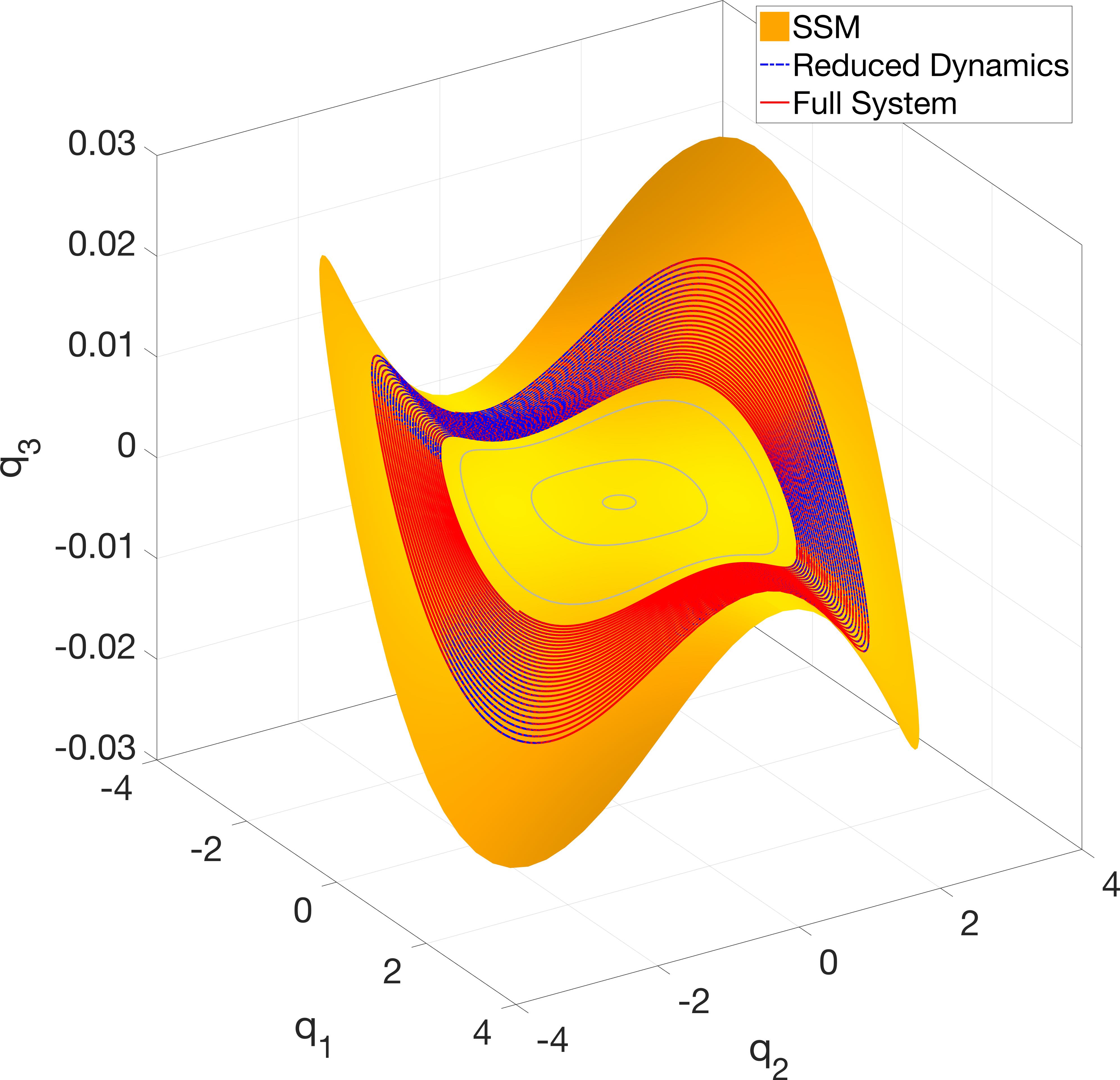}}\hfill{}\subfloat[\label{fig:beam_ssm_2}]{\begin{centering}
\includegraphics[scale=0.05]{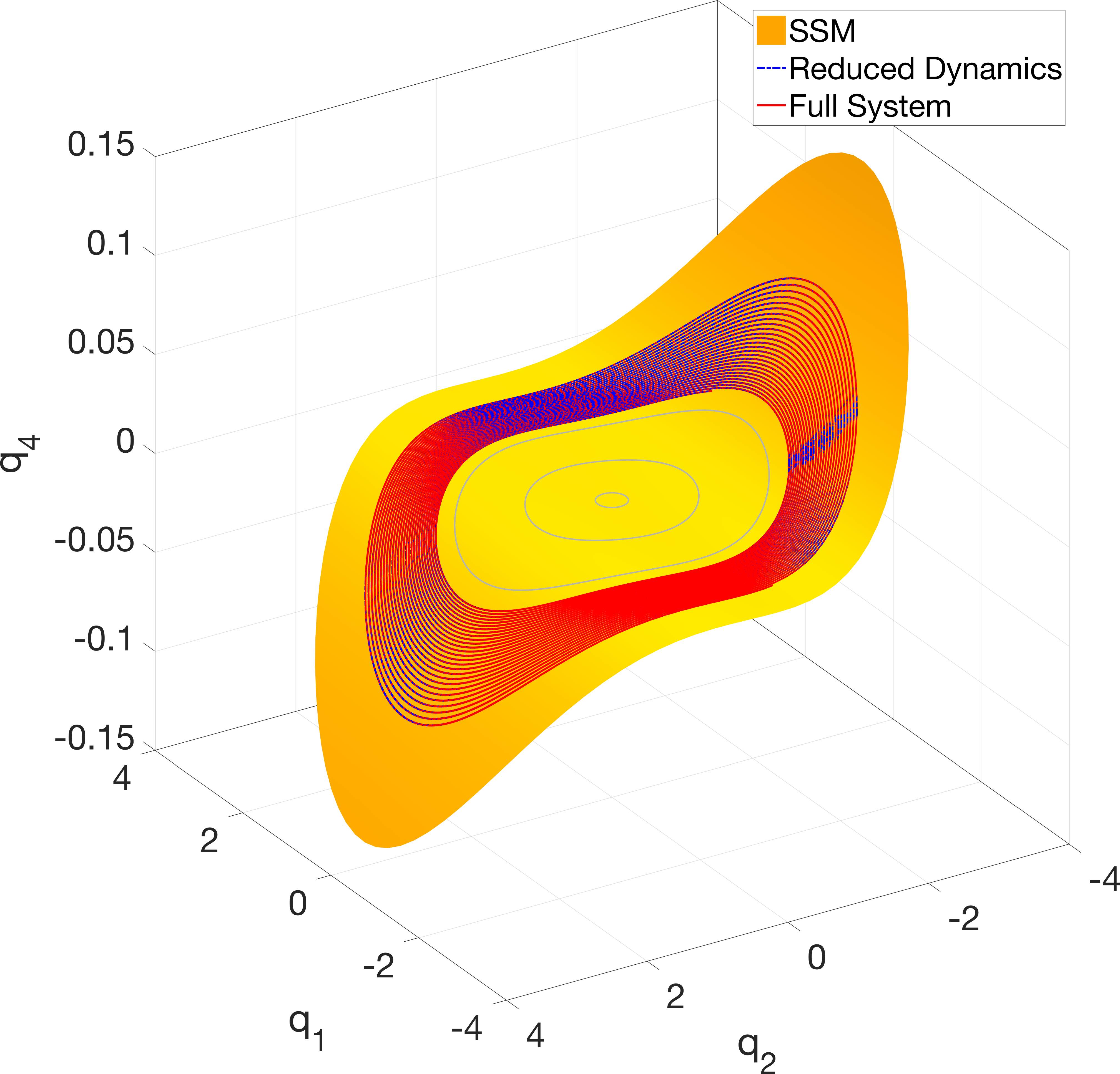}
\par\end{centering}
}
\par\end{centering}
\caption{Lower-dimensional projections of the 32-dimensional phase space, of
the discretized nonlinear Timoshenko beam, showing the $10^{\text{th}}$
order approximation of $\vec{W}(\mathcal{E})$. Figures \ref{fig:beam_ssm_1},
\ref{fig:beam_ssm_2} show the spectral submanifold $\vec{W}(\mathcal{E})$
tangent to slowest modal subspace $\mathcal{E}$. The dashed curve
corresponds to a trajectory of the reduced two-dimensional system
$\vec{R}(z)$ corresponding to the initial positions $r=1.5$ and $\theta=3$.
The solid curve represents a trajectory of the full system for the
same initial position, with $t_{\text{end}}=15\text{ s}$. \label{fig:SSM_beam}}
\end{figure}

We observe that the $4^{\text{th}}$ order approximation of $\vec{W}(\mathcal{E})$
is already accurate up to the chosen value of $\rho_{0}=1.5$, which
corresponds to a maximum physical vertical displacement of the endpoint
of the beam of $160\text{ mm}$. By increasing the order of approximation
to $10$, the invariance error is reduced further by approximately
three orders of magnitude. Figure \ref{fig:SSM_beam} displays two
lower-dimensional projections of the 32-dimensional phase space, showing
the $10^{\text{th}}$ order approximation of $\vec{W}(\mathcal{E})$, with
the the modal coordinates $q_{3}$ and $q_{4}$ plotted over the coordinates
$q_{1}$ and $q_{2}$. 

\begin{figure}[H]
\begin{centering}
\includegraphics[scale=0.25]{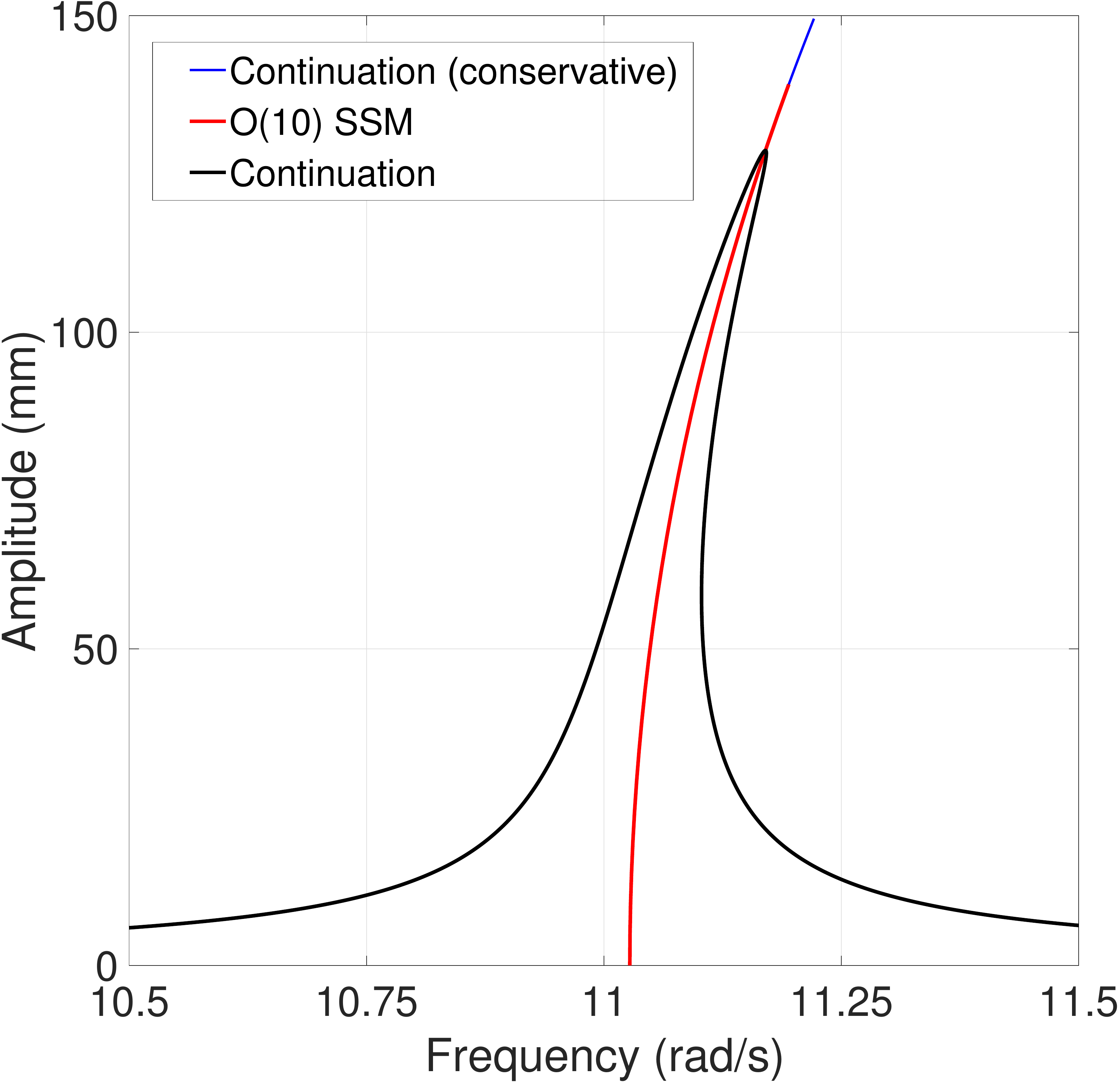}
\par\end{centering}
\caption{Backbone curve and periodically forced responses of the discretized
Timoshenko beam having a 32-dimensional phase space. The $\mathcal{O}(10)$
(red) approximation of the backbone curve is computed up to $\rho=1.3$.
The black line corresponds to periodic orbits of the periodically
forced system for varying forcing frequency. The blue line represents the backbone curve extracted from COCO in the conservative limit of the beam, without any forcing. \label{fig:Backbone-curve-beam}}
\end{figure}

As $\lambda_{1}$ and $\bar{\lambda}_{1}$ have small negative real
parts, the near-inner-resonances conditions related to $\mathcal{O}(\left|\vec{z}\right|^{i})$
with $i=3,5,7,9$ are satisfied within the spectral subspaces $\mathcal{E}$.
This in turn leads to the following expressions for reduced dynamics
on $\vec{W}(\mathcal{E})$, obtained from SSMtool:

\begin{align}
\dot{\rho}= & -0.022856\rho-0.00017033\rho^{3}-4.9542\cdot10^{-6}\rho^{5}\label{eq:instant_amp_beam}\\
 & +8.5365\cdot10^{-8}\rho^{7}-3.0348\cdot10^{-9}\rho^{9},\nonumber \\
\omega= & 11.027+0.099097\rho^{2}-0.000020843\rho^{4}-2.8625\cdot10^{-6}\rho^{6}\label{eq:instant_freq_beam}\\
 & +1.729\cdot10^{-7}\rho^{8}\nonumber 
\end{align}
Using the definition of the instantaneous amplitude (\ref{eq:instant_amp})
and the corresponding instantaneous frequency (\ref{eq:instant_freq_beam}),
SSMtool computes the parameterized backbone curve $\mathcal{B}$ (shown
in figure \ref{fig:Backbone-curve-beam}) in less than 4 minutes time.
The continuation curve, shown in black, has been computed using the
numerical continuation software COCO \cite{Dankowicz2013}, after
applying a periodic force, $F=A\text{ cos }\omega t$, to the vertical
displacement coordinate $w$ at the free end of the beam with a forcing
amplitude of $A=300\text{ N}$. 

\begin{figure}[H]
\begin{center}
\subfloat[\label{fig:beam_ssm_OT}]{
\begin{centering}
\includegraphics[scale=0.06]{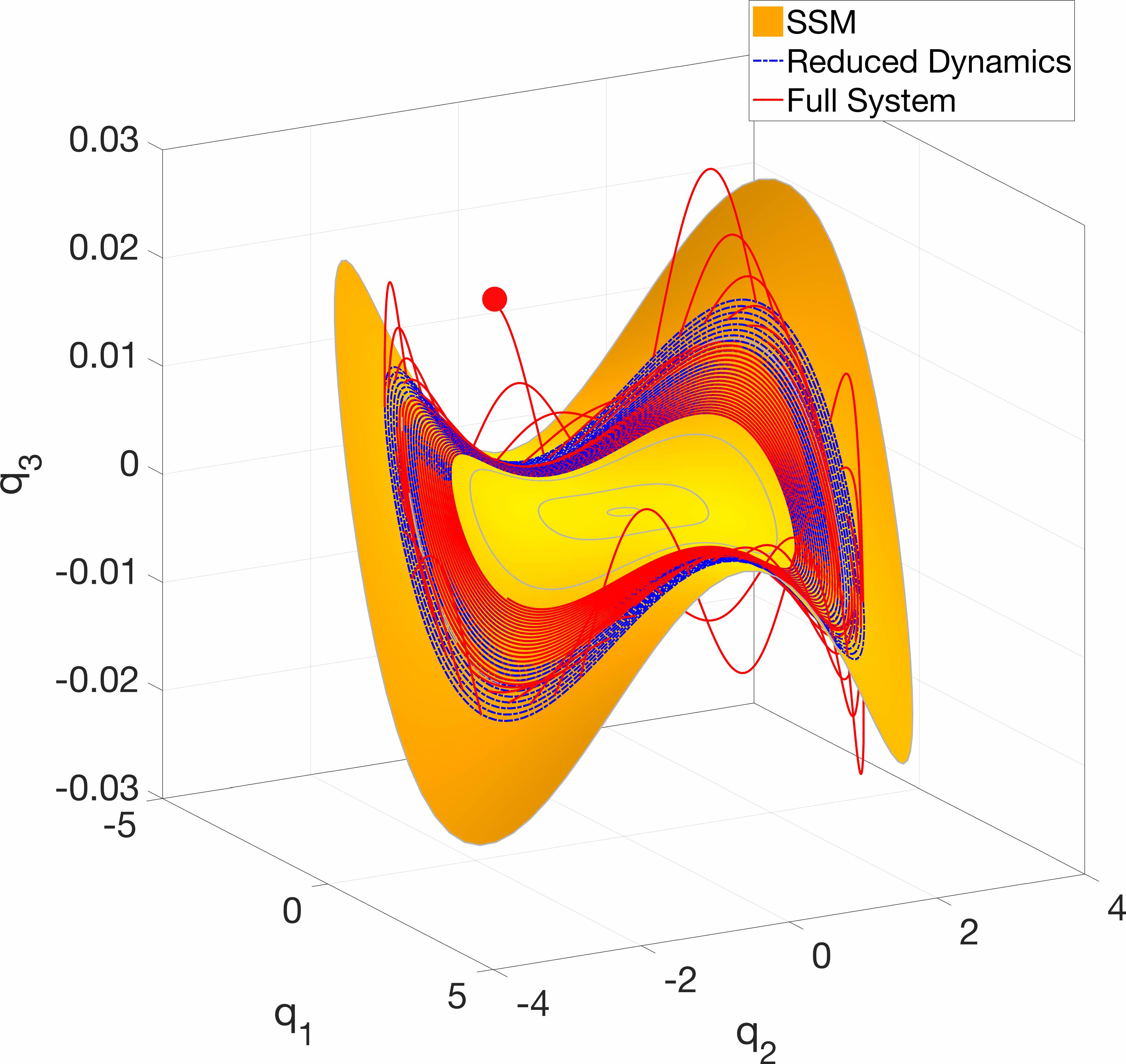}
\end{centering}
}
\vfill{}
\subfloat[\label{fig:beam_SQ}]{
\begin{centering}
\includegraphics[scale=0.4]{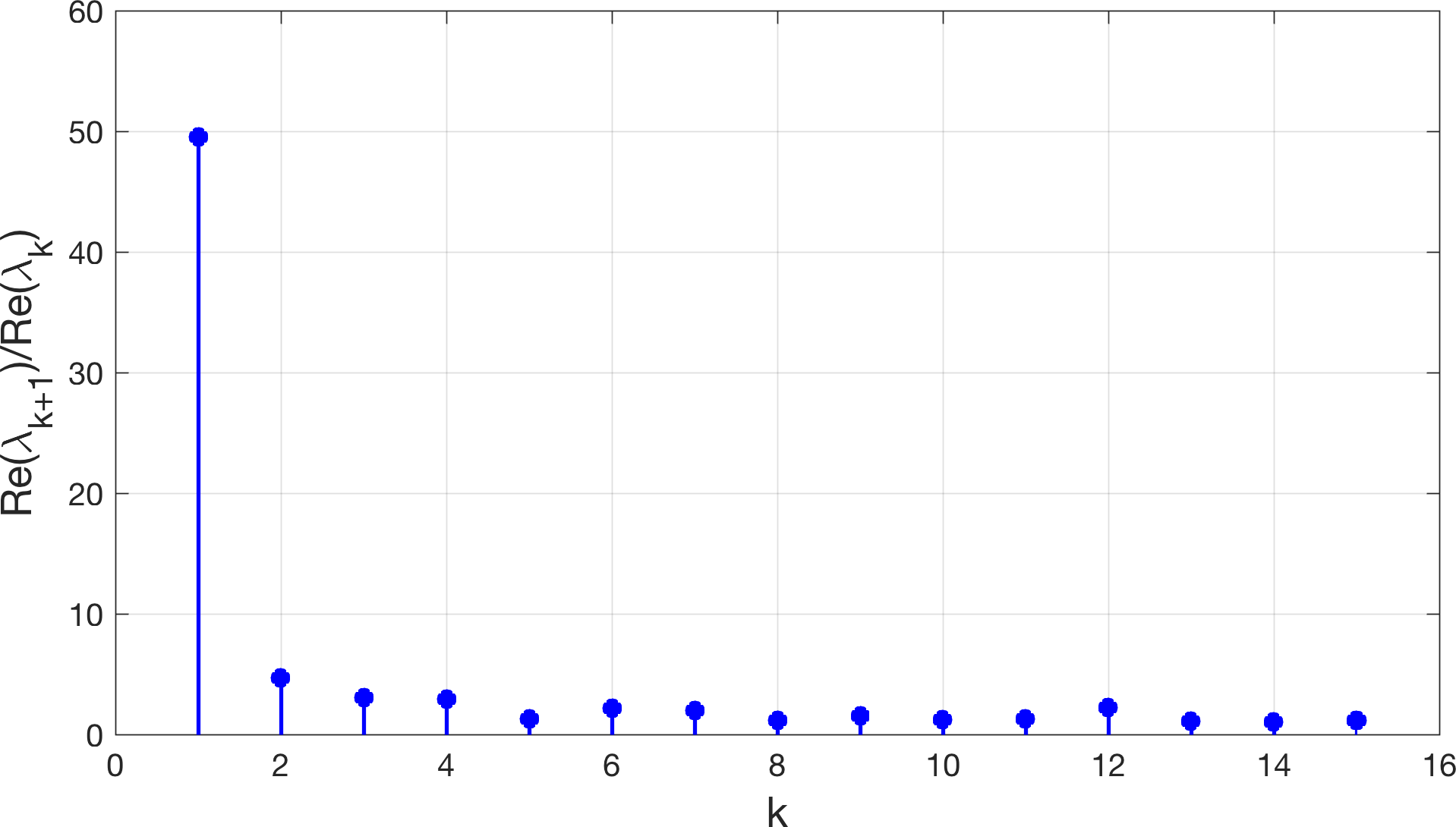}
\end{centering}
}
\end{center}

\caption{Lower-dimensional projections of the 32-dimensional phase space, of
the discretized nonlinear Timoshenko beam, showing the $10^{\text{th}}$
order approximation of $\vec{W}(\mathcal{E})$. Figure \ref{fig:beam_ssm_OT}
shows the spectral submanifold $\vec{W}(\mathcal{E})$ tangent
to slowest subspace $\mathcal{E}$. The dashed (blue) curve corresponds
to a trajectory of the reduced two-dimensional system $\vec{R}(\vec{z})$ corresponding
to the initial positions $r=1.5$ and $\theta=3$. The solid (red)
curve represents a trajectory of the full system having an initial
position off the manifold, showing the convergence towards $\vec{W}(\mathcal{E}_{1})$,
with $t_{\text{end}}=15\text{ s}$. The collapsing nature of the trajectories
onto $\vec{W}(\mathcal{E})$, is a direct consequence of the high
spectral ratio between the slowest eigenspace and the second slowest
eigenspace shown in figure \ref{fig:beam_SQ}. \label{fig:OT_SQ}}
\end{figure}

The continuation algorithm takes 4 hours 42 minutes and 17 seconds
to compute. Additionally, in figure \ref{fig:OT_SQ} we show that
trajectories starting off $\vec{W}(\mathcal{E}_{1})$, will converge towards
$\vec{W}(\mathcal{E}_{1})$ as a consequence of having the high spectral
quotient between the slowest and the remaining eigenspaces. 

\section{Future work}

Currently, SSMtool is capable of computing two-dimensional spectral
submanifolds in autonomous nonlinear mechanical systems of arbitrary
degrees of freedom. As future work, we can extend the current setting
to higher-dimensional spectral submanifolds and, additionally, allow
for possible forcing, as the underlying theory is already worked out
in Haller and Ponsioen \cite{Haller2016}.

\section{Conclusions }

We have developed and tested an automated computational algorithm
for two-dimensional autonomous SSMs that extends modal subspaces of
linear systems to nonlinear systems. Implemented in the MATLAB package
SSMtool, the algorithm can handle non-conservative mechanical systems
of arbitrary (finite) degrees of freedom, subject to numerical memory
limitations only. We used a systematic approach, the parameterization
method, allowing us to construct the SSMs, their reduced dynamics
and corresponding backbone curves up to any required order of precision.

Because the SSMs are embedded using the parameterization method, the
construction of the SSMs does not break down when the SSM folds over
its underlying spectral subspace, as opposed to constructing the SSMs
as graphs over a set of coordinates. The implementation (a MATLAB
based graphical user interface called SSMtool) detects near-outer
and near-inner resonances. In case of an exact outer resonance, the
SSM construction will break down, whereas the presence of near-inner
resonances in general leads to nonlinear terms in the reduced dynamics
on the SSMs. 

Szalai et al. \cite{Szalai2017} have exactly shown how backbone curves
can be extracted from SSMs. We computed the backbone curves, in this
fashion, up to $15^{\text{th}}$-order for a two-degree-of-freedom
non-conservative mechanical system and used the numerical continuation
software COCO \cite{Dankowicz2013} to find periodic orbits of the
periodically forced system for different forcing amplitudes while
varying the forcing frequency to verify the accuracy of the backbone
curves. 

Under an approximate outer resonance, folding of the SSMs is likely
to occur. We demonstrated such folding on the same two-degree-of-freedom
non-conservative mechanical system by varying the system parameters
accordingly to create a near-outer resonance. 

Finally, we have used the developed numerical method to construct
a reduced-order model for a discretized nonlinear Timoshenko beam.
We have computed the reduced dynamics on the two-dimensional slow
SSM arising from the slowest modal subspace of the linearized system.
The backbone curve obtained from the SSMtool shows close agreement
with a single amplitude frequency sweep computed from COCO. While
substantially limited in its scope relative to COCO's, SSMtool has
returned backbone curves in a fraction of the times required by COCO
to construct the response curve for one forcing amplitude. The spectral
quotient between the slowest eigenspace and the second slowest eigenspace
indicates that trajectories transverse to the slow SSM die out fast,
making the slowest SSM an optimal choice for reducing our beam model
to a two-dimensional system of ordinary differential equations.

\section*{Acknowledgments}

We are grateful to Hinke Osinga and Harry Dankowicz for useful remarks
and technical explanations on the numerical continuation toolbox COCO
\cite{Dankowicz2013}. We are also thankful to Shobhit Jain, Thomas
Breunung and Zsolt Verasztó for helpful discussions. 

\appendix

\section{Properties of the Kronecker product\label{sec:Appendix-A}}

In this section, we list several useful properties of the Kronecker
product. For further reading, we refer the reader to A.J. Laub \cite{Laub2005}.
\begin{description}
\item [{(i)}] The Kronecker product is associative, i.e.,

\begin{gather}
(\vec{A}\otimes \vec{B})\otimes \vec{C}=\vec{A}\otimes(\vec{B}\otimes \vec{C}),\label{eq:kp_asso}\\
\vec{A}\in\mathbb{C}^{m\times n},\quad \vec{B}\in\mathbb{C}^{p\times q},\quad \vec{C}\in\mathbb{C}^{r\times s}.\nonumber 
\end{gather}

\item [{(ii)}] The Kronecker product is right-distributive, i.e., 

\begin{gather}
(\vec{A}+\vec{B})\otimes \vec{C}=\vec{A}\otimes \vec{C}+\vec{B}\otimes \vec{C},\label{eq:kp_rdis}\\
\vec{A},\vec{B}\in\mathbb{C}^{m\times n},\quad \vec{C}\in\mathbb{C}^{p\times q}.\nonumber 
\end{gather}

\item [{(iii)}] The Kronecker product is left-distributive, i.e.,

\begin{gather}
\vec{A}\otimes(\vec{B}+\vec{C})=\vec{A}\otimes \vec{B}+\vec{A}\otimes \vec{C},\label{eq:kp_ldis}\\
\vec{A}\in\mathbb{C}^{m\times n},\quad \vec{B},\vec{C}\in\mathbb{C}^{p\times q}.\nonumber 
\end{gather}

\item [{(iv)}] The product of two Kronecker products yields another Kronecker
product:

\begin{gather}
(\vec{A}\otimes \vec{B})(\vec{C}\otimes \vec{D})=\vec{A}\vec{C}\otimes \vec{B}\vec{D}\in\mathbb{C}^{mr\times pt},\label{eq:kp_prod}\\
\vec{A}\in\mathbb{C}^{m\times n},\quad \vec{B}\in\mathbb{C}^{r\times s},\quad \vec{C}\in\mathbb{C}^{n\times p},\quad \vec{D}\in\mathbb{C}^{s\times t}.\nonumber 
\end{gather}

\end{description}

\section{Equations of motion for the nonlinear Timoshenko beam\label{sec:Equations-of-motion-beam}}

We derive the equations of motion for the nonlinear Timoshenko beam,
based on the kinematical and constitutive assumptions made in section
\ref{subsec:application-beam}. Applying Hamilton's principle, we
require that the true evolution of the displacement field between
two specified time instances, $t_{1}$ and $t_{2}$, is a stationary
point of the action functional $S$. Consequently, the variation of
the action functional under a virtual displacement, of our system,
should be identically zero at $t_{1}$ and $t_{2}$, i.e.,

\begin{equation}
\delta S=\intop_{t_{1}}^{t_{2}}-\delta K+\delta U+\delta V\ dt=0,\label{eq:action-functional}
\end{equation}
where $\delta K$ and $\delta U$ are the variations in kinetic and
strain energy due to an arbitrary virtual displacement, and $\delta V$
is the virtual work done by the external forces. For simplicity, we
assume that the virtual work done by internal forces acting on the
beam's cross section is large relative to the work done by internal
forces due to out-of-plane stresses. This corresponds to assuming
a state of either plane stress or plane strain, causing all terms
related to out-of-plane stresses to drop out from the expression for
the internal strain energy. This assumption can be justified by the
fact that the respective neglected quantities would only contribute
to the governing equations nonlinear terms of the fourth order in
the stiffest (and hence typically the fastest decaying) degrees of
freedom, $\phi_{y}$, while all other nonlinearities are of order
three or lower. Thus, those terms can be considered small compared
to the rest. This can also be shown by nondimensionalizing the system
and treating the ratio of the beam's length to its height, $\frac{h}{L}$,
as a small parameter. 

In addition to our viscoelastic constitutive law, we allow for external
damping by introducing a simple damping model similar to the ones
discussed in Lesieutre and Kauffman \cite{Lesieutre2013} and Lesieutre
\cite{Lesieutre2010}. Our model assumes a body force acting proportional
to the displacement velocity:

\begin{equation}
\vec{f}_{c}=c\dot{\vec{u}}(x,y,z).
\end{equation}
The virtual work done by this force is

\begin{align}
\delta V & =\intop_{L}\intop_{A}\ \vec{f}_{c}\delta \vec{u}\ dAdx=\intop_{L}\ \lambda(\dot{u}_{0}\delta u_{0}+\dot{w}\delta w+\frac{I_{2}}{I_{0}}\dot{\phi}_{y}\delta\phi_{y})\ dx,\label{eq:virt_work_damp}
\end{align}
with

\begin{equation}
I_{k}=\intop_{A}z^{k}dA,\quad\lambda=cI_{0}.
\end{equation}
Equation (\ref{eq:virt_work_damp}) can be interpreted as the work
done by two line-distributed forces, $f_{x}=\lambda\dot{u}_{0}$ and
$f_{z}=\lambda\dot{w}$ , proportional to the time derivative of the
transverse displacement and the axial displacement, respectively.
Additionally, we have a line-distributed force couple, $T_{y}=\lambda\frac{I_{2}}{I_{0}}\dot{\phi}_{y}$,
proportional to the time derivative of the rotation of the cross section.
The pre-factor $I_{2}/I_{0}$ in equation (\ref{eq:virt_work_damp})
ensures that the contribution of the external damping to the damping
matrix of the FEM model derived below is proportional to the mass
matrix. As a consequence, the entire system will be subjected to Rayleigh
damping, a damping model which is frequently used in FEM simulations,
see e.g. Takács et al. \cite{Takacs2012}. For more extensive discussions
about how to model damping in beams, we refer the reader to the work
of Lesieutre and Kauffman \cite{Lesieutre2013} and Lesieutre \cite{Lesieutre2010}.
The preceding considerations lead to the following expressions

\begin{align}
\delta K & =\intop_{L}\ m_{0}\dot{u}_{0}\delta\dot{u}_{0}+m_{0}\dot{w}\delta\dot{w}+m_{2}\dot{\phi}_{y}\delta\dot{\phi}_{y}\ dx,\label{eq:elastic_en}\\
\delta U & =\intop_{L}\intop_{A}\ \sigma_{xx}\delta\varepsilon_{xx}+\sigma_{xz}\delta\gamma_{xz}\ dAdx,\label{eq:strain_en}\\
\delta V & =\intop_{L}\ f_{x}\delta u_{0}+f_{z}\delta w+T_{y}\delta\phi_{y}\ dx,\label{eq:ext_en}
\end{align}
with

\begin{equation}
m_{k}=\intop_{A}\ \varrho z^{k}\ dA,
\end{equation}
Substituting (\ref{eq:elastic_en}), (\ref{eq:strain_en}) and (\ref{eq:ext_en})
into (\ref{eq:action-functional}), plugging in the kinematical assumptions
and using integration by parts with respect to $t$ and $x$ we obtain

\begin{align}
\delta S & =\intop_{t_{1}}^{t_{2}}\intop_{0}^{L}\left(m_{0}\ddot{u}_{0}+\lambda\dot{u}_{0}-\partial_{x}\left(M_{xx}^{0}+\phi_{y}M_{xz}^{0}\right)\right)\delta u_{0}+\left(m_{0}\ddot{w}+\lambda\dot{w}-\partial_{x}\left(M_{xz}^{0}+\partial_{x}wM_{xx}^{0}\right)\right)\delta w\nonumber \\
 & +\left(m_{2}\ddot{\phi}_{y}+\lambda\frac{I_{2}}{I_{0}}\dot{\phi}_{y}-\partial_{x}\left(M_{xx}^{1}+\phi_{y}M_{xz}^{1}\right)+\left(M_{xz}^{0}+\partial_{x}u_{0}M_{xz}^{0}+\partial_{x}\phi_{y}M_{xz}^{1}\right)\right)\delta\phi_{y}\ dx\label{eq:functional_S}\\
 & +\left[\left(M_{xx}^{0}+\phi_{y}M_{xz}^{0}\right)\delta u_{0}\right]_{0}^{L}+\left[\left(M_{xz}^{0}+\partial_{x}wM_{xx}^{0}\right)\delta w\right]_{0}^{L}+\left[\left(M_{xx}^{1}+\phi_{y}M_{xz}^{1}\right)\delta\phi_{y}\right]_{0}^{L}\ dt=0,\nonumber 
\end{align}
where we have defined

\begin{equation}
M_{ij}^{k}=\intop_{A}\ \sigma_{ij}z^{k}\ dA.
\end{equation}
Equating the variational derivative of this functional with zero,
we obtain the Euler-Lagrange equations

\begin{align}
 & m_{0}\ddot{u}_{0}+\lambda\dot{u}_{0}-\partial_{x}\left(M_{xx}^{0}+\phi_{y}M_{xz}^{0}\right)=0,\label{eq:eom_1}\\
 & m_{0}\ddot{w}+\lambda\dot{w}-\partial_{x}\left(M_{xz}^{0}+\partial_{x}wM_{xx}^{0}\right)=0,\label{eq:eom_2}\\
 & m_{2}\ddot{\phi}_{y}+\lambda\frac{I_{2}}{I_{0}}\dot{\phi_{y}}-\partial_{x}\left(M_{xx}^{1}+\phi_{y}M_{xz}^{1}\right)+\left(M_{xz}^{0}+\partial_{x}u_{0}M_{xz}^{0}+\partial_{x}\phi_{y}M_{xz}^{1}\right)=0,\label{eq:eom_3}
\end{align}
along with the corresponding boundary conditions

\begin{align}
 & \left.\left(M_{xx}^{0}+\phi_{y}M_{xz}^{0}\right)\right|_{L}\delta u_{0}(L)=0,\quad\left.\left(M_{xx}^{0}+\phi_{y}M_{xz}^{0}\right)\right|_{0}\delta u_{0}(0)=0,\\
 & \left.\left(M_{xz}^{0}+\partial_{x}wM_{xx}^{0}\right)\right|_{L}\delta w(L)=0,\quad\left.\left(M_{xz}^{0}+\partial_{x}wM_{xx}^{0}\right)\right|_{0}\delta w(0)=0,\\
 & \left.\left(M_{xx}^{1}+\phi_{y}M_{xz}^{1}\right)\right|_{L}\delta\phi_{y}(L)=0,\quad\left.\left(M_{xx}^{1}+\phi_{y}M_{xz}^{1}\right)\right|_{0}\delta\phi_{y}(0)=0.
\end{align}
The $M_{ij}^{k}$ terms can be written as a function of our displacement
field by using the kinematical and constitutive relations, i.e.,

\begin{align}
M_{xx}^{0} & =I_{0}\left(E(\partial_{x}u_{0}+\frac{1}{2}(\partial_{x}w)^{2})+\eta(\partial_{x}\dot{u}_{0}+\partial_{x}\dot{w}\partial_{x}w)\right),\\
M_{xx}^{1} & =I_{2}\left(E\partial_{x}\phi_{y}+\eta\partial_{x}\dot{\phi}_{y}\right),\\
M_{xz}^{0} & =I_{0}\left(G(\phi_{y}+\partial_{x}w+\phi_{y}\partial_{x}u_{0})+\mu(\dot{\phi}_{y}+\partial_{x}\dot{w}+\dot{\phi}_{y}\partial_{x}u_{0}+\phi_{y}\partial_{x}\dot{u}_{0})\right),\\
M_{xz}^{1} & =I_{2}\left(G\phi_{y}\partial_{x}\phi_{y}+\mu(\dot{\phi}_{y}\partial_{x}\phi_{y}+\phi_{y}\partial_{x}\dot{\phi}_{y})\right).
\end{align}

We discretize equations (\ref{eq:eom_1})-(\ref{eq:eom_3}) using
a finite-element discretization (cf. Reddy \cite{Reddy2014} for a
more detailed description). We use cubic shape functions to approximate
$u_{0}$, quadratic shape functions for $w$ and linear shape functions
for $\phi_{y}$. A beam element with three equally spaced nodes, situated
at the beginning, the middle and at the end of the element is used.
The node in the middle of the element is only needed for the interpolation
of the transverse displacement $w$. To avoid shear and membrane locking,
the $\varepsilon_{xx}^{0}$ and $\gamma_{xz}^{0}$ terms should be
approximated by shape functions of the same order. After discretization,
we obtain a set of $n$ ordinary differential equations (ODEs) governing
the dynamics the nonlinear Timoshenko beam:

\begin{equation}
\vec{M}\ddot{\vec{y}}+\vec{C}\dot{\vec{y}}+\vec{K}\vec{y}+\vec{f}(\vec{y},\dot{\vec{y}})=\vec{0}
\end{equation}
where we have defined the vector

\begin{equation}
\vec{y}=\left[\begin{array}{c}
\tilde{u}_{0}\\
\tilde{w}\\
\tilde{\phi}_{y}
\end{array}\right]
\end{equation}
representing the discretized degrees of freedom corresponding to the
unknowns $(u_{0},w,\phi_{y})$. The quantities $\vec{M}\in\mathbb{R}^{n\times n}$,
$\vec{C}\in\mathbb{R}^{n\times n}$, $\vec{K}\in\mathbb{R}^{n\times n}$ are the
mass, damping and stiffness matrices of our discretized model, respectively,
and the nonlinear force vector $\vec{f}\in\mathbb{R}^{n}$ is of the form

\begin{equation}
f_{i}=D_{ijk}y_{j}y_{k}+G_{ijk}y_{j}\dot{y}_{k}+H_{ijkl}y_{j}y_{k}y_{l}+L_{ijkl}y_{j}y_{k}\dot{y}_{l}
\end{equation}
where $i\in\left\{ 1,\ldots,n\right\} $, and the Einstein summation
convention is followed. 

\section{Multiple representations for the nonlinear coefficient matrices \label{sec:Multiple-representations}}

The vector $\vec{q}^{\otimes i}$ contains all possible combinations of
its own elements up to order $i$, and therefore will contain all
the monomial terms related to a homogeneous multivariate polynomial
of degree $i$ in the variables $\vec{q}$. The number of unique monomial
terms $S(2n,i)$ in a multivariate polynomial of degree $i$ with
the variables $\vec{q}\in\mathbb{C}^{2n}$ is equal to the number of multisets
of cardinality $i$, with elements taken from the set $\{1,2,\dots,2n\}\in\mathbb{N}^{2n}$
\cite{stanley1997enumerative}, i.e,

\begin{equation}
S(2n,i)=\left(\begin{array}{c}
i+2n-1\\
i
\end{array}\right)=\frac{\left(i+2n-1\right)!}{(2n-1)!i!}.
\end{equation}
To illustrate this, we now give an example.
\begin{example}
\label{ex:Unique-monomial-terms} {[}\textit{Unique monomial terms
representing a multivariate polynomial of degree two} {]} Assume that
$\vec{q}=(q_{1},q_{2},q_{3},q_{4})^{T}\in\mathbb{C}^{4}$, then the unique
monomial terms related to the multivariate polynomial of degree two
in the $\vec{q}$ variables are 

\begin{equation}
\left(\begin{array}{c}
q_{1}^{2}\\
q_{1}q_{2}\\
q_{1}q_{3}\\
q_{1}q_{4}\\
q_{2}^{2}\\
q_{2}q_{3}\\
q_{2}q_{4}\\
q_{3}^{2}\\
q_{3}q_{4}\\
q_{4}
\end{array}\right)\rightarrow\left.\begin{array}{c}
\left\{ 1,1\right\} \\
\left\{ 1,2\right\} \\
\left\{ 1,3\right\} \\
\left\{ 1,4\right\} \\
\left\{ 2,2\right\} \\
\left\{ 2,3\right\} \\
\left\{ 2,4\right\} \\
\left\{ 3,3\right\} \\
\left\{ 3,4\right\} \\
\left\{ 4,4\right\} 
\end{array}\right\} \quad S=\left(\begin{array}{c}
5\\
2
\end{array}\right)=10.
\end{equation}
Here we see the equivalence with the ten multisets of cardinality
two, with elements taken from the set $\left\{ 1,2,3,4\right\} $.
Indeed, the Kronecker product $\vec{q}\otimes \vec{q}$ will result in a $16$-dimensional
vector containing $6$ extra cross terms that are contained in the
ten multisets. 
\end{example}
As a direct consequence of this redundancy for a representation of
$2n$ multivariate polynomials of degree $i$ in the variables $\vec{q}$,
there will be infinitely many possible representations for the nonlinear
coefficient matrix $\vec{G}_{i}$ corresponding to the $i^{\text{th}}$
order in equation (\ref{eq:ds_diag}). We show this in Example \ref{ex:non-unique}.
\begin{example}
\label{ex:non-unique}{[}\textit{Multiple representations for the
matrices $\vec{G}_{i}$}{]} If we assume, for simplicity, that $\vec{q}=(q_{1},q_{2})^{T}\in\mathbb{C}^{2}$,
and $\vec{G}(\vec{q})$ is only of $\mathcal{O}\left(\left|\vec{q}\right|^{2}\right)$,
i.e.\ $\Gamma=2$,

\begin{align}
\vec{G}(\vec{q})=\sum_{i=2}^{2}\vec{G}_{i}\vec{q}^{\otimes i}=\vec{G}_{2}\vec{q}\otimes \vec{q} & =\left[\begin{array}{cccc}
g_{11} & g_{12} & g_{13} & g_{14}\\
g_{21} & g_{22} & g_{23} & g_{24}
\end{array}\right]\left[\begin{array}{c}
q_{1}^{2}\\
q_{1}q_{2}\\
q_{2}q_{1}\\
q_{2}^{2}
\end{array}\right]\label{eq:ex_G2}\\
 & =\left[\begin{array}{c}
g_{11}q_{1}^{2}+(g_{12}+g_{13})q_{1}q_{2}+g_{14}q_{2}^{2}\\
g_{21}q_{1}^{2}+(g_{22}+g_{23})q_{1}q_{2}+g_{24}q_{2}^{2}
\end{array}\right].\nonumber 
\end{align}
Assume that the quadratic nonlinearities of the underlying system
are modeled as follows 

\begin{equation}
\vec{P}(\vec{q})=\left[\begin{array}{c}
a_{1}q_{1}^{2}+b_{1}q_{1}q_{2}+c_{1}q_{2}^{2}\\
a_{2}q_{1}^{2}+b_{2}q_{1}q_{2}+c_{2}q_{2}^{2}
\end{array}\right].\label{eq:ex_P2}
\end{equation}
If we want to transform $\vec{P}(\vec{q})$ into the form of $\vec{G}(\vec{q})$, equating
$\vec{G}(\vec{q})$ and $\vec{P}(\vec{q})$ and collecting terms of equal power in $q_{1}$
and $q_{2}$ gives

\begin{align*}
 & g_{11}=a_{1},\quad(g_{12}+g_{13})=b_{1},\quad g_{14}=c_{1},\\
 & g_{21}=a_{2},\quad(g_{22}+g_{23})=b_{2},\quad g_{24}=c_{2}.
\end{align*}
The presence of the redundant term $q_{2}q_{1}$ in $\vec{q}\otimes \vec{q}$
introduces two extra coefficients $g_{13}$ and $g_{23}$, giving
us the freedom to introduce a constraint between $g_{12}$ and $g_{13}$,
and between $g_{22}$ and $g_{23}$. We then have two independent
equations, each containing two independent variables. For each equation,
an independent constraint can be introduced such that $g_{12}$, $g_{13}$,
$g_{22}$ and $g_{23}$ are uniquely determined and the product $\vec{G}_{2}\vec{q}\otimes\vec{q}=\vec{G}(\vec{q})$
precisely represents $\vec{P}(\vec{q})$. Each monomial term in the vector $\vec{q}^{\otimes i}$ has a unique location in the vector itself, which in turn points to a unique location in the matrix $\vec{G}_{i}$ for each row. In this way, the constraints are automatically satisfied when SSMtool identifies the nonlinearities of the underlying mechanical system.
\end{example}
\section{Memory requirements for the coefficient matrices \label{sec:scaling}}
In our current setting, the most computationally demanding terms in the SSM construction are the summation terms in equation (\ref{eq:Bmat}), which are shown below for the $i^{\text{th}}$ order
\begin{align}
 & \sum_{m=2}^{i-1}\vec{W}_{m}\sum_{\left|\vec{p}\right|=1}\vec{R}_{i+1-m}^{p_{1}}\otimes\ldots\otimes \vec{R}_{i+1-m}^{p_{m}}, \label{eq:sumWR}\\
 & \sum_{m=2}^{i-1}\vec{G}_{m}\sum_{\left|\vec{r}\right|=i}\vec{W}_{r_{1}}\otimes\ldots\otimes \vec{W}_{r_{m}} \label{eq:sumGW}.
\end{align}
If we assume that all matrices in equations (\ref{eq:sumWR}) and (\ref{eq:sumGW}) are densely filled with doubles, where each double has an allocated memory of 8 bytes in MATLAB, the total amount of memory needed (in bytes) to store the matrices in equations (\ref{eq:sumWR}) and (\ref{eq:sumGW}) corresponding to a  mechanical system of n degrees of freedom at the $i^{\text{th}}$ order is equal to  

\begin{equation}
M(n,i)=8\cdot\sum\limits_{m=2}^{i-1}\left((2n)2^m + 2^{m+i}m +(2n)^{m+1} + (2n)^{m}2^{i}c(m,i)\right), \label{eq:mem}
\end{equation}
where $c(m,i)$ is the number of all possible combinations of $m$ positive integers $l_1,\ldots,l_m\in\mathbb{N}^+$, with $\left|\vec{l}\right|=i$.
\begin{example}
\label{ex:memory}{[}\textit{Memory requirements for different orders}{]} We consider a mechanical system of two degrees of freedom ($n=2$) with a single cubic nonlinearity and arbitrary near-inner-resonances. The cubic nonlinear spring will only contribute to the $\vec{G}_3$ coefficient matrix, therefore the only contribution from equation (\ref{eq:sumGW}) to equation (\ref{eq:mem}) is for $m=3$. In figure \ref{fig:memory}, we show the output of equation (\ref{eq:mem}) for different orders of expansion. 

\begin{figure}[H]
\begin{centering}
\includegraphics[scale=0.35]{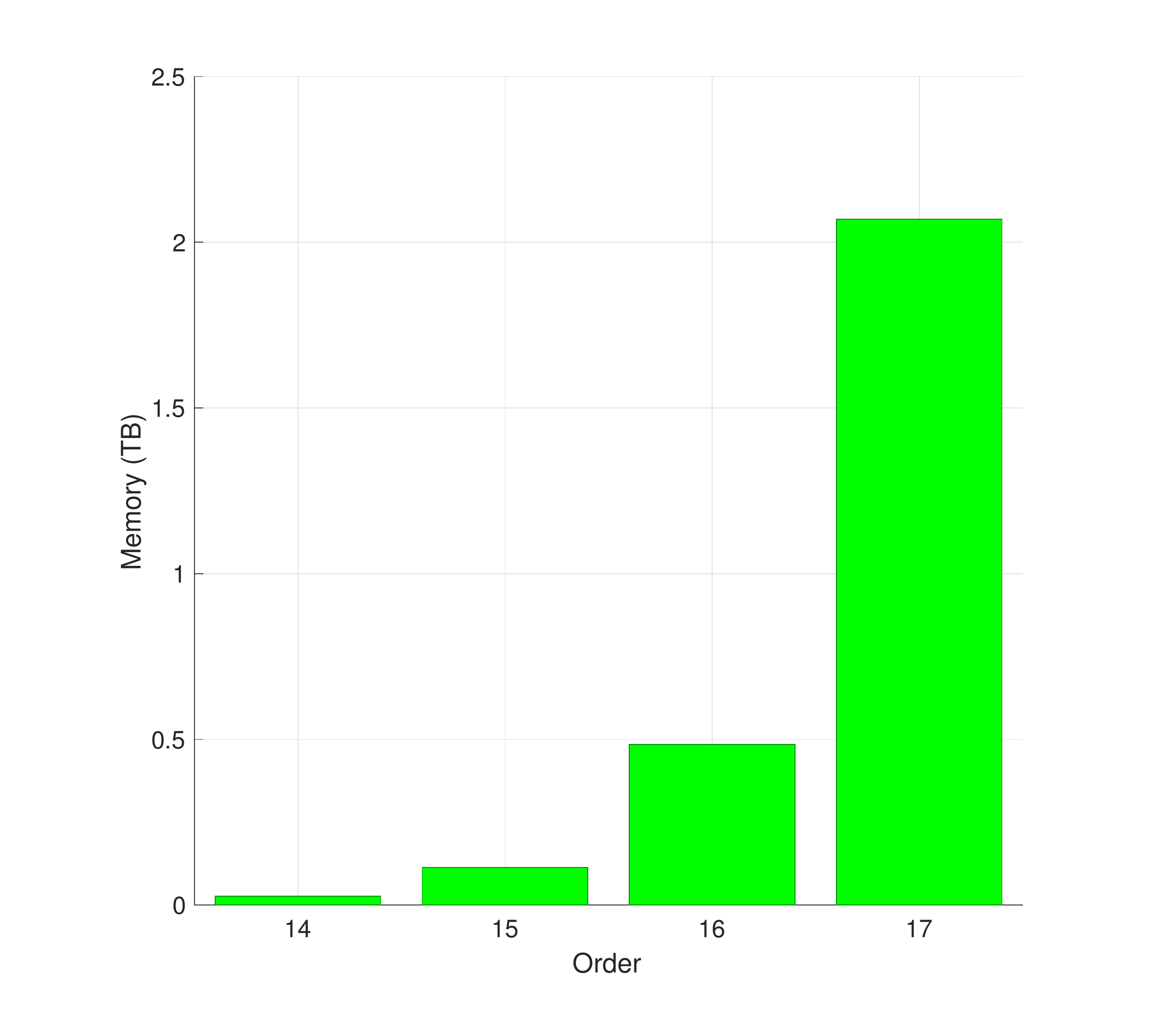}
\par\end{centering}
\caption{Memory requirements in terabytes for equations (\ref{eq:sumWR}) and (\ref{eq:sumGW}) for different orders of the two-degree-of-freedom mechanical system. The amount of memory needed drastically increases from $0.4846 \text{ TB}$    for order $16$ to $2.0696 \text{ TB}$ for order $17$.  \label{fig:memory}}
\end{figure}
\end{example}
\section*{References}
\bibliography{ref}

\end{document}